\documentclass[12pt]{article}
\usepackage{german}
\usepackage{amsmath}
\usepackage{amssymb}
\usepackage{latexsym}

\newcommand{\ang}{\raisebox{0.2ex}{\scriptsize$\triangleright$}}
\newcommand{\mn}{\medskip\noindent}
\newcommand{\bn}{\bigskip\noindent}
\newcommand{\sn}{\smallskip\noindent}

\newcommand{\E}{{\mathcal{E}}}
\newcommand{\F}{{\mathcal{F}}}

\newcommand{\W}{{\mathcal{W}}}
\newcommand{\cB}{{\mathcal{B}}}

\newcommand{\G}{{\mathcal{G}}}
\newcommand{\A}{{\mathcal{A}}}

\newcommand{\J}{{\mathcal{J}}}
\newcommand{\U}{{\mathcal{U}}}
\newcommand{\Z}{{\mathcal{Z}}}
\newcommand{\cP}{{\mathcal{P}}}
\newcommand{\Q}{{\mathcal{Q}}}

\newcommand{\sS}{{\mathcal{S}}}

\newcommand{\Hh}{{\mathcal{H}}}
\newcommand{\cO}{{\mathcal{O}}}
\newcommand{\cD}{{\mathcal{D}}}

\newcommand{\cL}{{\mathcal{L}}}
\newcommand{\cK}{{\mathcal{K}}}

\newcommand{\gr}{{\mathfrak{r}}}
\newcommand{\gA}{{\mathfrak{A}}} 
\newcommand{\gB}{{\mathfrak{B}}}

\newcommand{\kd}{{\mathrm d}}

\newcommand{\Lin}{{\mathrm{Lin}}}

\newcommand{\id}{{\mathrm{id}}}
\newcommand{\rd}{{\mathrm{d}}}

\newcommand{\sx}{{\mathsf x}}
\newcommand{\sy}{{\mathsf y}}

\newcommand{\fF}{{\mathsf F}}
\newcommand{\fK}{{\mathsf K}}
\newcommand{\sK}{{\mathsf K}}
\newcommand{\sE}{{\mathsf E}}
\newcommand{\sF}{{\mathsf F}}

\newcommand{\bbbr}{{\mathbb{R}}}
\newcommand{\bbbc}{{\mathbb{C}}}

\newcommand{\bbbzz}{{\mathbb{Z}}}
\newcommand{\bbbn}{{\mathbb{N}}}

\begin{document}
\begin{center} 
{\LARGE On the Quantum Quarter Plane and the Real Quantum Plane\par}
\end{center}

\centerline{Konrad Schm\"udgen}
\centerline{Fakult\"at f\"ur Mathematik und Informatik}
\centerline{Universit\"at Leipzig, Augustusplatz 10, 04109 Leipzig, Germany}
\centerline{E-mail: schmuedg@mathematik.uni-leipzig.de}

\begin{abstract} 
Suppose $q\ne \pm 1$ is a complex number of modulus one. 
Let $\cO(\bbbr^2_q)$ be the $\ast$-algebra with two hermitean generators
$x$ and $y$ satisfying the relation $xy=qyx$. Using operator 
representations of the $\ast$-algebra $\cO(\bbbr^2_q)$ on Hilbert space 
and the Weyl calculus of pseudodifferential operators we construct 
$\ast$-algebras of ``functions'' on the quantum quarter plane 
$\bbbr^{++}_q$ and on the real quantum plane $\bbbr^2_q$ which are 
left module $\ast$-algebras for the Hopf $\ast$-algebra $\U_q(gl_2(\bbbr))$. We define a family $h_k$, $k \in \bbbzz^2$, of covariant positive linear functionals on these $\ast$-algebras  and study the actions of the $\ast$-algebras $\cO(\bbbr^2_q)$ and $\U_q(gl_2(\bbbr))$ on the associated Hilbert spaces. Quantum analogs of the partial Fourier transforms and the Fourier transform are found. A differential calculus on the ``function'' $\ast$-algebras is also developed and investigated.
\end{abstract}

\mn
Mathematics Subject Classifications (1991), 17B37,81R50,47D40

\mn
{\bf 0. Introduction}

\sn
Suppose that $q\ne \pm 1$ is a complex number of modulus one. 
Let $\cO(\bbbr^2_q)$ be the $\ast$-algebra with two hermitean 
generators $x$ and $y$ satisfying the relation
$$
xy=qyx.
$$
In quantum group theory  $\cO(\bbbr^2_q)$ is usually called the 
coordinate $\ast$-algebra of the {\it real quantum  plane}. It is 
well-known that $\cO(\bbbr^2_q)$ is a left module $\ast$-algebra of 
the Hopf $\ast$-algebra $\U_q(gl_2(\bbbr))$ with action given by formulas 
(\ref{A3})--(\ref{A4}) below. However these structures are not sufficient in order to study analytic properties of the real quantum plane. In the undeformed case $q=1$ the $\ast$-algebra $\cO(\bbbr^2_q)$ is just the polynomial algebra $\bbbc[x,y]$ in two hermitean indeterminates $x$ and $y$ equipped with the usual action of the Lie algebra $gl_2(\bbbr)$. In this situation we can replace the polynomial algebra by the larger $\ast$-algebra $\A(\bbbr^2):=\bbbc[x,y]+C^\infty_0(\bbbr^2)$ of functions on $\bbbr^2$ and extend the action of $gl_2(\bbbr)$ to $\A(\bbbr^2)$. 
On the algebra $C^\infty_0(\bbbr^2)$ we can study differential and integral calculus and so we can develop analysis on $\bbbr^2$. Roughly speaking, for the real quantum plane we try to proceed in a similar way.

An interesting approach to the quantum space $\bbbr^d_q$ has 
been developed by M. Rieffel [R] in the framework of his theory of 
deformation quantization. This approach is essentially based on 
$C^\ast$-algebras.

In the undeformed case the points of $\bbbr^2$ are in one-to-one 
correspondence to the well-behaved irreducible $\ast$-representations 
(see [S1]) of the polynomial algebra $\bbbc[x,y]$ in Hilbert space. Thus 
we are lead to look for well-behaved irreducible $\ast$-representations
of the $\ast$-algebra $\cO(\bbbr^2_q)$. This problem has been studied 
in [S2]. In this paper we consider four irreducible well-behaved $\ast$-representations of $\cO(\bbbr^2_q)$. They are defined as follows. We fix a real number $\gamma$ such that $q=e^{2\pi i\gamma}$ and two reals $\alpha$ and $\beta$ such that $\alpha\beta=\gamma$. Let $\cP$ and $\Q$ be the self-adjoint operators on the Hilbert space $L^2(\bbbr)$ given by 
$(\cP f)(x)=(2\pi i)^{-1} f^\prime(x) $ and $\Q f=xf(x)$. Then, for $\epsilon, \epsilon^\prime \in\{{+},{-}\}$ there exists an irreducible $\ast$-representation $\rho_{\epsilon\epsilon^\prime}$ of the $\ast$-algebra $\cO(\bbbr^2_q)$ on $L^2(\bbbr^2)$ such that 
$$
\rho_{\epsilon\epsilon^\prime}(x)=\epsilon e^{2\pi\alpha\Q},~ \rho_{\epsilon\epsilon^\prime}(y)=\epsilon^\prime e^{2\pi\beta\cP}.
$$
Because of the spectra of operators $\rho_{\epsilon\epsilon^\prime}(x)$ and $\rho_{\epsilon\epsilon^\prime}(y)$, we think of the $\ast$-representation $\rho_{\epsilon\epsilon^\prime}$ as realization of the algebra $\cO(\bbbr^2_q)$ on the {\it quantum quarter plane} $\bbbr^{\epsilon\epsilon^\prime}_q$.

Let us sketch the main ideas of our investigations and begin with the 
$\ast$-representation $\rho_{++}$ corresponding to the (open) quantum quarter  plane $\bbbr^{++}_q$. We want to construct ``functions on the quantum  quarter plane which are vanishing at the boundary''. In order to do so,  we define these ``functions'' as pseudodifferential operators by means  of the Weyl calculus [Fo], [St], that is, 
$$
Op(a)=\gamma\iint\hat{a}(\alpha s,\beta t)e^{2\pi i(s\alpha\Q+t\beta\cP)} ds dt.
$$
As symbol class we take the set of all functions $a(x_1,x_2)$ on 
$\bbbr^2$ which are in the intersection of domains of operators 
$e^{2\pi c_1\Q_1} e^{2\pi d_1\cP_1}\otimes e^{2\pi c_2\Q_2} 
e^{2\pi d_2\cP_2}$, where $c_1,c_2,d_1,d_2\in\bbbr$. This set 
is a $\ast$-algebra, denoted $\gA(\bbbr^2)$, with respect to the 
twisted product and involution of pseudodifferential operators. 
The first aim of our construction is to extend the algebraic structure 
and the left action of $\U_q(gl_2)$ to the direct sum 
$$
\A(\bbbr^{++}_q):=\cO(\bbbr^2_q)+\gA(\bbbr^2)
$$
such that $\A(\bbbr^{++}_q)$ becomes a left $\U_q(gl_2(\bbbr))$-module 
$\ast$-algebra. We think of $\A(\bbbr^{++}_q)$ as counterpart of 
the $\ast$-algebra $\A(\bbbr^{++}):=\bbbc[x,y]+C^\infty_0(\bbbr^{++})$ of 
functions on the quarter plane $\bbbr^{++}$ equipped with the action of 
Lie algebra $gl_2(\bbbr)$. For each $k\in\bbbzz^2$ there exists a 
faithful linear functional $h_k$ on the $\ast$-algebra $\gA(\bbbr^2)$ 
which is covariant with respect to the left action of $\U_q(gl_2(\bbbr))$. 
These functionals $h_k$ can be viewed as quantum anologs of the state 
on the $\ast$-algebra $C^\infty_0(\bbbr^{++})$ given by the Lebesgue 
measure. Further, there are two $U_q(gl_2)$-covariant differential 
calculi on the algebra $\cO(\bbbr^2_q)$ invented in [PW] and [WZ]. 
We extend one of these calculi to a differential calculus on the larger 
algebra $\A(\bbbr^{++}_q)$. Thus the key ingredients
 for a differential and integral calculus on the quantum quarter 
 plane $\bbbr^{++}_q$ are developed.  Similar considerations 
 can be done for three other quantum quarter planes. We strongly believe 
 that this approach, with some technical modifications, could serve as a guide-line 
 for the constructions of other non-compact quantum spaces as well.

The next main step is the construction of the function $\ast$-algebra 
for the real quantum plane $\bbbr^2_q$. The idea is to obtain the quantum 
plane by ``gluing together the four quantum quarter planes at the coordinate 
axis''. In order to do so, it is natural to begin with the direct sum 
$\gA_0(\bbbr^2)_4$ of the four $\ast$-algebras $\gA(\bbbr^2)$ 
corresponding to the quantum quarter planes. The elements 
of $\gA_0(\bbbr^2)_4$ are interpreted as ``functions on the quantum 
plane which are vanishing at the coordinate axis''. As in the case of 
the ``ordinary'' plane we consider $\gA_0(\bbbr^2)_4$ as subspace 
of the Hilbert space obtained from a covariant positive linear functional. 
Then the generators $E^\prime :=q^{1/2} (q-q^{-1})E~~~{\rm and}~~~
F^\prime := q^{-1/2}(q-q^{-1})F$ of $U_q(gl_2[\bbbr))$ act on $\gA_0(\bbbr^2)_4$ 
as symmetric operators which are not essentially self-adjoint. In order 
to remedy this defect, we extend $\gA_0(\bbbr^2)_4$ to a larger 
$\ast$-algebra $\gA(\bbbr^2)_4$ by allowing, roughly speaking, 
symbols having singularities. All hermitean generators $E^\prime, 
F^\prime, q^{-1/4} K_1, q^{-1/4} K_2$ of $\U_q(gl_2(\bbbr))$ and $x,y$ of 
$\cO(\bbbr^2_q)$ act on the larger function $\ast$-algebra $\gA(\bbbr^2)_4$ 
by essentially self-adjoint operators. Because of the singularities of 
symbols, we do not get actions of the {\it whole} $\ast$-algebras 
$\U_q(gl_2(\bbbr))$ and $\cO(\bbbr^2_q)$ on $\gA(\bbbr^2)_4$.

The aims and main steps of our construction are explained, at least to some 
extent, by the preceding discussion. However, the rigorous undertaking of 
this program requires a number of technical lemmas on unbounded operator 
theory, on quantum groups, and on the Weyl calculus. We have 
collected these results in a rather long preliminary Section 1. 
Moreover, notation and terminology are fixed in Section 1. The reader might start 
at Section 2.

Let us describe the content of this paper more in detail. In 
Section 2 we investigate the left $\U_q(gl_2(\bbbr))$-module 
$\ast$-algebra  $\cO (\bbbr^2_q)$ and the covariant first order 
differential calculus on $\cO(\bbbr^2_q)$. In Section 3 the corresponding 
formulas and structures are extended to a larger auxilary 
$\ast$-algebra $\W$. This $\ast$-algebra contains in an operator 
representation on the Hilbert space $L^2(\bbbr^2)$ the operators 
$e^{2\pi i(s\alpha\Q +t\beta\cP)}, s,t\in\bbbr$. In Section 4 the 
left action of $\U_q(gl_2(\bbbr))$ on these operators is used in 
order to derive a left action on operators $Op(a)$ and so on symbols 
$a\in\gA(\bbbr^2)$. In this manner, the symbol algebra $\gA(\bbbr^2)$ 
and the direct sum $\A(\bbbr^{++}_q)=\cO(\bbbr^2_q)+\gA(\bbbr^2)$ of 
vector spaces become left $\U_q(gl_2(\bbbr))$-module $\ast$-algebras. 
This is the function algebra of the quantum quarter plane $\bbbr^{++}_q$ 
mentioned above. We also extend the differential calculus 
of $\cO(\bbbr^2_q)$ to the function algebra $\A(\bbbr^{++}_q)$. In 
Section 5 we define for each $k\in\bbbzz^2$ a $\U_q(gl_2(\bbbr))$-covariant 
faithful positive linear functional $h_k$ on the $\ast$-algebra 
$\gA(\bbbr^2)$ by a weighted integral over the symbol. The $\ast$-representations 
$\psi_k$ of the $\ast$-algebras $\cO(\bbbr^2_q)$ and 
$\U^{tw}_q(gl_2(\bbbr))$ on 
the associated unitary space $\gA_k:=(\gA(\bbbr^2),\langle\cdot,
\cdot\rangle_{k})$ are described on the generators. The 
$\ast$-representations $\psi_k$ and product and involution of 
the $\ast$-algebra $\A(\bbbr^{++}_q)$ are transformed by a 
unitary transformation 
to the Hilbert space $L^2(\bbbr^2)$. These transformed structures are 
essentially used for the construction of the quantum plane in Section 6.
In the last subsection of Section 5 a uniqueness theorem 
for the covariant functional $h_k$ is proved. Section 6 is 
devoted to the construction of the real quantum plane from four 
quantum quarter planes. The function algebra of $\bbbr^2$ can be thought 
as direct sum of function algebras of the four quarter planes with 
boundary  conditions $f(+0,y)=f(-0,y)$ and $f(x,+0)=f(x,-0)$. 
We first give an equivalent formulation of this picture and use then the 
corresponding formulas as motivation for the definitions of structures 
of the real quantum plane. We also define three unitary operators $\F_x^q$, 
$\F_y^q$ and $\F^q$ which interchange up to some powers of the generators 
$K_1$ and $K_2$ the coordinate operators $x$, $y$ and the corresponding 
$q$-deformed partial derivatives $\cD^q_x$, $\cD_y^q$, respectively. 
The unitaries $\F_x^q$, $\F_y^q$, and $\F^q$ can be considered as quantum 
analogs of the partial Fourier transforms and the Fourier transform, respectively.

\bn
\renewcommand{\baselinestretch}{1.0}
{\bf 1. Preliminaries}

\sn
{\it 1.1 Algebraic Preliminaries}
\sn

All algebras in this paper are over the complex field. All the 
notions and facts on Hopf algebras and quantum groups used in this paper 
can be found in [Mo] and [KS], see also [FRT].

Let $\U$ be a Hopf algebra. We use the Sweedler notation $\Delta (f) =
f_{(1)} \otimes f_{(2)}$ for the comultiplication $\Delta (f)$ of $f \in \U$.
A {\it left $\U$-module algebra} $\Z$ is an algebra (without unit in general)
which is a left $\U$-module with left action $\ang$ such that  
\begin{equation}\label{A8}
f\ang (zz^\prime)= (f_{(1)}\ang z)(f_{(2)}\ang z^\prime),~ 
z,z^\prime\in\Z, f\in\  \U.
\end{equation}
A {\it dual pairing} of two Hopf algebras
${\cal U}$ and ${\cal A}$ is a bilinear mapping $\langle \cdot ,
\cdot \rangle : {\cal U}\times {\cal A} \rightarrow \bbbc$ such that
$$
\langle \Delta (f),a_1\otimes a_2\rangle
=\langle f, a_1a_2\rangle , \ \ \
\langle f_1f_2,a\rangle =\langle f_1\otimes f_2,
\Delta (a)\rangle ,  $$
$$
\langle f,1 \rangle =\varepsilon (f),\ \ \
\langle 1,a\rangle =\varepsilon (a),\ \ \
\langle S(f),a\rangle =\langle f,S(a)\rangle  
$$
for all $f,f_1,f_2\in {\cal U}$ and $a,a_1,a_2\in {\cal A}$. 
By a dual pairing 
of two Hopf $\ast$-algebras
${\cal U}$ and ${\cal A}$ we mean a dual pairing $\langle \cdot ,
\cdot \rangle $ of the Hopf algebras
${\cal U}$ and ${\cal A}$ which has the additional property that
\begin{equation}\label{A14}
\langle f^\ast,a\rangle =\overline{\langle f,S(a)^\ast\rangle}~~
\mbox{and}~~ \langle f,a^\ast\rangle = 
\overline{\langle S(f)^\ast,a\rangle},~~
f \in {\cal U},~a \in {\cal A}. 
\end{equation}
Let $\langle .,.\rangle$ be a dual 
pairing of Hopf algebras $\U$ and $\gA$. Any  right 
$\gA$-comodule algebra $\Z$ is a left 
$\U$-module algebra with left action $\ang$ defined by
\begin{equation}\label{A14b}
f \ang z= \langle f,z_{(1)}\rangle z_{(0)},~ f\in\U,z\in\Z,
\end{equation}
where  $\phi(z)= z_{(0)}\otimes z_{(1)}$ is the Sweedler notation 
for the right coaction $\phi$.  

\mn
{\bf Lemma 1.} {\it Suppose that $\langle .,.\rangle$ is a dual pairing 
of Hopf $\ast$-algebras $\U$ and $\gA$. 
If $\Z$ is a right $\gA$-comodule $\ast$-algebra with 
right coaction $\phi: \Z\rightarrow \Z\otimes \gA$, then 
the associated left action of $\U$ on $\Z$ 
satisfies the condition}
\begin{equation}\label{A15}
(f\ang z)^\ast=(S(f)^\ast) \ang z^\ast,~ f\in\U,z\in\Z .
\end{equation}
{\bf Proof.}  Since $\Z$ is a $\gA$-comodule $\ast$-algebra, 
the coaction $\phi$ preserves the involution, so that
$$\phi (z^\ast)\equiv (z^\ast)_{(0)}\otimes 
(z^\ast)_{(1)}=(z_{(0)})^\ast\otimes 
(z_{(1)})^\ast\equiv\phi(z)^\ast.$$ 
Using this condition and the second relation of (\ref{A14})
we conclude that
\begin{equation*}
\begin{split}
S(f)^\ast \ang z^\ast =  \langle S(f)^\ast,(z^\ast)_{(1)}
\rangle (z^\ast)_{(0)} 
=\langle
\overline{f,(z_{(1)})^{\ast\ast}}\rangle (z_{(0)})^\ast =
(f\ang z)^\ast. \hspace{1,5cm}\Box
\end{split}
\end{equation*}
From now on we suppose that $\U$ is a Hopf $\ast$-algebra. 
Equation (\ref{A15}) in Lemma 1 gives the motivation for the following 
definition: A $\ast$-algebra $\Z$ is called 
a {\it left $\U$-module $\ast$-algebra} 
if $\Z$ is a left $\U$-module algebra with 
left action $\ang$ such that equation (\ref{A15}) holds.
Then Lemma 1 says any 
right $\gA$-comodule $\ast$-algebra $\Z$ is a left 
$\U$-module $\ast$-algebra with respect to the associated left action.

Suppose that $\Z$ is a left $\U$-module $\ast$-algebra with left 
action $\ang$ and let $\chi$ be a linear functional on the Hopf 
$\ast$-algebra $\U$. We shall say that a linear functional 
$h$ on $\Z$ is {\it covariant with respect to} $\chi$ if
\begin{equation}\label{G1}
h(f\ang z)=\chi(f)h(z), f\in\U, z\in\Z.
\end{equation}
Suppose for a moment that $h\not\equiv 0$ is covariant with respect to 
$\chi$. Then it follows from the conditions of a left action 
that $\chi$ is a character, that is, 
\begin{equation}\label{G2}
\chi(fg)=\chi(f)\chi(g),~ f,g\in\U,~{\rm and}~ \chi(1)=1.
\end{equation}
If in addition $h$ is hermitian (that is, $h(z^\ast)=\overline{h(z)}$ 
for $z\in\Z$), then we conclude from (\ref{A14}) and (\ref{A15}) that
\begin{equation}\label{G3}
\overline{\chi(f)}=\chi(S(f)^\ast), f\in\U.
\end{equation}
Note that a linear function $h$ on the $\U$-module algebra $\Z$ is invariant if and only if $h$ is covariant with respect to the counit $\varepsilon$.


\mn
{\bf Lemma 2.} {\it Let $h$ be a linear functional on the left 
$\U$-module $\ast$-algebra $\Z$ and set}
\begin{equation}\label{G4}
\langle y,x\rangle:=h(x^\ast y),~ x,y\in\Z.
\end{equation}
{\it Consider the following three conditions:}\\
(i) $h$ {\it is covariant with respect to} $\chi$.\\
(ii) $\chi(f)\langle y,x\rangle=\langle f_{(2)}\ang y, 
S(f_{(1)})^\ast\ang x\rangle$ {\it for} $f\in\U, x, y\in\Z.$\\
(iii) $\langle f\ang y,x\rangle =\chi (f_{(2)})\langle y, 
f_{(1)}^{~~\ast} \ang x\rangle$ {\it for} $ f\in\U, x,y\in\Z.$\\
{\it Then we have} (i)$\rightarrow$(ii)$\rightarrow$(iii). {\it If $\Z$ 
has a unit, then} (iii)$\rightarrow$(i) {\it and so all three 
conditions are equivalent.}

\mn
{\bf Proof.} (i)$\rightarrow$(ii): Using the formulas (\ref{G4}),
(\ref{A15}) and (\ref{A8}) and the fact that 
$S\circ\ast\circ S\circ\ast=\id$ in any Hopf $\ast$-algebra we get
\begin{align*}
\langle f_{(2)}\ang y, S(f_{(1)})^\ast\ang x\rangle &= 
h((S(f_{(1)})^\ast\ang x)^\ast (f_{(2)}\ang y))=
h((S(S(f_{(1)})^\ast)^\ast\ang x^\ast)(f_{(2)}\ang y))\\
&=h((f_{(1)}\ang x^\ast)(f_{(2)}\ang y))=h(f\ang x^\ast y)=
\chi(f)\langle y,x\rangle.
\end{align*}
(ii)$\rightarrow$(iii): Using once more the relation 
$S\circ\ast\circ S\circ\ast =\id$ and condition (ii) we compute
\begin{eqnarray*}
\langle f\ang y, x\rangle &=\langle f_{(2)}\ang y, \varepsilon ((f^{~~\ast}_{(1)})x\rangle=\langle f_{(3)}\ang y, S^{-1}(f_{(2)}^{~~\ast})f_{(1)}^{~~\ast}\ang x\rangle\\
&=\langle f_{(3)}\ang y, S(f_{(2)})^\ast\ang (f_{(1)}^{~~\ast}\ang x)\rangle =\chi(f_{(2)})\langle y, f_{(1)}^{~~\ast}\ang x\rangle.
\end{eqnarray*}
(iii)$\rightarrow$(i): Applying condition (iii) with $x=1$ we obtain
\begin{align*}
h(f\ang y)&=\langle f\ang y, 1\rangle =
\chi (f_{(2)})\langle y, f_{(1)}^{~~\ast}\ang 1\rangle= 
\chi(f_{(2)})\langle y, \varepsilon (f_{(1)}^{~~\ast})1\rangle\\
& =\chi (f_{(2)})\overline{\varepsilon(f_{(1)}^{~~\ast})} h (1^\ast y)=
\chi (f) h(y).\hspace{5cm}\Box
\end{align*}
The special case where $\chi$ is the counit $\varepsilon$ and 
$\Z$ has a unit will be 
stated separately as

\mn
{\bf Corollary 3.} {\it A linear functional $h$ on the left $\U$-module 
$\ast$-algebra $\Z$ with unit is invariant} (that is, 
$h(f\ang z)=\varepsilon (f)h(z)$ for $f\in\U$ and $z\in\Z$) 
{\it if and only if} $\langle f\ang y, x\rangle = \langle y, 
f^\ast\ang x\rangle$ {\it for all} $x,y\in\Z$ {\it and }$f\in\U$.

\mn
Suppose $h$ is an invariant linear functional on the left 
$\U$-module $\ast$-algebra $\Z$ such that the form 
(\ref{G4}) is a scalar product. 
Then, by the implication (i)$\rightarrow$(iii) of Lemma 2, the left action 
of $\U$ on $\Z$ is a $\ast$-representation of the $\ast$-algebra on 
the unitary space $(\Z,\langle\cdot,\cdot\rangle$).

Let $\sigma_1$ and $\sigma_2$ be automorphisms of an algebra $\Z$.
Recall that a linear mapping $\cD$ of $\Z$ is called at 
{\it $(\sigma_1,\sigma_2)$-derivation} if
$$
\cD(z_1z_2) = \sigma_1(z_1)\cD(z_2) + \cD(z_1) \sigma_2(z_2),~~
z_1,z_2 \in \Z.
$$
A {\it first order differential calculus} (briefly, a FODC) over an
algebra $\Z$ is a $\Z$-bimodule $\Gamma$ equipped with a linear mapping
$\kd:\Z\rightarrow\Gamma$ such that $\Gamma$ is the linear 
span of elements $z_1{\cdot}\kd z_2, z_1,z_2 \in \Z$, and $\kd(z_1z_2)=z_1{\cdot}\kd z_2+
 \kd z_1{\cdot} z_2$ for $z_1,z_2 \in \Z$.

\mn

Next we recall the definitions of the Hopf algebras 
$\U_q(gl_2)$, $\U_q(sl_2)$ and $\cO(GL_q(2))$ 
as used in what follows. Let $\U_q(gl_2)$ be the complex 
unital algebra with generators $E$, $F$, $K_1$, $K_2$, $K_1^{-1}$, 
$K_2^{-1}$ and defining relations
$$
K_1K_2=K_2K_1,~ K_jK_j^{-1}=K_j^{-1}K_j=1 ~{\rm for}~ j=1,2,$$
$$
K_1EK_1^{-1}=q^{1/2}E,~K_2EK_2^{-1}=q^{-1/2}E,~
K_1FK_1^{-1}=q^{-1/2}F,~K_2FK_2^{-1}=q^{1/2}F, $$
$$
EF-FE=\lambda^{-1}(K^2-K^{-2}),$$
where we set
$$
K:=K_1K_2^{-1} ~~{\rm and}~~ \lambda := q-q^{-1}.
$$
The algebra $\U_q({\rm gl}_2)$ is a Hopf algebra
with structure maps given on the generators by
$$
\Delta (K_j)=K_j\otimes K_j,~
\Delta (E)=E\otimes K +K^{-1}\otimes E,~
\Delta (F)=F\otimes K +K^{-1} \otimes F,$$
$$
\varepsilon (K_j)=1,~
\varepsilon (E)=\varepsilon (F)=0,~
 S(K_j)=K_j^{-1},~ S(E)=-qE,~S(F)=-q^{-1}F
$$
for $j=1,2.$ Note that the element $L:=K_1K_2$ is group-like and central. 
Let $\U_q(sl_2)$ denote the subalgebra of the algebra $\U(gl_2)$ generated 
by the elements $E, F, K, K^{-1}$. Clearly, 
$\U_q(sl_2)$ is a Hopf subalgebra of $\U(gl_2)$.

Let ${\cal O}(GL_q(2))$ be the coordinate Hopf algebras of the 
quantum group $GL_q(2)$ and let $u_{jk}, j,k=1,2,$ be the entries of 
the corresponding fundamental matrix.
There exists a dual pairing $\langle .,.\rangle$ of the 
Hopf algebras $\U_q(gl_2)$ and $\cO(GL_q(2))$. It is 
determined 
by the values on the generators $K_1, K_2,E,F$ and 
$u_{11},u_{12},u_{21},u_{22},$ respectively, which are given by
\begin{equation}\label{A2}
\langle K_1,u_{11}\rangle=\langle K_2,u_{22}\rangle =q^{-1/2},
\langle K_1,u_{22}\rangle= \langle K_2,u_{11}\rangle =
\langle E,u_{21}\rangle=\langle F,u_{12}\rangle=1
\end{equation}
and zero otherwise.

The algebra with generators $x$ and $y$ satisfying the relation
$xy=qyx$
is called the coordinate algebra ${\cal O}(\bbbc_q^2)$ of 
the {\it quantum plane} $\bbbc_q^2$. It 
is a right $\cO(GL_q(2))$-comodule algebra with 
right coaction $\varphi$ given by 
\begin{equation}\label{A1}
\varphi(x)=x\otimes u_{11}+y\otimes u_{21},~\varphi(y)=x\otimes u_{12}
+y\otimes u_{22}.
\end{equation}
Let $\check{\cO}(\bbbc^2_q)$ denote the algebra with 
geneators $x, x^{-1}, y, y^{-1}$ and relations 
$$
xy =qyx,~ xx^{-1} = x^{-1}x = 1,~ yy^{-1} = yy^{-1} = 1.
$$
That is, $\check{\cO}(\bbbc^2_q)$ is the localization of $\cO(\bbbc^2_q)$ with 
respect to the elements $x$ and $y$ and the algebra $\cO(\bbbc^2_q)$ 
can be considered as a subalgebra of $\check{\cO}(\bbbc^2_q)$.

Assume now  that the deformation parameter $q$ is of 
modulus one. Then there exists an involution $f \rightarrow f^\ast$ on the algebra 
$\U_q(gl_2)$ determined by 
\begin{equation}\label{A13}
E^\ast :=-qE,~ F^\ast := -q^{-1} F,~ K_1^\ast:=K_1,~ K_2^\ast := K_2.
\end{equation}
Equipped with this involution the Hopf algebra $\U_q (gl_2)$ is a Hopf $\ast$-algebra 
denoted by $\U_q(gl_2(\bbbr))$. We often work with the hermitean elements 
$$
E^\prime :=q^{1/2}\lambda E~{\rm and}~F^\prime := q^{-1/2}\lambda F
$$
of the $\ast$-algebra $\U_q(gl_2(\bbbr))$. Further, there is an 
involution $f \rightarrow f^\dagger$
given by 
\begin{equation}\label{A13a}
E^\dagger := -q^{-1}E,~ F^\dagger := -qF,~ 
K_1^\dagger := q^{-1/2} K_1,~K_2^\dagger := q^{-1/2} K_2
\end{equation}
such that $\U_q(gl_2)$ becomes a $\ast$-algebra. It will be denoted by 
$\U^{tw}_q(gl_2(\bbbr))$. In Section 5 we study covariant linear functionals 
with respect to the character $\chi$ on $\U_q(gl_2(\bbbr))$ 
defined by $\chi (K_1)=\chi(K_2)=q^{1/2}$ and $\chi(E)=\chi(F)=0$.
Then, by Lemma 2, the corresponding left action of $\U_q(gl_2)$ is a 
$\ast$-representation of the $\ast$-algebra $\U^{tw}_q(gl_2(\bbbr))$.

The Hopf algebra $\cO(GL_q(2))$ is a Hopf $\ast$-algebra, 
denoted $\cO(GL_q(2,\bbbr))$, with 
involution determined on the generators by 
$u_{jk}^\ast =u_{jk}, j,k=1,2$. The dual pairing $\langle.,.\rangle$ of 
the Hopf algebras $\U_q(gl_2)$ 
and $\cO(GL_q(2))$ given by (\ref{A2}) is also a dual pairing of the 
Hopf $\ast$-algebras $\U_q(gl_2(\bbbr))$ and $\cO(GL_q(2,\bbbr))$.

Further, there exist an involution of the algebra $\cO(\bbbc^2_q))$ given by
\begin{equation} \label{A12}
x^\ast=x, y^\ast=y
\end{equation}
such that this algebra is a $\ast$-algebra. It is denoted by 
$\cO(\bbbr^2_q)$ and called the coordinate $\ast$-algebra of the 
{\it real quantum plane}. From  the preceding formulas we see at once that 
the right coaction $\varphi$ of $\cO(GL_q(2))$ on $\cO(\bbbc^2_q)$ 
preserves the corresponding involutions, that is, $\cO(\bbbr^2_q)$ is a 
right comodule $\ast$-algebra of the Hopf $\ast$-algebra 
$\cO(GL_q(2,\bbbr))$. Hence, by Lemma 1, $\cO(\bbbr^2_q)$ is a left module $\ast$-algebra for the Hopf $\ast$-algebra $\U(gl_2(\bbbr))$. 

\mn
Remark 1. In the literature the involution of $\U_q (gl_2(\bbbr))$ is often defined by the requirements 
$E^\ast =-E, F^\ast=-F, K_1^\ast =K_1, K_2^\ast = K_2$. The latter defines 
indeed an involution which makes $\U_q(gl_2)$ into a 
Hopf $\ast$-algebra. However, with respect to this involution the dual 
pairing with $\cO(GL_q(2,\bbbr))$ is not a dual pairing of Hopf 
$\ast$-algebras and the $\ast$-algebra 
$\cO(\bbbr^2_q)$ is not a left module $\ast$-algebra.

\mn
{\it 1.2. Operator-theoretic Preliminaries}

\sn
First we fix some notation. Let $\J(a,b)$ be the strip 
$\{ z \in \bbbc: a < {\rm Im} z < b\}$. 
The Fourier transform $\F$ and its 
inverse are used in the form
\begin{equation}\label{fourier}
(\F f)(x) = {\hat f}(x) = \int e^{-2\pi itx} f(t) dt,~~
(\F^{-1}f)(x)  = \int e^{2\pi itx} f(t) dt.
\end{equation}
Throughout, we denote the domain of an operator $T$ by $\cD(T)$ and the 
scalar product of $L^2(\bbbr^n)$ by $(\cdot,\cdot)$. Let $\cP$ and $\Q$ be 
the self-adjoint operators on the Hilbert space $L^2(\bbbr)$ 
defined by
$$
(\cP f)(x)=\tfrac{1}{2\pi i}~ f^\prime(x) ~\text{and}~(\Q f)(x)=xf(x).
$$
The operators $\cP$ and $\Q$ are unitarily equivalent by the 
Fourier transform
\begin{equation}\label{fpq}
\F \Q \F^{-1} = - \cP,~~ \F \cP \F^{-1} = \Q.
\end{equation}
The first  assertion of the following lemma describes the 
domain and the action of the operators 
$e^{-2\pi\beta \cP}$ for real $\beta$. We formulate the result for 
$\beta> 0$; the case $\beta < 0$ is completely similar.

\mn
{\bf Lemma 4.} (i): {\it Suppose that} $\beta > 0$. {\it Let $g(z)$ be a
holomorphic function on the strip $\J(0,\beta)$ such that}
\begin{equation*}
\sup_{0<y<\beta}\;\int\limits_{-\infty}^\infty\;\vert g(x{+}iy)\vert^2\,
dx\;<\;\infty\;.
\end{equation*}
{\it Then there exist functions} $g(x)\in L^2(\bbbr)$ {\it and} 
$g_{\beta}(x)\in L^2(\bbbr)$ {\it such that} 
$\lim_{y\downarrow 0}\,g_y=g$ {\it and} 
$\lim_{y\uparrow \beta}\,g_y=g_{\beta}$ in $L^2(\bbbr)$, 
{\it where} $g_y(x):=g(x{+}iy)$ {\it for} $x\in\bbbr$ {\it and} 
$y\in(0,\beta)$. {\it Setting} $g(x{+}i\beta):=g_{\beta}(x), 
x\in\bbbr$, {\it we have}
\begin{equation*}
\lim_{n\rightarrow\infty}\,g(x{+}n^{-2}i )=g(x)~{\rm and}~
\lim_{n\rightarrow\infty}g(x{+}(\beta{-}n^{-2})i)=
g(x{+}\beta i)~{\rm a.e.~on}~ \bbbr \;.
\end{equation*}
{\it The function} $g$ {\it belongs to the domain} 
${\cal D}(e^{-2 \pi \beta \cP})$ {\it and we have}
\begin{equation}\label{C111}
(e^{-2 \pi \beta \cP} g)(x)=
g(x{+}\beta i).
\end{equation} 
{\it Conversely, each function $g$ in the domain} ${\cal D}(e^{-2 \pi\beta \cP})$ {\it arises in this way.}

(ii): {\it For any $\beta\in \bbbr$ and $\delta > 0$, the vector 
space ${\cal D}_\delta :={\rm Lin}\{ e^{-\delta x^2{+} c x}{:}
~ c\in\bbbc\}$ is a core for the self-adjoint operators 
$e^{2 \pi \beta \cP}$ and $e^{2 \pi \beta \Q}$.}

\mn
{\bf Proof.} [S2], Lemma 1--3. \hfill $\Box$

\mn
Throughout this paper we assume that the deformation parameter 
$q\ne  \pm 1$ of 
modulus one and that $\gamma$ is a {\it fixed} real number such that
$$
q=e^{2\pi i\gamma}.
$$
Further, let $\alpha$ and $\beta$ denote real numbers such that 
$\alpha \beta=\gamma$ and put 
\begin{equation}\label{B1}
X=e^{2\pi\alpha \Q}~\text{and}~Y=e^{2\pi\beta \cP}~.
\end{equation}
From (\ref{C111}) it follows that the operators $X$ and $Y$ defined by (\ref{B1}) satisfy the relation $XY\eta=qYX\eta$ for  each vector $\eta \in \cD(XY)\cap\cD(YX)$. Therefore, for each $\epsilon,\epsilon^\prime \in \{+,-\}$, there exist a unique faithful $\ast$-representation $\rho_{\epsilon\epsilon^\prime}$ of the $\ast$-algebra $\cO(\bbbr^2_q)$ on the domain $\gA(\bbbr)$ such that 
\begin{equation}\label{B1a}
\rho_{\epsilon\epsilon^\prime} (x) = \epsilon X{\upharpoonright}  \gA(\bbbr), ~~\rho_{\epsilon\epsilon^\prime} (y) 
= \epsilon^\prime Y{\upharpoonright} \gA(\bbbr).
\end{equation}

Let $\gA(\bbbr)$ be the set of entire holomorphic functions $a(x)$ on the complex plane satisfying

\begin{equation}\label{C2}
{\sup_{\delta_1< y<\delta_2}} {\int\limits^\infty_{-\infty}} |a(x+iy)|^2 
e^{2sx} dx<\infty
\end{equation}

for all $s,\delta_1,\delta_2\in\bbbr$, $\delta_1< \delta_2$. From Lemma 4 we easily derive that 
$$
\gA(\bbbr) = \bigcap^{+\infty}_{n,m=-\infty}\cD(X^n Y^m)
= \bigcap^{+\infty}_{n,m=-\infty}\cD(Y^n X^m).
$$
Clearly, $\gA(\bbbr)$ is invariant under the Fourier transform and its inverse. 

Throughout, we denote by $f_\alpha$ the function 
\begin{equation}\label{f}
f_\alpha(x) = - 2\sinh\pi\beta(2x{+}\alpha i)
\end{equation} 
and by $L_\alpha$ and $R_\alpha$ the operators on the 
Hilbert space $L^2(\bbbr)$ given by
\begin{equation}\label{lrform}
L_\alpha=\overline{f_\alpha} (\cP) e^{-2\pi\alpha Q},~R_\alpha=e^{-2\pi\alpha\Q} f_\alpha(\cP).
\end{equation}
Some properties of these operators are collected in the next lemma.

\mn
{\bf Lemma 5.} (i) {\it $L_\alpha$ is a closed symmetric operator. }\\
(ii) {\it $R_\alpha$ is the adjoint operator of $L_\alpha$.}\\
(iii) {\it $\gA(\bbbr)$ is a core for the operator $L_\alpha$.}\\
(iv) {\it $f_\alpha(\cP)^{-1}\gA(\bbbr)$ is a core for the operator 
$R_\alpha$.}

\mn
{\bf Proof.} By formula (\ref{fpq}), we can 
replace $\cP$ by $\Q$ and 
$\Q$ by $-\cP$. But then the assertions (i)--(iii) have been stated 
in [S2] and [S4].

It remains to prove assertion (iv). First note that 
$f_{\alpha}(\cP)^{-1}$ is a bounded 
normal operator on the Hilbert space $L^2(\bbbr)$, 
so $\cB_{\alpha}:=f_\alpha(\cP)^{-1}\gA(\bbbr)$ 
is a dense linear subspace of $L^2(\bbbr)$. We show that
\begin{equation}\label{F6a}
(R_\alpha{\upharpoonright} \cB_{\alpha})^\ast\subseteq L_\alpha.
\end{equation}
Indeed, suppose that $\zeta\in\cD((R_\alpha{\upharpoonright}\
\cB_{\alpha})^\ast)$.  
Then there exists a vector $\xi\in L^2(\bbbr)$ such that 
$\langle R_\alpha\eta,\zeta\rangle=\langle\eta,\xi\rangle$ 
for all  $\eta\in\cB_{\alpha}$. Writing $\eta$ as 
$\eta=f_\alpha(\cP)^{-1}\eta^\prime$ with 
$\eta^\prime\in\gA(\bbbr)$, we obtain $\langle 
e^{-2\pi\alpha\Q}\eta^\prime,\zeta\rangle = \langle\eta^\prime, 
\overline{f_\alpha}(\cP)^{-1}\xi\rangle$. Since $\gA(\bbbr)$ 
is a core for $e^{-2\pi\alpha\Q}$ by Lemma 4(ii), the latter implies 
that $\zeta\in\cD(e^{-2\pi\alpha\Q})$ and 
$e^{-2\pi\alpha\Q}\zeta=\overline{f_\alpha}(\cP)^{-1}\xi$. 
Thus, we have $\zeta\in\cD(\overline{f_\alpha}(\cP)e^{-2\pi\alpha\Q})=
\cD(L_\alpha)$ which proves (\ref{F6a}).

By the assertion of (i), (\ref{F6a}) implies that $(R_\alpha{\upharpoonright}\cB_{\alpha})^{\ast\ast}
\supseteq (L_\alpha)^\ast= R_\alpha$. But the latter 
means that $\cB_{\alpha}$ is a core for $R_\alpha$.
\hfill$\Box$

\mn
Next we essentially use some results obtained in [S4]. We 
restate them here using a slight different notation. 
For $\delta_1,\delta_2\in\bbbr$, $\delta_1>\delta_2$, let 
$\Hh(\delta_1,\delta_2)$ denote the set of all holomorphic functions 
$f$ on the strip $\J(\delta_1,\delta_2)$ satisfying
$$
\sup_{\delta_1{<}y{<}\delta_2} 
{\int\limits^\infty_{-\infty}} |f(x+iy)|^2 e^{-sx^2} dx<\infty
$$ 
for all $s>0$. By Lemma 2 in [S4], each function 
$f\in\Hh(\delta_1,\delta_2)$ has a.e. boundary limits $f(x+i\delta_1)$ and 
$f(x+i\delta_2)$ on $\bbbr$. For notational simplicity we assume that 
$\alpha>0$. With some obvious modifications all results remain valid for 
$\alpha<0$. 

We apply Theorem 1 in [S4] to the function $f(x):=-2\sinh 2\pi\beta x$ and 
with $\alpha$ replaced by $\alpha /2$. Note that 
$f(x{-}i\alpha/2)=f_{-\alpha} (x)$, where $f_\alpha$ is defined by 
$(\ref{f})$. Then Theorem 1 in [S4] can be restated as follows.

\mn
{\bf Lemma 6.} {\it There exist holomorphic functions 
$w_\alpha\in\Hh(\alpha,0)$ and $v_\alpha\in\Hh(\alpha,-\alpha)$ such that
\begin{align}\label{absolute}
&|w_\alpha(x)|=|v_\alpha (x)|=1\quad a.e. ~{\rm on}~\bbbr,\\
\label{vwr}
&w_\alpha(x)=f_{-\alpha}(x) v_\alpha(x{-}\alpha i), v_\alpha (x)=f_{-\alpha}(x)w_\alpha(x{-}\alpha i)\quad a.e.~{\rm on}~\bbbr.
\end{align}
The functions $w_\alpha,v_\alpha$ are uniquely determined up to a constant factor of modulus one by these properties.}

\sn
Let $W_\alpha$ and $A_\alpha$ denote the operator matrices
\begin{equation}
W_\alpha (\cP)=\Big( \begin{matrix} 
w_\alpha (\cP) & 0\\ 
              0&v_\alpha (\cP)\end{matrix} \Big) , 
A_\alpha =\left(\begin{matrix} 
              0&L_\alpha\\ 
      R_\alpha &0\end{matrix}\right), 
B_\alpha=\left(\begin{matrix}  
             0   &e^{2\pi\alpha Q}\\ 
e^{2\pi\alpha Q} &0\end{matrix}\right).
\end{equation}
Since $|w_\alpha|=|v_\alpha|=1$ a.e. on $\bbbr$ by (\ref{absolute}) and $L^\ast_\alpha=R_\alpha$ by Lemma 5, $W_\alpha(\cP)$ is a unitary operator and $A_\alpha$ and $B_\alpha$ are self-adjoint operators on the Hilbert space $L^2(\bbbr)\otimes L^2(\bbbr)$.

\mn
{\bf Lemma 7.} $W_\alpha (\cP)^\ast A_\alpha W_\alpha(\cP)= B_{-\alpha} ~and~ W_\alpha(\cP)^\ast B_\alpha W_\alpha (\cP)=A_{-\alpha}$. 

\mn
{\bf Proof.}  By (\ref{fpq}), we have $\F W_\alpha(\cP)\F^{-1}= 
W_\alpha(\Q)$, $\F e^{-2\pi\alpha\Q} \F^{-1}=e^{2\pi\alpha \cP}$ and 
$\F L_\alpha\F^{-1}=-2\sin\pi\beta(2x{-}\alpha i) e^{2\pi\alpha\cP}$, 
that is, $\F L_\alpha\F^{-1}$ is the operator $L_f$ with 
$f(x)=-2\sinh 2\pi\beta x$ and $\alpha$ replaced by $\alpha / 2$ 
in the notation of [S4]. Thus, under the Fourier transform the 
assertion $W_\alpha (\cP)^\ast A_\alpha W_\alpha(\cP)=B_{-\alpha}$ is 
just equation (\ref{vwr}) in [S4].

Next we prove that $W_\alpha (\cP)^\ast B_\alpha W_\alpha(\cP)=A_{-\alpha}$. 
Since the self-adjoint operator $A_{-\alpha}$ has no proper self-adjoint 
extension in the same Hilbert space, it suffices to show that 
$W_\alpha(\cP)^\ast B_\alpha W_\alpha (\cP)\supseteq A_{-\alpha}$ which 
means that
\begin{align*}
&w_\alpha(\cP)^\ast e^{2\pi\alpha \Q} v_\alpha (\cP)
\supseteq L_{-\alpha}\equiv \overline{f_{-\alpha}}(\cP) e^{2\pi\alpha \Q},\\
&v_\alpha(\cP)^\ast e^{2\pi\alpha \Q} w_\alpha (\cP)\supseteq R_{-\alpha}\equiv e^{-2\pi\alpha \Q} f_{-\alpha}(\cP) .
\end{align*}
Note that $\overline{f_{-\alpha}(x)}=f_\alpha(x)$. Applying the unitary 
transformation $\F\cdot\F^{-1}$ and using (\ref{fpq}) it follows that the 
latter relations are equivalent to 
\begin{align}\label{wvr}
&\overline{w_\alpha(x)} e^{-2\pi\alpha\cP} v_\alpha(x)\supseteq f_\alpha 
(x)e^{-2\pi\alpha\cP},\\
\label{vwr1}
&v_\alpha(x) e^{-2\pi\alpha\cP} w_\alpha(x)\supseteq e^{-2\pi\alpha\cP} 
f_{-\alpha}(x).
\end{align}
Recall that $f(x\pm {\frac{\alpha}{2} i})=f_{\pm\alpha}(x)$. Therefore, 
formula (\ref{vwr}) can be rewritten as
\begin{align}\label{lf}
&f_{-\alpha} (x)e^{2\pi\alpha\cP}= w_\alpha e^{2\pi\alpha\cP }{\overline{v_\alpha}},\\
\label{rf}
&e^{2\pi\alpha\cP} f_\alpha(x)= v_\alpha e^{2\pi\alpha\cP} \overline{w_\alpha}.
\end{align}
Let $\eta\in\gA (\bbbr)$. Since $\overline{\eta}\in 
\cD(f_{-\alpha}(x)e^{2\pi\alpha\cP})$ and hence 
$\overline{v_\alpha}\overline{\eta} \in\cD(e^{2\pi\alpha\cP})$ by 
(\ref{lf}), we have $v_\alpha\eta\in \cD
(e^{-2\pi\alpha\cP})$. The relations (\ref{wvr}) combined with the 
facts that $w_\alpha\in\Hh(\alpha,0)$ and $v_\alpha\in\Hh(\alpha,-\alpha)$ imply that $v_\alpha(x{+}\alpha i)=f_{-\alpha}(x{+}\alpha i)w_\alpha(x)$ 
a.e. on $\bbbr$ (see e.g. formula (\ref{G3}) in [S4]). Since 
$f_{-\alpha} (x{+}\alpha i)=f_\alpha(x)$, we obtain
$$
\overline{w_\alpha(x)} e^{-\pi\alpha\cP} v_\alpha(x)\eta=\overline{w_\alpha (x)} v_\alpha(x{+}\alpha i)e^{-2\pi\alpha\cP}\eta=f_\alpha(x)e^{-2\pi\alpha\cP}\eta.
$$
That is, the operators $\overline{w_\alpha} e^{-2\pi\alpha\cP} v_\alpha$ and $f_\alpha e^{-2\pi\alpha\cP}$ coincide on the domain $\gA(\bbbr)$. By Lemma 5, $\gA(\bbbr)$ is a core for the closed operator $L_{-\alpha}$ and so for $f_\alpha(x)e^{-2\pi\alpha\cP}{=}$\\ $\F L_{-\alpha} \F^{-1}$. Thus we conclude that 
$$
\overline{w_\alpha} e^{-2\pi\alpha\cP} v_\alpha\supseteq f_\alpha e^{-2\pi\alpha\cP}.
$$

Next we verify the second relation (\ref{vwr}). By (\ref{rf}), $\overline{w_\alpha(x)}f_\alpha(x)^{-1} \overline{\eta(x)}\in\cD(e^{2\pi\alpha\cP})$ and hence $w_\alpha (x)f_{-\alpha}(x)^{-1} \eta(x)\in\cD (e^{-2\pi\alpha\cP})$. From (\ref{wvr}) we derive that $w_\alpha(x{+}\alpha i)=
f_{-\alpha} (x{+}\alpha i)v_\alpha(x)$ (see formula (\ref{G2}) in [S4]). 
Therefore, for $\varphi(x)=f_{-\alpha}(x)^{-1}\eta(x)$ with 
$\eta\in\gA(\bbbr)$ we obtain
\begin{align*}
\overline{v_\alpha (x)} e^{-2\pi\alpha\cP} w_\alpha(x)\varphi&= 
\overline{v_\alpha (x)} e^{-2\pi\alpha\cP} x_\alpha 
(x)f_{-\alpha}(x)^{-1}\eta\\
&=\overline{v_\alpha (x)}w_\alpha(x{+}\alpha i)f_{-\alpha}(x{+}\alpha i)^{-1}e^{-2\pi\alpha\cP}\eta\\
&=\overline{v_\alpha (x)}v_\alpha(x)e^{-2\pi\alpha\cP}\eta=
e^{-2\pi\alpha\cP}f_{-\alpha}(x)\varphi.
\end{align*}
Thus, the operators $\overline{v_\alpha} e^{-2\pi\alpha\cP}w_\alpha$ 
and $e^{-2\pi\alpha\cP} f_{-\alpha}$ coincide on the dense domain 
$f_{-\alpha} (x)^{-1}\gA(\bbbr)$. Since $f_{-\alpha}(\cP)^{-1}\gA(\bbbr)$ 
is a core for $R_{-\alpha}$ by Lemma 5 and 
$\F R_{-\alpha} \F^{-1}=
e^{-2\pi\alpha\cP} f_{-\alpha}(x)$, it follows that 
$\overline{v_\alpha} e^{-2\pi\alpha\cP} w_\alpha\supseteq 
e^{-2\pi\alpha\cP}f_{-\alpha}(x)$. This proves (\ref{vwr}) and 
completes the proof of Lemma 7.\hfill$\Box$

\mn
{\bf Lemma 8.} {\it Let $c,d\in\bbbr$ and $\delta_0>0$. Suppose that 
$8|dc|<1$. Then the vector space $\cL_{\delta_0} = \Lin \{ e_{t,\delta}(x) = e^{2\pi(it x-\delta x^2)} ;t\in\bbbr, 0 < \delta <\delta_0 \}$ is dense 
in $\gA(\bbbr)$ with respect to the norm 
$\parallel{\cdot}\parallel_{c,d}=
\parallel (e^{2\pi c\Q}+e^{-2\pi c\Q})(e^{2\pi d\cP}+e^{-2\pi d\cP})\cdot
\parallel$.}

\mn
{\bf Proof.}  Since both operators 
$e^{2\pi c\Q} + e^{-2\pi c\Q}$ and $e^{2\pi d\cP}+e^{-2\pi d\cP}$ on the 
Hilbert space $L^2(\bbbr)$ are self-adjoint and greater than the identity, 
the operator domain $\E_{c,d}:=\cD((e^{2\pi c\Q}+e^{-2\pi c\Q})(e^{2\pi d\cP}+e^{-2\pi d\cP}))$ equipped with the norm $\|{\cdot}\|_{c,d}$ is a Hilbert space.
Assume to the contrary that the assertion of the lemma is not true. 
Then there exists a non-zero vector $\psi_0\in\E_{c,d}$ such that
\begin{equation}\label{F13}
((e^{2\pi c\Q}{+}e^{-2\pi c\Q})(e^{2\pi d\cP}{+}e^{-2\pi d\cP})\psi_0, 
(e^{2\pi c\Q}{+}e^{-2\pi c\Q})(e^{2\pi d\cP}{+}e^{-2 d\cP})e_{t,\delta})=0
\end{equation}
for $t\in\bbbr$ and $0<\delta<\delta_0$. In order to write this 
relation in another form, we set
\begin{equation}\label{F14}
\psi= (e^{2\pi c\Q}{+}e^{-2\pi c\Q})(e^{2\pi d\cP}{+}e^{-2\pi d\cP})\psi_0,~~
\xi_\delta=e^{-2\pi \delta x^2}(e^{2\pi cx}{+}e^{-2\pi cx})\psi .
\end{equation}
Since $\psi\in L^2(\bbbr)$, it follows that 
$e^{n|x|}\xi_\delta (x)\in L^1(\bbbr)$ for all 
$n\in\bbbn$. Further, we compute
\begin{equation}\label{F15}
\big( (e^{2\pi d\cP}{+}e^{-2\pi d\cP})e_{t,\delta}\big) (x)=
e^{-2\pi \delta(x^2-d^2)}\big( e^{2\pi ix(t+2d\delta){+}\pi t\delta}{+}
e^{2\pi ix(t-2d\delta)-2\pi t\delta}\big).
\end{equation}
Therefore, condition (\ref{F13}) means that the Fourier transform 
$\hat{\xi}_\delta$ of the $L^1$-function $\xi_\delta$ satisfies the relation
\begin{equation}\label{F16}
e^{2\pi t d}\hat{\xi}_\delta(t{+}2\delta d){+}
e^{-2\pi td}\hat{\xi}_\delta(t-2\delta d)=0,~t\in\bbbr~.
\end{equation}
Setting
\begin{equation}\label{F17}
\eta_\delta (t):=e^{\pi(2\delta)^{-1} t^2-
\pi(4\delta d)^{-1} ti}\hat{\xi}_\delta (t)~,
\end{equation}
equation (\ref{F16}) is obviously equivalent to the relation
\begin{equation}\label{F18}
\eta_\delta(t{{+}}2\delta d)=\eta_\delta (t{-}2\delta d), ~t\in\bbbr~.
\end{equation}
Since $e^{n|x|}\xi_\delta (x)\in L^1(\bbbr)$ for $n\in\bbbn$, 
$\hat{\xi}_\delta$ is an entire holomorphic function on the complex 
plane and so is $\eta_\delta$ by (\ref{F17}). Therefore, by (\ref{F18}),
the function $\eta_\delta$ is bounded on 
the real axis and hence $e^{n|t|}\hat{\xi}_\delta (t)\in 
L^1(\bbbr)\cap L^2(\bbbr)$ for $n\in\bbbn$ by (\ref{F17}). Consequently, 
$\xi_\delta (x)$ is an entire holomorphic function and 
$\xi_\delta (x)\in\cD(e^{n|\cP|})$ for $n\in\bbbn$. Hence, by (\ref{F14}),
\begin{equation}\label{F19}
\zeta (x):=(e^{2\pi cx}{+}e^{-2\pi cx}) \psi (x)
=e^{2\pi\delta x^2}\xi_\delta (x)
\end{equation}
is also an entire function.

Computing $\xi_\delta (x)$ by the inverse Fourier transform from 
$\hat{\xi}_\delta (t)$ and using equation (\ref{F16}), we derive that
\begin{equation}\label{F20}
\xi_\delta (x)+e^{-8\pi\delta d^2+8\pi i\delta x}\xi_\delta(x{+}2 di)=0,~ 
x\in\bbbr~.
\end{equation}
Inserting (\ref{F19}) into (\ref{F20}) we conclude that
\begin{equation}\label{F21}
\zeta (x)+\zeta (x{+}2 di)=0~.
\end{equation}
Since $8|dc|<1$ by assumption, the function 
$(e^{2\pi cx}{+}e^{-2\pi cx})^{-1}$ is holomorphic on the strip  
$\J_{d,\varepsilon}:= \{x\in\bbbc :|{\rm Im}~x|<2d{+}\varepsilon\}$ for small $\varepsilon>0$ and 
\begin{equation}\label{F22}
{\rm inf}~ \{ \big|(e^{2\pi cx}{+}e^{-2\pi cx})^{-1}\big|; 
x\in\J_{d,\varepsilon}\}>0~.
\end{equation}
Therefore, since $\xi_\delta (x)\in\cD(e^{-4\pi d\cP})$ as noted above, 
we conclude from Lemma 4 that the function $e^{-2\pi \delta x^2}\psi (x)=
(e^{2\pi cx}{+}e^{-2\pi cx})^{-1}\xi_\delta (x)$ belongs to the domain 
$\cD(e^{-4\pi d\cP})$ and 
\begin{equation}\label{F23}
e^{-4\pi d\cP} ( e^{-2\pi\delta x^2}\psi (x) )=e^{-2\pi\delta(x-2di)^2} 
\psi(x{+}2di)~.
\end{equation}
Note that $\psi (x)\in L^2(\bbbr)$ by construction and 
$\psi(x{+}2di)\in L^2(\bbbr)$ by (\ref{F19}), (\ref{F21}) and (\ref{F22}). 
Since $e^{-2\pi\delta x^2}\psi (x)\rightarrow\psi(x)$ and 
$e^{-4\pi d\cP}(e^{-2\pi\delta x^2}\psi)\rightarrow\psi(x{+}2di)$ 
as $\delta \rightarrow 0$ by (\ref{F23}) and the operator $e^{-4\pi d\cP}$ 
is closed, it follows that $\psi\in\cD(e^{-4\pi d\cP})$ and 
$(e^{-4\pi d\cP}\psi)(x)=\psi(x{+}2di)$.  
Applying this fact and formula (\ref{F21}) we obtain
\begin{equation}\label{F24}
\begin{split}
(e^{-4\pi d\cP}&\psi,\psi) = \int \psi (x{+}2di) \overline{\psi(x)} dx\\
&= -\int ( e^{2\pi cx} {+} e^{-2\pi cx})(e^{2\pi c(x{+}2di)} {+} 
e^{-2\pi c (x{+}2di)})^{-1} |\psi(x)|^2 dx.
\end{split}
\end{equation}
Because of the assumption $8|cd|<1$, we have $\cos 4\pi cd>0$. 
Hence the function under the integral sign in (\ref{F24}) is non-negative. 
On the other hand, since $\psi \neq 0$ by construction, 
we have $(e^{-4\pi d\cP}\psi,\psi)>0$. Thus we arrived at a 
contradiction and the assertion of Lemma 8 is proved.\hfill$\Box$
\sn

Next we turn to some facts on tensor products of certain operators. If 
$T_1$ and $T_2$ are closed operators on a Hilbert space $\Hh$, then the 
symbol $T_1 \otimes T_2$ means the closure of the linear operator on the 
domain $\cD(T_1) \otimes \cD(T_2)$ in the Hilbert space $\Hh \otimes \Hh$ 
defined by $(T_1 \otimes T_2)(\eta_1 \otimes \eta_2)= T_1\eta_1 \otimes 
T_2\eta_2$.

Let $\cP_j$ and $\Q_j,j=1,2$, be the self-adjoint operators on the Hilbert 
space $L^2(\bbbr^2)$ given by
$$
(\cP_jf)(x_1,x_2)={\tfrac{1}{2\pi i}}~{\tfrac{\partial f}
{\partial x_j}}~(x_1,x_2)~{\rm and}~ (\Q_jf)(x_1,x_2)=x_j f(x_1,x_2).
$$
Let $\bbbr^{++}:= \{\mu=(\mu_1,\mu_2)\in\bbbr^2: \mu_1 > 0, \mu_2 > 0 \}$. 
For $\mu=(\mu_1,\mu_2)\in\bbbr^2$ we set
\begin{align*}
&S(e^{2 \pi\mu \Q}) :=(e^{2\pi\mu_1\Q_1}\otimes I+
e^{-2 \pi\mu_1 \Q_1}\otimes I)(I\otimes e^{2 \pi\mu_2 \Q_2}+
I\otimes e^{-2 \pi\mu_2 \Q_2}),\\
&e^{2 \pi\mu \Q}:= e^{2 \pi\mu_1 \Q_1}\otimes e^{2 \pi \mu_2 \Q_2}.
\end{align*}
The operators $S(e^{2\pi\mu\cP}), e^{2\pi\mu\cP}$, $e^{2 \pi\mu|\cP|}$ and 
$e^{2 \pi\mu|\Q|}$ are defined in 
a similar manner. Then we have
\begin{equation}\label{ds}
\cD_{\mu,\nu} := \cD(S(e^{2\pi\mu \Q})S(e^{2\pi\nu\cP}))=
{\bigcap\limits_{\varepsilon,\delta\in\bbbzz^2_2}} 
\cD(e^{2\pi\varepsilon\mu \Q} e^{2\pi\delta\nu \cP})
= {\bigcap\limits_{\varepsilon,\delta\in\bbbzz^2_2}} 
\cD(e^{2\pi\varepsilon\nu \cP} e^{2\pi\delta\mu \Q}),
\end{equation}
where $\bbbzz^2_2=\{\varepsilon=(\varepsilon_1,\varepsilon_2):
\varepsilon_1,\varepsilon_2\in\{1,-1\}\}$ and 
$2\pi\varepsilon\mu:=(2\pi\varepsilon_1\mu_1,2\pi\varepsilon_2\mu_2)$. 
If $\nu \in \bbbr^{++}$, then 
$\cD_{\mu,\nu}$ is the vector space of all holomorphic 
functions on $\{(z_1,z_2) \in \bbbc^2:|{\rm Im}~z_j |< |\nu_j|, j = 1,2 \}$
satisfying 
\begin{equation}
{\sup_{|y_j| < |\nu_j|}} 
\iint |a(x_1{+}iy_1,x_2{+}iy_2)|^2 e^{4\pi (|\mu_1 x_1|+|\mu_2 x_2|)} 
dx_1dx_2 <\infty~.
\end{equation}
The latter fact can be proved in a similar manner as Lemma 1.1 in 
[S2] using the Paley-Wiener Theorem.

\mn
{\bf Lemma 9.} (i) {\it Suppose that $\mu=(\mu_1,\mu_2)\in\bbbr^2$ and
$\nu=(\nu_1,\nu_2)\in\bbbr^{++}$.\\
If $f\in\cD(S(e^{2\pi \nu \cP}) e^{2\pi\mu \Q}) 
\cap \cD(e^{2\pi\mu \Q} S(e^{2\pi\nu\cP}))$, then}
\begin{equation}\label{F2}
|f(x_1{+}iy_1, x_2{+}iy_2)|\le {\frac{1}{2\pi}} ((\nu_1{-}|y_1|)
(\nu_2{-}|y_2|))^{-1/2} e^{-2\pi(\mu_1 x_1{+}\mu_2 x_2)}
\parallel e^{2\pi\nu |\cP |}e^{2\pi\mu \Q} f\parallel
\end{equation}
{\it for} $x_1,x_2,y_1,y_2\in\bbbr, |y_1|< \nu_1, |y_2|< \nu_2.$

\sn 
(ii) {\it Let $\mu=(\mu_1,\mu_2),\nu=(\nu_1,\nu_2)\in\bbbr_{++}$. If
$f\in\cD_{\mu,\nu}$, then}
\begin{align}\label{F4a}
&|f(x_1{+}iy_1, x_2{+}iy_2)|\le\nonumber\\
&\quad {\frac{1}{2\pi}}((\nu_1{-}|y_1|)(\nu_2{-}|y_2|))^{-1/2} e^{-2\pi(\mu_1|x_1|{+}\mu_2|x_2|)}{\sum\limits_{\varepsilon,\delta\in\bbbzz^2_2}}\parallel e^{2\pi\varepsilon\nu\cP} e^{2\pi\delta\mu \Q} f
\parallel
\end{align}
{\it for} $x_1,x_2,y_1,y_2\in\bbbr, |y_1|< \nu_1, |y_2|< \nu_2$. {\it The vector space $\cD_{\mu,\nu}$ is contained in the Schwartz space $\sS(\bbbr^2)$.}

\mn
{\bf Proof.} (i): Setting $g=\F f$ and $\varepsilon_j=\nu_j{-}|y_j|,j=1,2$, 
and using formulas (\ref{fpq}), we get
\begin{align*}
&\left| e^{2\pi(\mu_1x_1{+}\mu_2x_2)} f(x_1{+}iy_1,x_2{+}iy_2)\right|=
\left| (e^{-2\pi y\cP} e^{2\pi\mu \Q} \F^{-1}g) (x_1,x_2)\right|\\
&\quad =\left| (\F^{-1} e^{-2\pi y \Q} e^{-2\pi\mu\cP} g)(x_1,x_2)\right|\\
&\quad =\left|\iint e^{2\pi i(x_1t_1+x_2t_2)} (e^{-2\pi y\Q} e^{-2\pi \mu\cP} g)(t_1,t_2) dt_1 dt_2\right|\\
&\quad \le \left(\iint e^{-4\pi(\varepsilon_1 |t_1|+\varepsilon_2|t_2|)} dt_1 dt_2\right)^{1/2}\parallel e^{2\pi\varepsilon |\Q|} e^{-2\pi y \Q} e^{-2\pi\mu\cP} g\parallel_{L^2(\bbbr^2)}\\
&\quad =((2\pi\varepsilon_1)(2\pi\varepsilon_2))^{-1/2}\parallel e^{2\pi\varepsilon |\Q| -2\pi y \Q} e^{-2\pi\mu\cP} g\parallel\\
&\quad = {\frac{1}{2\pi}} (\varepsilon_1\varepsilon_2)^{-1/2}\parallel e^{2\pi((\nu_1-|y_1|) |t_1|-y_1t_1+(\nu_2-|y_2|)|t_2|-y_2t_2)} e^{2\pi\mu\cP} g\parallel\\
&\quad\le {\frac{1}{2\pi}} (\varepsilon_1\varepsilon_2)^{-1/2}\parallel 
e^{2\pi (\nu_1|t_1|+\nu_2|t_2|)} e^{-2\pi\mu\cP} \F f\parallel\\
&\quad ={\frac{1}{2\pi}} ((\nu_1{-}|y_1|)(\nu_2{-}|y_2|))^{-1/2}\parallel e^{2\pi\nu|\cP|} e^{2\pi\mu \Q} f\parallel,
\end{align*}
which proves (\ref{F2}). Note that by the domain assumptions on $f$ the 
function $g=\F f $ belongs to the corresponding operator domains.

\sn
(ii): Since obviously $\parallel e^{2\pi\nu |\cP|} 
e^{2\pi\mu \Q}f \parallel\le \sum_{\varepsilon,\delta}\parallel 
e^{2\pi\varepsilon\nu\cP} e^{2\pi\delta\mu \Q} f\parallel$, 
inequality (\ref{F4a}) follows at once from (\ref{F2}) applied 
with $\mu$ replaced by $\varepsilon\mu$.

Finally, we prove that  $\sS(\bbbr^2) \supseteq \cD_{\mu,\nu}$. 
Let $a \in \cD_{\mu,\nu}$. By (\ref{ds}), we have 
$a \in \cD(e^{\mu |\Q|})$ which implies that $a \in \cD(\Q_1^n \otimes \Q_2^m)$ 
for all $n,m \in \bbbn_0$. Similarly, since $\F(a) \in \cD_{\nu,\mu}$ by 
(\ref{fpq}), we have $\F(a) \in \cD(\Q_1^n \otimes \Q_2^m)$ and hence 
$a \in \cD(\cP_1^n \otimes \cP_2^m)$ for $n,m \in \bbbn_0$. Both conditions implies that $a$ belongs to the Schwartz space $\sS(\bbbr^2)$ 
(see, for instance, Example 10.2.14 in [S1] for this 
apparantly weaker characterization of the Schwartz space). \hfill$\Box$

\mn
{\bf Lemma 10.} {\it Let $c=(c_1,c_2), d=(d_1,d_2) \in \bbbr^2$, 
$\delta_1 {>} 0$, $\delta_2 {>} 0$. Suppose that $8|c_jd_j| {<} 1$ for $j=1,2$. 
Then the vector space $\cL_{\delta_1} \otimes \cL_{\delta_2}$ is dense 
in $\gA(\bbbr^2)$ with respect to the norm}
\begin{equation}\label{F2a} 
\parallel \cdot \parallel_{c,d} 
:=\parallel S(e^{2\pi c\Q})S(e^{2\pi d\cP}) \cdot \parallel .
\end{equation}

\mn
{\bf Proof.} Assume the contrary. Then there exists a vector $\psi \neq 0$ which 
is orthogonal in $L^2(\bbbr^2)$ to 
$S(e^{2\pi c\Q})S(e^{2\pi d\cP})(\cL_{\delta_1} \otimes \cL_{\delta_2})$. 
From the assertion of Lemma 8 it follows that 
$(e^{2\pi c_j\Q_j}+e^{-2\pi c_j\Q_j})(e^{2\pi d_j\cP_j}+
e^{-2\pi d_j\cP_j}) \cL_{\delta_j}$ is dense in $L^2(\bbbr)$ for $j=1,2$. 
But this in turn implies that $\psi =0$. \hfill$\Box$ 

\sn
Let $\gA(\bbbr^2)$ be the intersection of all 
domains $\cD(S(e^{c \Q})S(e^{d\cP}))$, $c,d \in \bbbr^2$, or 
equivalently the vector space of all holomorphic functions on $\bbbc^2$ 
satisfying condition (\ref{F2}) for all $\mu, \nu \in \bbbr^2$. Let $\tau$ denote the locally convex topology on $\gA(\bbbr^2)$ defined by the family of 
norms (\ref{F2a}), $c,d, \in \bbbr^2$. Since it obviously suffices to take 
a countable 
subfamily of such norms, the topology $\tau$ is metrizable. 
Since $\gA(\bbbr^2)$ is the intersection of domains 
$\cD(e^{2 \pi c \Q} e^{2 \pi d \cP})$, $\gA(\bbbr^2)$ is complete with 
respect to this topology. Thus $\gA(\bbbr^2)[\tau]$ is a Frechet space.
The space $\gA(\bbbr^2)$ will play a crucial role as symbol algebra 
for the Weyl calculus.

\mn
{\it 1.3. The Weyl Calculus}

\sn
In this subsection we shall be concerned with pseudodifferential 
operators on the Hilbert space $L^2(\bbbr)$ defined by means of the 
Weyl calculus. Our standard references in this matter are the books 
[Fo] and [St], see also [GV] and [H].
The Weyl correspondence assigns an operator $Op(a)$
to any function $a$ on $\bbbr^2$ such that $\hat{a} \in L^1(\bbbr^2)$ by 
\begin{equation}\label{C3}
Op(a)=\gamma\iint \hat{a}(\alpha s,\beta t) 
e^{2\pi i(s\alpha \Q+t\beta\cP)}~ ds dt.
\end{equation}
Recall that $\hat{a}$ is the Fourier transform (\ref{fourier}) of the 
function $a$. $\alpha$ and $\beta$ are real numbers such that $\alpha\beta=\gamma$ and $q=e^{2\pi i\gamma}$. Since $\hat{a} \in L^1(\bbbr^2)$, the integral (\ref{C3}) can be understood as a Bochner integral and it defines a bounded operator $Op(a)$ on the Hilbert space $L^2(\bbbr)$.

Let us restate some well-known facts on the Weyl calculus 
(see [Fo], Chapter 2). The operator $Op(a)$ acts by the formula
\begin{equation}\label{C4}
(Op(a)f)(x)={\iint} a({\tfrac{1}{2}} (x{+}y),t) e^{2\pi 
i(x{-}y)t} f(y) dy dt.
\end{equation}
For the operator product $Op(a)Op(b)$ and the adjoint operator 
$Op(a)^\ast$ we have
\begin{equation}\label{oppro}
Op(a)Op(b) = Op(a{\scriptstyle{\#}} b)~~{\rm and}~~ Op(a)^\ast=Op(a^\ast),
\end{equation}
where the symbols $a {\scriptstyle{\#}} b$ and $a^\ast$ are defined by 

\begin{align}
&(a{\scriptstyle{\#}} b)(x_1,x_2): =\nonumber\\
\label{C5}
& 4 {\iiiint} a(u_1,u_2)b(v_1,v_2) 
e^{4\pi i[(x_1-u_1)(x_2-v_2)-(x_1-v_1)(x_2-u_2)]} du_1 du_2 dv_1 dv_2,\\
\label{C5a}
&a^\ast(x_1,x_2):=\overline{a(x_1,x_2)},~ x_1,x_2\in\bbbr.
\end{align}

\mn
{\bf Lemma 11.} {\it Let }$\mu = (\mu_1,\mu_2)$, $\nu = (\nu_1,\nu_2)$, 
$\mu^\prime = (\mu^\prime_1, \mu^\prime_2)$, 
$\nu^\prime = (\nu^\prime_1, \nu^\prime_2) \in \bbbr^{++}$, $a\in \cD_{\mu,\nu}, b\in\cD_{\nu^\prime,\mu^\prime}$. {\it For $t\in\bbbc$, we have}

\begin{alignat}{2}\label{F5}
&e^{2\pi t\Q_1}(a{\scriptstyle{\#}} b)=(e^{2\pi t\Q_1}a){\scriptstyle{\#}}(e^{\pi t\cP_2}b)~~
{\rm if}~ |{\rm Re}~t|<\mu_1, |{\rm Re}~t|<2\nu^\prime_2,\\
\label{F6}
&e^{2\pi t\Q_1}(a{\scriptstyle{\#}} b)=(e^{-\pi t\cP_2}a){\scriptstyle{\#}}(e^{2\pi t\Q_1}b)~~
{\rm if}~ |{\rm Re}~t|<2\nu_2, |{\rm Re}~t|<\mu^\prime_1,\\
\label{F7}
&e^{2\pi t\Q_2}(a{\scriptstyle{\#}} b)=(e^{2\pi t\Q_2}a){\scriptstyle{\#}}(e^{-\pi t\cP_1}b)~~
{\rm if}~ |{\rm Re}~t|<\mu_1, |{\rm Re}~t|<2\nu^\prime_1,\\
\label{F8}
&e^{2\pi t\Q_2}(a{\scriptstyle{\#}} b)=(e^{\pi t\cP_1}a){\scriptstyle{\#}}(e^{2\pi t\Q_2}b)~~
{\rm if}~ |{\rm Re}~t|<2\nu_1, |{\rm Re}~t|<\mu_2^\prime,\\
\label{F9}
&e^{2\pi t\cP_1} (a{\scriptstyle{\#}} b)=
(e^{2\pi\cP_1}a){\scriptstyle{\#}} (e^{2\pi t\cP_1} b)~
{\rm if}~|{\rm Re}~t|<\nu_1,|{\rm Re}~t|<\nu_1^\prime.\\
\label{F10}
&e^{2\pi t\cP_1}(a{\scriptstyle{\#}} b)=
(e^{4\pi t\Q_2}a){\scriptstyle{\#}} (e^{-4\pi t\Q_2} b)~~
{\rm if}~|{\rm Re}~t|< \nu_2 /2,|{\rm Re}~t|< \nu_2^\prime /2.\\
\label{F11}
&e^{2\pi t\cP_2} (a{\scriptstyle{\#}} b)=
(e^{2\pi t\cP_2}a){\scriptstyle{\#}} (e^{2\pi t\cP_2} b)~~
{\rm if}~|{\rm Re}~t|<\nu_2,|{\rm Re}~t|<\nu_2^\prime.\\
\label{F12}
&e^{2\pi t\cP_2} (a{\scriptstyle{\#}} b)=
(e^{-4\pi t\Q_1}a){\scriptstyle{\#}} (e^{4\pi t\Q_1} b)~~
{\rm if}~|{\rm Re}~t|< \nu_1/2,|{\rm Re}~t|< \nu_1^\prime/2.
\end{alignat}
{\bf Proof.} As samples, we carry out the proofs of formulas (\ref{F5}) 
and (\ref{F10}). The other equations are proved by a similar reasoning.

First we prove formula (\ref{F10}) for real $t$. It is well-known 
(see [Fo], p. 104) that the Fourier transform of the product $a{\scriptstyle{\#}} b$ 
is the twisted convolution of the Fourier transform $\F(a)$ and $\F(b)$, 
that is\, $\F(a{\scriptstyle{\#}} b)=\F(a)\ast_t \F(b)$, where

$$
(c\ast_t d)(x_1,x_2)=\iint c(u_1, u_2) d(x_1{-}u_1, x_2{-}u_2)
e^{\pi i(x_1 u_2{-}x_2 u_1)} du_1 du_2.
$$
Using the preceding fact and formula (\ref{fpq}) we compute
\begin{align*}
&\F(e^{2\pi t\cP_1} (a{\scriptstyle{\#}} b))(x_1,x_2)=\left( e^{2\pi t\Q_1} 
\F (a{\scriptstyle{\#}} b)\right) (x_1,x_2)\\
&= e^{2\pi t x_1}(\F(a)\ast_t \F(b))(x_1,x_2)\\
&=\iint \F(a)(u_1,u_2) \F(b)(x_1{-}u_1, x_2{-}u_2) 
e^{\pi i(x_1(u_2{-}2t i)-x_2u_1)}du_1 du_2\\
&=\iint \F(a)(u_1, u_2{+}2 ti)\F(b)(x_1{-}u_1,x_2{-}u_2{-}2ti) 
e^{\pi i(x_1u_2-x_2u_1)} du_1 du_2\\
&=\iint \left( e^{-4\pi t\cP_2} \F(a)\right)(u_1,u_2)
(e^{4\pi t\cP_2} \F(b))(x_1{-}u_1, x_2{-}u_2) 
e^{\pi i(x_1u_2-x_2u_1)} du_1 du_2\\
&=\left( e^{-4\pi t\cP_2} \F(a)\ast_t e^{4\pi t\cP_2} \F(b)\right)(x_1,x_2)\\
&=(\F\left( e^{4\pi t\Q_2}a\right)\ast_t \F(e^{-4\pi t\Q_2} b))(x_1,x_2)\\
&=\F\left(e^{4\pi t Q_2} a{\scriptstyle{\#}} e^{-4\pi t\Q_2} b\right)(x_1,x_2)
\end{align*}

which in turn implies (\ref{F10}). It remains to justify the fourth 
equality sign which follows by the formal substitution 
$u_2\rightarrow u_2{+}2ti$. First we note that the assumptions 
$a\in\cD_{\mu,\nu}$ and $b\in\cD_{\mu^\prime,\nu^\prime}$ imply 
that  $\F(a)\in\cD_{\nu,\mu}$ and $\F(b)\in\cD_{\nu^\prime,\mu^\prime}$, 
so that $\F(a)\in\cD(e^{\pm 2\pi\nu_2\cP_2})$ and 
$\F(b)\in\cD(e^{\pm 2\pi\nu^\prime_2\cP_2})$. Therefore, since $2|t|<\nu_2$ 
and $2|t|<\nu^\prime_2$, the function
$$
\F(a)(u_1,u_2) \F(b)(x_1{-}u_1, x_2{-}u_2)e^{\pi i(x_1(u_2-2ti)-x_2u_1)}
$$
of $u_2$ is holomorph on a strip 
$-\varepsilon<{\rm Im}~ u_2<2|t|{+}\varepsilon$ of the complex 
$u_2$-plane for some small $\varepsilon > 0$. 
Hence the integral of this function along the boundary of the rectangle 
with corners ${-}R,R,R{+}2ti,{-}R{+}2ti$ vanishes. In order to justify the 
substitution $u_2\rightarrow u_2{+}2ti$, it is sufficient to show that 
the corresponding integrals from $\pm R$ to $\pm R{+}2ti$ tend to zero 
as $R\rightarrow +\infty$. Using formula (\ref{F4a}) we estimate 
\begin{align*}
&\left|{\int\limits^{2t}_0}{\int\limits^\infty_{-\infty}} 
\F(a)(u_1,\pm R+si) \F(b)(x_1{-}u_1,x_2{-}(\pm R{+}si))
e^{\pi i(x_1(\pm R{+}si - 2ti)-x_2u_1)} ds du_1 \right|\\
&\le C\left|{\int\limits^{2t}_0}{\int\limits^\infty_{-\infty}} 
e^{-2\pi(\nu_1|u_1|+\nu_2 R)-2\pi(\nu^\prime_1|x_1-u_1|+
\nu^\prime_2|x_2\mp R|)+\pi x_1 (2t-s)} ds du_1 \right|\\
&\le C_{x_1,x_2} e^{-2\pi\nu_2 R},
\end{align*}
where $C$ and $C_{x_1,x_2}$ are not depending on $R$. Since 
$\nu_2 > 0$, the integral goes to zero if $R\rightarrow+\infty$. 
This proves formula (\ref{F10}) for real $t$. 

Next we prove (\ref{F5}) for real $t$. From the 
definition (\ref{C5}) of the product $\scriptstyle{\#}$ we obtain
\begin{align*}
&(e^{2\pi\cdot t Q_1} (a{\scriptstyle{\#}} b))(x_1, x_2) =\\
&\iiiint e^{2\pi t u_1} a(u_1,u_2) b(v_1, v_2) 
e^{4\pi i[(x_1-u_1)(x_2-v_2- ti/2)-(x_1-v_1)(x_2-u_2)]} du_1 du_2 dv_1 dv_2.
\end{align*}
Recall that $a,b\in\cD_{\mu,\nu}$ by assumption. Hence we have 
$a\in\cD(e^{\pm 2\pi\mu_1 \Q_1})$ and 
$b\in\cD(e^{\pm 2\pi\nu^\prime_2\cP_2})$. Since $|t|<\mu_1$ and 
$|t|<2\nu^\prime_2$, the latter implies that  
$a\in\cD(e^{2\pi t\Q_1})$ and $b\in\cD(e^{\pi t\cP_2})$. In fact, 
we even have that $e^{2\pi t\Q_1} a, e^{\pi t\cP_2}a \in
\cD_{\tilde{\mu},\tilde{\nu}}$ for certain $\tilde{\mu},
\tilde{\nu}\in\bbbr^{++}$. Equation (\ref{F5}) follows 
from the preceding formula by the formal substitution 
$v_2\rightarrow v_2+ it/2$.  In order to show that 
this formal replacement is justified we integrate in the 
complex $v_2$-plane along the boundary of the rectangle with corners 
$-R,R,R + it/2,-R+it/2$, where  $R>0$. To complete 
the proof, it suffices to show that the integrals from $\pm R$ to 
$\pm R+ it/2$ tend to zero as $R\rightarrow +\infty$. 
Indeed, using formula (\ref{F4a}) and the assumptions $|t|<\mu_1$ 
and $|t|<2\nu^\prime_2$, we estimate
\begin{align*}
&\biggl| {\int\limits^{t/2}_0} \iiint ds du_1 du_2 dv_1~   e^{2\pi t u_1} a(u_1,u_2)b(v_1, {\pm} R{+}si)\\ 
&\hspace{4.1cm}e^{4\pi i[(x_1-u_1)(x_2-(\pm R{+}si)-ti/2)-(x_1-v_1)(x_2-u_2)]} \biggl|\\
&\le C\biggl|{\int\limits^{t/2}_0} 
\iiint e^{2\pi t u_1} e^{-2\pi(\mu_1|u_1|+\mu_2|u_2|+
\mu^\prime_1|v_1|+\mu^\prime_2 R)} e^{4\pi (x_1-u_1)(t/2-s)} 
ds du_1 du_2 dv_1 \biggl|\\
&\le C^\prime_{x_1} e^{-2\pi\mu^\prime_2 R}\biggl| 
{\int\limits^\infty_{-\infty}} 
{\int\limits^{t/2}_0} e^{2\pi(2su_1 -\mu_1|u_1|} du_1 ds\biggl|\\
&\le C^{\prime\prime}_{x_1} e^{-2\pi\mu^\prime_2 R} {\int\limits^\infty_{-\infty}} 
e^{2\pi|t||u_1|-2\pi\mu_1|u_1|} du_1\le C^{\prime\prime\prime}_{x_1} 
e^{-2\pi\mu^\prime_2 R},
\end{align*}
where $C,C^\prime_{x_1},C^{\prime\prime}_{x_1}, 
C^{\prime\prime\prime}_{x_1}$ are numbers not depending on $R$.
Since $\mu^\prime_2{>}0$, the integral goes to zero if $R{\rightarrow}{+}\infty$. This completes the proof of (\ref{F5}) for real $t$.

For imaginary $t$ the above reasoning works as well. In this 
case we are lead to {\it real} translations of $u_2$ and $v_2$,
respectively, which are possible by the translation invariance 
of the Lebesgue measure. The case of general $t\in\bbbc$ follows 
by combining the real and the imaginary cases. \hfill$\Box$

\mn
For $\mu, \nu \in \bbbr^{++}$, let $\cD^{\mu,\nu}$ denote 
the intersection of domains $\cD_{\mu^\prime,\nu^\prime}$ (see (\ref{ds})), 
where $\mu^\prime ,\nu^\prime \in \bbbr^{++}$, 
$\mu^\prime_j < \mu_j$, $\nu^\prime_j <\nu_j$ for $j=1,2$. 

\mn
{\bf Corollary 12.} {\it Let $\mu, \nu \in \bbbr^{++}$. If 
$\mu_1 < 2 \nu_2$ and $\mu_2 < 2 \nu_1$, then $\cD^{\mu,\nu}$ is a 
$\ast$-algebra with product ${\scriptstyle{\#}}$ and involution $\ast$ 
defined (\ref{C5}) and (\ref{C5a}), respectively. In particular, 
$\gA(\bbbr^2)$ is a $\ast$-algebra. }

\mn
{\bf Proof.} Since $\mu_1 < 2 \nu_2$ and $\mu_2 < 2 \nu_1$, we conclude 
from formulas (\ref{F5}), (\ref{F7}), (\ref{F9}), (\ref{F11}) and (\ref{ds}) 
that $a,b \in \cD^{\mu,\nu}$ implies that 
$a{\scriptstyle{\#}}b \in  \cD^{\mu,\nu}$.
By (\ref{ds}) it is obvious that $a^\ast \in  \cD^{\mu,\nu}$ for 
$a \in \cD^{\mu,\nu}$. Thus $\cD^{\mu,\nu}$ is a 
$\ast$-algebra. Since $\gA(\bbbr^2)$ is the intersection of all domains 
$\cD^{\mu,\nu}$, $\gA(\bbbr^2)$ is a $\ast$-algebra as well.  \hfill$\Box$

\mn 
{\bf Lemma 13.} {\it Suppose that $\mu, \nu \in \bbbr^{++}$. 
Let $\parallel \cdot \parallel$ denote the norm of $L^2(\bbbr^2)$.
If $a,b \in \cD_{\mu,\nu} = \cD(S(e^{2\pi\mu \Q})S(e^{2\pi\nu\cP}))$, 
then $a, b, \bar{b}{\scriptstyle{\#}} a\in\sS(\bbbr)$ and we have}
\begin{align}\label{scalar}
\iint a(x_1,x_2) {\overline b(x_1,x_2)} dx_1 dx_2 &=
\iint ({\overline b} {\scriptstyle{\#}} a)(x_1,x_2) dx_1 dx_2,\\
\label{hsnorm}
\parallel a {\scriptstyle{\#}} b \parallel &\le \parallel a \parallel\, 
\parallel b \parallel.
\end{align}

\mn
{\bf Proof.} By Lemma 9(ii) and Corollary 12, we have $a, b \in \sS(\bbbr^2)$ and so $\bar{b}{\scriptstyle{\#}} a\in\sS(\bbbr^2)$. From Proposition 2 in [St], p. 555, it follows that 
$Op(a)$ and $Op(b)$ are Hilbert-Schmidt operators and 
\begin{equation}\label{hsscalar}
\langle Op(a),Op(b) \rangle_{HS} = {\rm Tr}~ Op(b)^\ast Op(a) 
= (a,b)~{\rm and}~
\parallel Op(a) \parallel_{HS} = \parallel a \parallel,
\end{equation}
where ${\rm Tr}$ is the trace and $\langle \cdot,\cdot \rangle_{HS}$ and 
$\parallel \cdot \parallel_{HS}$ denote scalar product and 
Hilbert-Schmidt norm of Hilbert-Schmidt operators, respectively.
Using formulas (\ref{hsscalar}) and (\ref{oppro}) and the 
submultiplicativity of the Hilbert-Schmidt norm we obtain
\begin{align*}
\parallel a {\scriptstyle{\#}} b \parallel &= \parallel Op(a{\scriptstyle{\#}} b) \parallel_{HS} = 
\parallel Op(a)Op(b) \parallel_{HS} \\
&\le \parallel Op(a) \parallel_{HS}\parallel Op(b) \parallel_{HS}=
\parallel a \parallel  \parallel b  \parallel.
\end{align*}
This proves (\ref{hsnorm}).

Put $c:= b^\ast {\scriptstyle{\#}} a$. By (\ref{C4}), the operator $Op(c)$ is an integral operator with kernel 
\begin{equation}\label{kern}
K_c(x_1,x_2) = \int c({\tfrac{1}{2}} (x_1{+}x_2),t) e^{2\pi i(x_1{-}x_2)t} dt.
\end{equation}
Since $c \in \sS(\bbbr^2)$, the function $d$ defined by 
$d(y_1,y_2) = \int c(y_1,t) e^{2\pi i y_2t} dt$ is in $\sS(\bbbr^2)$ and 
so is the function $K_c(x_1,x_2) = d({\tfrac{1}{2}} (x_1{+}x_2),x_1{-}x_2)$. 
It is well-known that any integral operator with kernel in the 
Schwartz space $\sS(\bbbr^2)$ is a trace class operator on $L^2(\bbbr)$ and 
that its trace is given by the integral over the diagonal. Using this 
fact and formulas (\ref{oppro}) and (\ref{kern}) we get
\begin{align*}
&{\rm Tr}~ Op(b)^\ast Op(a) = {\rm Tr}~ Op(b^\ast {\scriptstyle{\#}} a) = {\rm Tr}~ Op(c)\\
=&\int K_c(x_1,x_1) dx_1 = \iint c(x_1,x_2) dx_1 dx_2 \\
=&\iint a(x_1, x_2)\overline{b}(x_1,x_2) dx_1 dx_2.
\end{align*}
Comparing the latter with (\ref{hsscalar}), formula
(\ref{scalar}) follows. \hfill $\Box$
\sn 

Our next proposition says that $\gA(\bbbr^2)[\tau]$ is a 
Frechet $\ast$-algebra with approximate identity.

\mn
{\bf Proposition 14.} (i) {\it Provided with the product $\scriptstyle{\#}$, the 
involution $\ast$ and the locally convex topology $\tau,\gA(\bbbr^2)$ is a 
Frechet topological $\ast$-algebra.}

\sn
(ii) {\it Set $f_\varepsilon (x_1, x_2):=e^{-\pi\varepsilon(x^2_1+x^2_2)}.$ For each $a\in\gA(\bbbr^2)$, we have
\begin{equation}\label{apid}
{\lim\limits_{\varepsilon\rightarrow +0}} f_\varepsilon{\scriptstyle{\#}} a={\lim\limits_{\varepsilon\rightarrow +0}} a{\scriptstyle{\#}} f_\varepsilon =a
\end{equation}
in the locally convex space $\gA(\bbbr^2)[\tau]$.}

\mn
{\bf Proof.} (i): Recall that by definition the topology $\tau$ is generated by 
the family of norms $\parallel e^{2\pi c\Q} e^{2\pi d\cP}{\cdot}\parallel$, 
$c,d\in\bbbr^2$, where $\parallel\cdot\parallel$ is the norm of 
$L^2(\bbbr^2)$. Fix $c,d\in\bbbr^2$ and put
\begin{equation}\label{mu}
\mu=(\mu_1,\mu_2), \mu_1:=d_1+c_1/2, \mu_2 :=d_2-c_2/2.
\end{equation}
By (\ref{F6}), (\ref{F8}), (\ref{F9}), (\ref{F11}) and (\ref{hsnorm}), we 
obtain
\begin{equation}\label{mult}
\parallel e^{2\pi c\Q} e^{2\pi d\cP}(a{\scriptstyle{\#}} b)\parallel=
\parallel (e^{2\pi\mu\cP} a){\scriptstyle{\#}} (e^{2\pi c\Q} 
e^{2\pi d\cP} b)\parallel
\le \parallel e^{2\pi \mu\cP} a\parallel \parallel e^{2\pi c\Q} 
e^{2\pi d\cP} b\parallel
\end{equation}
for $a,b\in\gA(\bbbr^2)$. Since $\parallel e^{2\pi c\Q} e^{2\pi d\cP} 
a^\ast\parallel=\parallel e^{2\pi c\Q} e^{-2\pi d\cP} a\parallel$, 
product and involution are $\tau$-continuous, so $\gA(\bbbr^2)[\tau]$ is 
indeed a topological $\ast$-algebra. Since $\gA(\bbbr^2)[\tau]$ is a 
Frechet space as noted above, it is a Frechet topological $\ast$-algebra. 

(ii):  Let $b\in\gA(\bbbr^2)$ and $\mu=(\mu_1,\mu_2)\in\bbbr^2$. Our aim is to prove by explicit estimations that 
\begin{equation}\label{limit}
{\lim\limits_{\varepsilon\rightarrow +0}} (e^{2\pi\mu\cP} f_\varepsilon){\scriptstyle{\#}} b=b
\end{equation} 
in $L^2 (\bbbr^2)$. Note first that from the well-known equation
\begin{equation}\label{int}
\int e^{-\pi(s-c)^2\varepsilon-2\pi its} ds=
\varepsilon^{-1/2} e^{-\pi\varepsilon^{-1} t^2} e^{-2\pi itc},~ c\in\bbbc,
\end{equation}
we obtain that
$$
\F (e^{2\pi\mu\cP} f_\varepsilon)(x_1,x_2)=\varepsilon^{-1} e^{-\pi\varepsilon^{-1}(x^2_1+x^2_2)} e^{2\pi(\mu_1x_1+\mu_2x_2)}.
$$
Using the latter formulas and the definition (\ref{C5}) of the product ${\scriptstyle{\#}}$ we compute
\begin{align}\label{ab1}
&|b(x_1,x_2)-((e^{2\pi\mu\cP} f_\varepsilon){\scriptstyle{\#}} b)(x_1,x_2)|\nonumber\\
&=\Big| b(x_1,x_2)-4\iint \F(e^{2\pi\mu\cP} f_\varepsilon)(2x_2{-}2v_2, 2v_1-2x_1) b(v_1,v_2) e^{4\pi i(v_1x_2{-}x_1v_2)} dv_1 dv_2\Big|\nonumber\\
&=\Big| {\tfrac{4}{\varepsilon}}\iint (b(x_1,x_2)-b(v_1,v_2)) dv_1 dv_2\\
&\hspace{1.8cm} e^{-4\pi\varepsilon^{-1}((x_2-v_2)^2+(v_1-x_1)^2)+4\pi((x_2-v_2)\mu_1+(v_1-x_1)\mu_2)+4\pi i(v_1x_2-x_1v_2)} \Big| \nonumber\\
&\le {\tfrac{4}{\varepsilon}} \iint e^{-4\pi\varepsilon^{-1} ((x_1-v_1)^2+(x_2-v_2)^2)+4\pi((x_2-v_2)\mu_1-(x_1-v_1)\mu_2)}|b(x_1,x_2)-b(v_1,v_2)| dv_1 dv_2.
\end{align}
Fix a number $\delta >0$. 
Since $b\in\gA(\bbbr^2)$,(\ref{F4a}) holds for arbitrary $\nu,\mu\in\bbbr^{++}$. Hence there exists $M\in\bbbr$ such that for $(x_1,x_2)\in\bbbr^2$,
\begin{equation}\label{bound}
e^{(4\pi |\mu_1|+1)|x_1|+(4\pi|\mu_2|+1)|x_2|}|b(x_1,x_2)|\le M.
\end{equation}
Further, since $e^{|x_1|+|x_2|} b(x_1,x_2)\rightarrow 0$ as 
$|x_1|+|x_2|\rightarrow +\infty$ by (\ref{bound}), the function 
$c(x_1, x_2):=e^{|x_1|+|x_2|} b(x_1,x_2)$ is uniformly continuous on 
$\bbbr^2$. Thus there exists $\delta_1$ such that $1>\delta_1>0$ and
$$
(1+ eM) e^{4\pi(|\mu_1|+|\mu_2|)}|c(x_1,x_2)-c(v_1,v_2)|<\delta
$$
for $(x_1-v_1)^2 + (x_2-v_2)^2<\delta^2_1$. From this and (\ref{bound}) we easily derive that
\begin{equation}\label{ab2}
e^{|x_1|+|x_2|} e^{4\pi((x_2-v_2)\mu_2-(x_1-v_1)\mu_1)}|b(x_1,x_2)-b(v_1,v_2)|<\delta
\end{equation}
when $(x_1-v_1)^2+(x_2-v_2)^2<\delta^2_1$.

Next we turn to the domain where $(x_1-v_1)^2+(x_2-v_2)^2\ge \delta^2_1$.
Obviously, there exists $K\in\bbbr$ such that
\begin{equation}\label{Ka}
-\pi\varepsilon^{-1} (t^2_1+t^2_2)+ 
4\pi (|\mu_1|+1)|t_1|+(|\mu_2|+1)|t_2|\le K
\end{equation}
for all $(t_1,t_2)\in\bbbr^2$ and $1>\varepsilon>0$. If $(x_1-v_1)^2+(x_2-v_2)^2\ge \delta^2_1$, then by (\ref{bound}) and (\ref{Ka}) we obtain that
\begin{align}\label{ab3}
&e^{|x_1|+|x_2|} e^{-2\pi\varepsilon^{-1}((x_1-v_1)^2+(x_2-v_2)^2)} e^{4\pi((x_2-v_2)\mu_1-(x_1-v_1)\mu_2)}
|b(x_1,x_2)-b(v_1,v_2)|\nonumber\\
&\qquad \le e^{-\pi\varepsilon^{-1}
((x_1-v_1)^2+(x_2-v_2)^2)} e^K 2M\le e^{-\pi\varepsilon^{-1}\delta^2_1} 
e^K 2M<\delta
\end{align}
for sufficiently small $\varepsilon>0$. 

Using the relation 
$$
{{\tfrac{4}{\varepsilon}}}\iint e^{-2\pi \varepsilon^{-1}((x_1-v_1)^2+(x_2-v_2)^2)} dv_1 dv_2=2
$$
by (\ref{int}), it follows from estimates (\ref{ab1}), (\ref{ab2}) and (\ref{ab3}) that 
$$
|b(x_1,x_2)-((e^{2\pi\mu\cP} 
f_\varepsilon){\scriptstyle{\#}} b)(x_1,x_2)|\le 3\delta e^{-|x_1|-|x_2|}
$$
and so $\parallel b-(e^{2\pi\mu\cP} 
f_\varepsilon){\scriptstyle{\#}} b\parallel\le 12\delta $
for small $\varepsilon> 0$. This proves (\ref{limit}).

Now let $a\in\gA(\bbbr^2)$ and $c,d\in\bbbr^2$. Let $\mu$ be as in (\ref{mu}). Applying (\ref{limit}) with $b=e^{2\pi c\Q} e^{2\pi d\cP} a$, we get
\begin{equation*}
\parallel e^{2\pi c\Q} e^{2\pi d\cP} 
(f_\varepsilon{\scriptstyle{\#}} a-a)\parallel
=\parallel (e^{2\pi\mu\cP} f_\varepsilon){\scriptstyle{\#}}(e^{2\pi c\Q} 
e^{2\pi d\cP} a) - e^{2\pi c\Q} e^{2\pi d\cP} a\parallel\rightarrow 0
\end{equation*}
as $\varepsilon \rightarrow +0$. This proves that 
${\lim\limits_{\varepsilon\rightarrow +0}} f_\varepsilon{\scriptstyle{\#}} a
=a$ in $\gA(\bbbr^2)[\tau]$. Applying the involution we obtain the second 
equality in (\ref{apid}).\hfill $\Box$

\mn
Remark 2.  Upon scaling and multiplying by parameters, the operators 
$Op(f_\varepsilon)$, $\varepsilon>0$, form the so-called Hermite 
semigroup $e^{-2\pi t(\cP^2+\Q^2)}, t>0$, acting on the Hilbert space 
$L^2(\bbbr^2)$, see [Fo], pp. 236--238.

\sn
In this paper we shall mainly use the symbol algebra $\gA(\bbbr^2)$. 
However, for most considerations it suffices to work with the smaller symbol algebras
\begin{align*}
&\gA_{ex} (\bbbr^2):=\Lin\{ e^{-\varepsilon_1 x^2_1-\varepsilon_2 x^2_2+c_1x_1+c_2x_2}; \varepsilon_1>0,\varepsilon_2>0, c_1,c_2\in\bbbc\},\\
&\gA_{pex}(\bbbr^2):=\Lin\{x^{n_1}_1 x^{n_2}_2 e^{-\varepsilon_1x^2_1-\varepsilon_2x^2_2+c_1x_1+c_2x_2};\varepsilon_j>0, c_j\in\bbbc,n_j\in\bbbn_0\}.
\end{align*}
Both $\gA_{ex} (\bbbr^2)$ and $\gA_{pex} (\bbbr^2)$ are $\ast$-algebras with multiplication (\ref{C5}) and involution (\ref{C5a}). In order to prove this assertion it is sufficient to show that $a{\scriptstyle{\#}} b$ is in $\gA_{ex}(\bbbr^2)$ resp. $\gA_{pex}(\bbbr^2)$ when $a$ and $b$ are so. In the case of $\gA_{ex}(\bbbr^2)$ this can be verified by direct computation of the twisted product $a{\scriptstyle{\#}} b$ using formula (\ref{int}). From formula (3) in [GV] it follows at once that $a{\scriptstyle{\#}} b\in\gA_{pex}(\bbbr^2)$ for $a,b\in\gA_{pex}
(\bbbr^2)$.

\bn
{\bf 2. The Coordinate Algebra ${\bf \cO(\bbbc^2_q)}$ of the Quantum Plane} 

\sn
{\it 2.1 $ \cO(\bbbc^2_q)$ as a left module 
algebra of $\U_q(gl_2)$}

\sn
Let $\ang$ be the left action of $\U_g(gl_2)$ on $\cO(\bbbc^2_q)$ associated with the right coaction of $\cO(GL_q(2))$ defined by (\ref{A1}). From (\ref{A1}) and (\ref{A2}) we then obtain
\begin{align}\label{A3}
&K_1\ang x=q^{-1/2} x,K_1\ang y = y, K_2 \ang x =x, K_2\ang y=q^{-1/2}y,\\
\label{A3a}
&E \ang x =y, E\ang y=0,F\ang x=0,F\ang y=x.
\end{align}

Moreover, since $\varepsilon (K_1)= \varepsilon (K_2) =1$ and 
$\varepsilon (E)=\varepsilon (F)=0$, we also have
\begin{equation}\label{A4}
K_1 \ang 1=1, K_2 \ang 1 =1, E\ang 1=0, F\ang 1=0.
\end{equation}

The following proposition derives the action of the generators 
$K_1, K_2, E,F$ of $\U_q(gl_2)$ on general elements of the 
algebra $\cO(\bbbc^2_q)$. We set
$$
D_{q^{-2}}(f)(x) := \frac{f(x)-f(q^{-2}x)}{(1-q^{-2})x}.
$$

\mn
{\bf Proposition 15.} {\it If $g$ and $h$ are complex polynomials in a single variable, then we have}
\begin{align}\label{A5}
K_1\ang (g(x) h(y))&= g(q^{-1/2}x) h(y),~
K_2\ang (g(x) h(y))= g(x) h(q^{-1/2} y),\\
\label{A6}
E~~\ang (g(x) h(y))&= q^{-1/2} D_{q^{-2}} (g(q^{1/2}{\cdot}))(x) y h(q^{1/2} y),\\
\label{A7}
F~~\ang (g(x) h(y))&= q^{-1/2} g(q^{1/2} x) x D_{q^{-2}} (h(q^{1/2}{\cdot}))(y).
\end{align}
{\bf Proof.} Since $\cO(\bbbc^2_q)$ is a 
$\U_q(gl_2)$-comodule algebra, equation (\ref{A8}) holds.
The assertion follows from this equation combined with formulas 
(\ref{A3}) and (\ref{A4}). For the generators $K_1$ and $K_2$ 
this is obvious. We carry out the proof of formula (\ref{A7}). 
The proof of formula (\ref{A6}) is similar.

Since $\Delta (F)=F\otimes K +K^{-1} \otimes F$, 
it follows from (\ref{A8}) that
\begin{equation}\label{A9}
F\ang (zz^\prime)=(F\ang z)(K\ang z^\prime)+(K^{-1}\ang z)(F\ang z^\prime).
\end{equation}
Recall that $F\ang x=0$ and $F\ang 1=0$ by (\ref{A3}) and (\ref{A4}). Using these facts it follows from (\ref{A9}) by induction on $n$ that 
$F\ang x^n=0$ for $n\in\bbbn_0$. Thus we have $F\ang g(x)=0$. 

Since $K^{-1}\ang g(x)=K_1^{-1} \ang K_2 \ang g(x) = K_1^{-1} \ang g(x) =
g(q^{1/2} x)$, we derive from (\ref{A9}), applied to $z=g(x)$ and 
$z^\prime =h(y)$, that
\begin{equation}\label{A10}
F\ang (g(x) h(y))=g(q^{1/2} x)(F\ang h(y)).
\end{equation}
Therefore, in order to prove (\ref{A7}) it suffices to show that
\begin{equation}\label{A11}
F\ang y^n = q^{-1/2}  q^{n/2} x D_{q^{-2}}(y^n),~ n\in\bbbn.
\end{equation}
We prove (\ref{A11}) by induction on $n$. If $n=1$, then (\ref{A11}) 
is true by (\ref{A3}). If (\ref{A11}) is valid for a, then it follows 
from (\ref{A9}) and (\ref{A3}) that
\begin{align*}
F\ang y^{n+1} &=(F\ang y)(K\ang y^n)+(K^{-1}\ang y)(F\ang y^n)\\
&=x(q^{n/2}y^n)+(q^{-1/2} y)(q^{-1/2}(1-q^{-2})^{-1}(1-q^{-2n})q^{n/2} xy^{n-1})\\
&=q^{n/2}(1+q^{-2}(1-q^{-2})^{-1}(1-q^{-2n}))xy^n\\
&=q^{-1/2}(1-q^{-2})^{-1}(1-q^{-2(n+1)}) q^{(n+1)/2} xy^n\\
&=q^{-1/2} q^{(n+1)/2} x D_{q^{-2}}(y^{n+1}),
\end{align*}
which proves (\ref{A11}) in the case of $n+1$.\hfill$\Box$

\mn

For $z\in O(\bbbc^2_q)$, we define
\begin{equation}\label{dxy}
\cD^q_x (z)=Ky^{-1} E^\prime\ang z,~\cD^q_y(z)=K x^{-1} F^\prime\ang z,
\end{equation}
where the elements $y^{-1}$ and $x^{-1}$ of $\hat{\cO}(\bbbc^2_q)$ act 
by left multiplication on $\cO(\bbbc^2_q)$. From (\ref{A6}) and (\ref{A7}) 
we obtain
\begin{align}\label{daction}
&\cD^q_x (g(x)h(y))=q^{1/2}\lambda D_{q^{-2}}(g)(qx) h(q y),\\
&\cD^q_y(g(x)h(y))=q^{-1/2}\lambda g(x) D_{q^{-2}}(h)(qy).
\end{align}
for polynomials $g$ and $h$. 

For $r,s\in\bbbn_0$, let $\sigma_{rs}$ denote 
the automorphism of the algebra $\cO(\bbbc^2_q)$ defined by 
$\sigma_{rs}(z)=K^r_1K^s_2\ang z,z\in \cO(\bbbc^2_q)$. From formulas 
(\ref{A5})-(\ref{A7}) or (\ref{daction}) we easily derive that 
$\cD^q_x$ is a $(\sigma_{-2,0}, \sigma_{2,-2})$-derivation and 
$\cD^q_y$ is a $(\sigma_{0,-2},\sigma_{2,2})$-derivation of the 
algebra $\cO(\bbbc^2_q)$, that is, for $z_1,z_2\in \cO(\bbbc^2_q)$ we have
\begin{align*}
&\cD^q_x(z_1z_2)=(K^{-2}_1\ang  z_1)\cD^q_x(z_2)+\cD^q_x(z_1)
(K^2_1 K^{-2}_2\ang z_2),\\
&\cD^q_y(z_1z_2)=(K^{-2}_2\ang z_1)\cD^q_y(z_2)+\cD^q_y(z_1)
(K^2_1K_2^2\ang z_2).
\end{align*}
In the limit $q \to 1$ the preceding equations go into the Leibniz rule. 
We shall consider the linear mappings $\cD^q_x$ and $\cD^q_y$ as $q$-deformed partial 
derivatives of the algebra $\cO(\bbbc^2_q)$.

\mn
{\it 2.2 Covariant Differential Calculus on $\cO(\bbbc^2_q)$}
\sn

As shown in [PW] and [WZ], there are two distinguished first 
order differential calculi $\Gamma_+$ and $\Gamma_-$ on $\cO(\bbbc^2_q)$. For both calculi, the set of differentials $\{dx, dy\}$ is a basis for the right (and for the left) $\cO(\bbbc^2_q)$-module of first order forms. Therefore, for any $z\in\cO(\bbbc^2_q)$ there exist uniquely determined elements $\partial_x(z)$ and $\partial_y(z)$ of $\cO(\bbbc^2_q)$, called {\it partial derivatives} of $z$, such that
\begin{equation}\label{F0}
dz=dx{\cdot}\partial_x(z)+dy{\cdot} \partial_y(z).
\end{equation}
The bimodule structures of the calculi $\Gamma_+$ and $\Gamma_-$ are 
described by the following commutation relations:
\begin{align}\label{F1+}
\Gamma_+:~ &x\rd y=q \rd y{\cdot} x + (q^2-1)\rd x{\cdot} y,~ y\rd x= 
q\rd x{\cdot} y,\\
\label{F2+}
&x\rd x=q^2 \rd x{\cdot} x,~ y \rd y= q^2 \rd y{\cdot} y.\\
\label{F1-}
\Gamma_- :~ & y\rd x=q^{-1} \rd x{\cdot} y + (q^{-2} -1)\rd y{\cdot} x,~ 
x\rd y= q^{-1} \rd y{\cdot} x,\\
\label{F2-}
&x\rd x=q^{-2} \rd x{\cdot} x,~ y \rd y=q^{-2} \rd y{\cdot} y.
\end{align}
From these relations we see that $\eta_+:=y^{-2} xdx$ and 
$\eta_-:=x^{-2}y dy$ are non-zero central elements of the bimodules 
$\Gamma_+$ and $\Gamma_-$, respectively. Recall that an element 
$\eta$ of a bimodule 
over an algebra $\Z$ is called {\it central} if $\eta z=z\eta$ for all 
$z\in\Z$.\\ 
Note that the relations for $\Gamma_+$ go into the relations 
of $\Gamma_-$ if we interchange the coordinates $x$ and $y$ and the numbers 
$q$ and $q^{-1}$. The partial derivatives $\partial_x$ and $\partial_y$, 
considered as linear mappings of $\cO(\bbbc^2_q)$, and the coordinate 
functions $x$ and $y$, acting on $\cO(\bbbc^2_q)$ by left multiplication, 
satisfy the relations:
\begin{equation*}
\begin{split}
\Gamma_+ :~ &\partial_x y=qy\partial_x, ~\partial_y x=qx\partial_y,\\
&\partial_x x-q^2x\partial_x=1+(q^2-1)y\partial_y,
~\partial_y y-q^2y\partial_y=1.\\
\Gamma_-:~ &\partial_x y=q^{-1} y\partial_x,
~\partial_y x=q^{-1}x\partial_y,\\
&\partial_x x-q^{-2} x\partial_x=1,
~ \partial_y y-q^{-2} y\partial_y=1+(q^{-2}-1)x\partial_x.
\end{split}
\end{equation*}
From these formulas one derives by induction the expressions 
for the actions of $\partial_x$ and $\partial_y$ on general elements of 
$\cO(\bbbc^2_q)$. If $g$ and $h$ are complex polynomials in a single 
variable, then we have:
\begin{align}\label{F3}
&\Gamma_+:~\partial_x(g(y)h(x))=g(qy) D_{q^{-2}}(h)(x), 
\partial_y(g(y)h(x))=D_{q^{-2}}(g)(y) h(x),\\
&\Gamma_-:~\partial_x(g(x)h(y))=D_{q^{-2}}(g)(x)h(y),
\partial_y(g(x)h(y))=g(q^{-1}x)                 D_{q^{-2}}(h)(y).
\end{align}

All these facts and formulas are well-known. We now give another 
description of these calculi. Let $\Omega$ be the free bimodule of the 
localization algebra $\check{\cO}(\bbbc^2_q)$ generated by a central vector 
space $V$. 
That is, $\Omega$ is the vector space $\check{\cO}(\bbbc^2_q)\otimes V$ 
with bimodule structure given by 
\begin{equation}\label{bimod}
u\Big({\sum\nolimits_j} z_j\otimes e_j\Big)v:=
{\sum\nolimits_j} u z_jv\otimes e_j,
\end{equation}
 where $u,v,z_j\in 
\check{\cO}(\bbbc^2_q), e_j\in V$. For notational simplicity, we write $ze$ instead of $z\otimes e$, where $z\in\check{\cO}(\bbbc^2_q)$ 
and $e\in E$. Fix two elements $e_1, e_2\in V$ and put
\begin{equation*}\begin{split}
&\omega_+ := q^{-2}y^2x^{-2}  e_1+x^{-2} e_2,~ \omega_-=q^2x^2y^{-2} e_1+ y^{-2} e_2,\\
&d_+z :=\omega_+z-z\omega_+,~ d_-z:=\omega_-z-z\omega_-,~ z\in\cO(\bbbc^2_q).
\end{split}
\end{equation*}
Let us abbreviate $\Z:=\cO(\bbbc^2_q)$. Obviously,
 $\tilde{\Gamma}_\varepsilon := \Z{\cdot} d_\varepsilon \Z{\cdot}\Z$ 
is a $\Z$-bimodule and the mapping $d_\varepsilon 
:\Z\rightarrow\tilde{\Gamma}_\varepsilon$ satisfies the Leibniz rule 
for $\varepsilon{={+,-}}$. Thus, the pair $(\tilde{\Gamma}_\varepsilon, 
d_\varepsilon)$ is a first order differential calculus over 
the algebra $\Z=\cO(\bbbc^2_q)$. For the differentials of the 
coordinate functions we obtain
\begin{align}
&d_+ x=(q^{-2}-1)y^2x^{-1} e_1,~ d_+ y
=(q^{-2}-1)q^{-2} y^3 x^{-2} e_1+(q^{-2}-1)yx^{-2} e_2,\\
\label{F4}
&d_- x=(q^2-1)q^2x^3y^{-2} e_1+(q^2-1)xy^{-2} e_2,~
 d_-y=(q^2-1)x^2 y^{-1} e_1.
\end{align}

\mn
{\bf Lemma 16.} {\it Suppose that the elements $e_1$ and $e_2$ 
are linearly independent. Then the first order differential calculi 
$\Gamma_\varepsilon$ and $\tilde{\Gamma}_\varepsilon, \varepsilon={+,-}$, are 
isomorphic.}

\mn
{\bf Proof.} Since $\{dx, dy\}$ is a free left $\cO(\bbbc^2_q)$-module 
basis of $\Gamma_\varepsilon$, there is a well-defined left 
$\cO(\bbbc^2_q)$-module homomorphism $\psi_\varepsilon:\Gamma_\varepsilon\rightarrow\tilde{\Gamma}_\varepsilon$ such that
$$
\psi_\varepsilon (u dx + v dy)=u d_\varepsilon x + v d_\varepsilon 
y,~ u, v \in \cO (\bbbc^2_q).
$$
In order to prove that $\psi_\varepsilon$ is an $\cO(\bbbc^2_q)$-bimodule 
homomorphism, it suffices to show that the relations (\ref{F1+}) and
(\ref{F2+}) resp. (\ref{F1-}) and (\ref{F2-}) hold also in 
$\Omega_+$ and $\Omega_-$. As a sample, we verify the first relation of (\ref{F1-}). The other relations follow by similar straightforward computations. Using formulas (\ref{F4}) and the commutation rules in the algebra $\check{\cO} (\bbbc^2_q)$, we obtain
\begin{align*}
&q^{-1} d_-x{\cdot} y+(q^{-2}-1) d_- y{\cdot} x=\\
&(q^2-1)(qx^3y^{-2} e_1 y+q^{-1} xy^{-2} e_2y+(q^{-2}-1)x^2y^{-1} e_1 x)=\\
&(q^2-1)((qx^3y^{-1}+(q^{-2}-1)qx^3y^{-1})e_1+y xy^{-2} e_2)=yd_-x.
\end{align*}
From the construction it is clear that $\psi_\varepsilon$ is a surjective
 FODC homomorphism. We show that $\psi_-$ is injective and suppose that 
$ud_-x+v d_-y=0$ for some elements $u,v\in\cO(\bbbc^2_q)$. Inserting 
the expressions from (\ref{F4}) and using the assumption that $e_1$ and 
$e_2$ are linearly independent, we get
$$
uq^2x^3y^{-2}+vx^2y^{-1}=0, vxy^{-2}=0
$$
which in turn implies that $u=v=0$. The proof for $\psi_+$ is 
similar.\hfill$\Box$

We shall identify the isomorphic calculi $\Gamma_\varepsilon$ and 
$\tilde{\Gamma}_\varepsilon$. The above approach to the calculi $\Gamma_\varepsilon$ is convenient for many purposes. Among others, it allows us easily to extend these calculi to larger algebras.\\
The partial derivatives $\partial_x$ and $\partial_y$ can be also 
expressed in terms of the action of the generators of $\U_q(gl_2)$. 
Combining the formulas (\ref{A6}), (\ref{A7}) and (\ref{F3}) we 
obtain for the calculus $\Gamma_-$ the relations 
$$
\partial_x (z)=q^{\frac{3}{2}} y^{-1} EK^3_1K_2 \ang z,
~\partial_y (z) =q^{\frac{1}{2}} 
x^{-1} FK^3_1 K_2 \ang z
$$
for $z \in \cO(\bbbc^2_q)$, 
where the action of $E,F,K_1, K_2$ is given by Proposition 15 and the 
elements $y^{-1}$ and $x^{-1}$ act by left multiplication 
on $\cO(\bbbc^2_q)$.

\mn
{\bf 3. An Auxilary ${\bf \ast}$-algebra ${\bf \W}$ }

\sn
{\it 3.1} {\it The $\U_q(gl_2(\bbbr))$-module ${\bf \ast}$-algebra $\W$}
\sn

Let $\W$ denote the $\ast$-algebra generated by the operators
\begin{equation}\label{B4}
W(s,t):=e^{2\pi i(s\alpha \Q+t\beta\cP)},~~ s,t \in \bbbc.
\end{equation}
These operators satisfy the relations
\begin{align}\label{B5}
W(s_1,t_1)W(s_2,t_2)&=e^{\pi i\gamma (s_2t_1-s_1t_2)} W(s_1{+}s_2,t_1{+}t_2), \\
\label{B6}
W(s,t)^\ast&=W({-}\overline{s},{-}\overline{t})
\end{align}
for $s_1,t_1,s_2,t_2,s,t\in\bbbc$. By Lemma 4, the operator $W(s,t)$ acts as
\begin{equation}\label{B7}
(W(s,t)f)(x)=e^{2\pi is\alpha x+\pi ist\gamma} f(x{+}t).
\end{equation}
Equations (\ref{B5}) and (\ref{B7}) hold for vectors contained in the 
corresponding operator domains. For instance, they hold on the domains 
$\cD_\delta$, where $\delta > 0$, and $\gA(\bbbr^2)$ in the Hilbert space 
$L^2(\bbbr^2)$. Each of the dense subspaces is an invariant dense core for 
all operators $W(s,t)$. 
From (\ref{B1}) and (\ref{B7}) we see that  
\begin{equation}\label{B8}
W({-}i,0)=X, W(0,{-}i)=Y, W(s,0)=X^{is}, W(0,t)=Y^{is},~ s\in\bbbc.
\end{equation}

Our next aim is to define a left action of the 
Hopf algebra $\U_q(gl_2)$ on $\W$. Let us identify 
the generators $x$ with $X=W({-}i,0)$ and 
$y$ with $Y=W(0,{-}i)$. Then $\cO(\bbbr^2_q)$ becomes a 
$\ast$-subalgebra of $\W$. We now use formulas (\ref{A5})-(\ref{A7}) 
(which have been proved only for polynomials $g$ and $h$!) as a motivation and extend them formally to the functions $g(X)=X^{is}$ and 
$h(Y)=Y^{it},s,t\in\bbbc$, of the positive self-adjoint operators $X$ and 
$Y$ defined by (\ref{B1}). Throughout we interpret expressions $(q^{k/2}X)^{is}$ and $(q^{k/2}Y)^{it}$ as $e^{-\pi  k\gamma s} e^{2\pi is\alpha\Q}$ and $e^{-\pi k\gamma t} e^{2\pi it\beta\cP}$, respectively, for 
$k=0,\pm 1,3$ and $s,t\in\bbbc$. Recall that 
\begin{equation}\label{B9}
W(s,t)=e^{\pi i\gamma st} W(s,0) W(0,t)=e^{\pi i\gamma st} X^{is} Y^{it}
\end{equation}
by (\ref{B4}). Applying now (\ref{A5}) and (\ref{A7}) formally (!) and using 
(\ref{B9}) we derive
\begin{align*}
\lambda  F\ang W(s,t)&=\lambda e^{\pi i\gamma st} ((F\ang X^{is})(K\triangleright Y^{it})+(K^{-1}\ang X^{is})(F\ang YX^{it})\nonumber\\
&=\lambda e^{\pi i\gamma st}(0+(q^{1/2} X)^{is} q^{-1/2} X D_{q^{-2}}((q^{1/2} Y)^{it}))\nonumber\\
&=q^{1/2}(1-q^{-2})e^{\pi i\gamma st} e^{-\pi\gamma s} X^{is} X\left( \tfrac{(q^{1/2}Y)^{it}-(q^{1/2} q^{-2}Y)^{it}}{(1-q^{-2}) Y}\right)\nonumber\\
&=q^{1/2} e^{-\pi\gamma s} e^{\pi i\gamma st} X^{is} X\left(e^{-\pi\gamma t} - e^{3\pi\gamma t}\right) Y^{it} Y^{-1}\nonumber\\
&=q^{1/2} e^{-\pi\gamma s} \left( e^{-\pi\gamma t} - e^{3\pi\gamma t}\right) e^{\pi i\gamma st} X^{is + 1} Y^{it - 1}\nonumber\\
&=q^{1/2} e^{-\pi\gamma s}\left(e^{-\pi\gamma t} - 
e^{3\pi\gamma t}\right) e^{\pi i \gamma st} W(s{-}i,0)W(0,t{+}i)\nonumber\\
&=q^{1/2} e^{-\pi\gamma s} 
\left(e^{-\pi\gamma t}-e^{3\pi\gamma t}\right) e^{\pi i\gamma st} 
e^{-\pi i\gamma(s-i)(t+i)} W(s{-}i,t{+}i)\nonumber\\
&=\left(e^{-2\pi\gamma t}-e^{2\pi\gamma t}\right) W(s{-}i, t{+}i).\nonumber
\end{align*}
The formulas for the actions of the other generators $E, K_1$ and $K_2$ 
are derived by a similar formal reasoning. 
Replacing (\ref{A7}) by (\ref{A6}) 
and (\ref{A5}) we obtain
\begin{align*}
&\lambda E\ang W(s,t)=(e^{- 2\pi \gamma s}-e^{2\pi \gamma s})W(s{+}i,t{-}i),\\
&K_1\ang W(s,t)=e^{\pi\gamma s} W(s,t),~~K_2\ang W(s,t)=e^{\pi\gamma t} W(s,t).
\end{align*}

We now take the above formulas which have been
obtained by {\it formal} algebraic manipulations as the starting point 
for the {\it rigorous} definition of a left action of $\U_q(gl_2)$ on the 
$\ast$-algebra $\W$. That is, for $s,t\in\bbbc$ we define
\begin{align}\label{B13}
&E~~\ang W(s,t)=\lambda^{-1} (e^{-2\pi\gamma s}-e^{2\pi\gamma s}) W(s{+}i,t{-}i),\\
\label{B14}
&F~~\ang W(s,t)=\lambda^{-1} (e^{-2\pi\gamma t}-e^{2\pi\gamma t}) W(s{-}i, t{+}i),\\
\label{B15}
&K_1\ang W(s,t)= e^{\pi\gamma s} W(s,t),~K_2\ang W(s,t)= 
e^{\pi\gamma t} W(s,t).
\end{align}

\mn
{\bf Proposition 17.} {\it With definitions (\ref{B13})--(\ref{B15}), 
$\W$ is a left $\U_q(gl_2(\bbbr))$-module $\ast$-algebra.}

\mn
{\bf Proof.} Since the set of operators $W(s,t),s,t\in\bbbc$, is 
linearly independent as easily shown, the preceding definitions 
extend uniquely to well-defined linear mappings of $\W$ into 
itself. It is straightforward to check that the terms $K_1E-q^{1/2}EK_1, 
K_2E-q^{-1/2}EK_2, K_1F-q^{-1/2} FK_1, K_2F-q^{1/2}FK_2$ 
and $\lambda EF-\lambda FE-K^{2}+K^{-2}$ applied to an arbitrary basis 
element $W(s,t)$ of $\W$ vanish. Thus, formulas (\ref{B13})--(\ref{B15}) define indeed a left action of the algebra 
$U_q(gl_2)$ on $\W$.
That the left module $\W$ is a $\U_q(gl_2)$-module {\it algebra} 
means that (\ref{A8}) is satisfied. 
It suffices to check this condition for the 
generators $f=E,F,K_1,K_2$, $K_1^{-1}$, $K_2^{-1}$ and $z=W(s,t),z^\prime=
W(s^\prime,t^\prime),s,t,s^\prime, t^\prime\in\bbbc$. As a sample, we carry out this for the generator $f=E$. Using (\ref{B13}), 
(\ref{B15}) and (\ref{B7}) we compute 
\begin{align*}
&\lambda (E\ang W(s,t))(K\ang W(s^\prime, t^\prime))+\lambda 
(K^{-1}\triangleright W(s,t))(E\ang W(s^\prime, t^\prime))\\
&=(e^{-2\pi\gamma s}-e^{2\pi\gamma s}) e^{\pi\gamma (s^\prime - t^\prime)} W(s{+}i, t{-}i)W(s^\prime, t^\prime)\\
&\quad +e^{\pi\gamma (t-s})(e^{-2\pi\gamma s^\prime}-e^{2\pi\gamma s^\prime})W(s,t)
W(s^\prime+i t^\prime{-}i)\\
&=(e^{-2\pi\gamma s}- e^{2\pi\gamma s})e^{\pi\gamma(s^\prime-t^\prime)}
e^{\pi i\gamma(s^\prime(t-i)-(s+i)t^\prime)} W(s+s^\prime{+}i,t+t^\prime{-}i)\\
&\quad +e^{\pi\gamma (t-s})(e^{-2\pi\gamma s^\prime}-e^{2\pi\gamma s^\prime})
e^{\pi i\gamma((s^\prime+i)t-s(t^\prime-i))}W(s+s^\prime{+}i,t+t^\prime{-}i)\\
&=(e^{-2\pi\gamma(s+s^\prime)}-e^{2\pi\gamma (s+s^\prime)})e^{\pi 
i\gamma(s^\prime t-st^\prime)} W(s+s^\prime{+}i,t+t^\prime{-}i)\\
&=\lambda e^{\pi i \gamma (s^\prime t-st^\prime)}E\ang 
W(s{+}s^\prime,t{+}t^\prime)\\
&=\lambda E\ang (W(s,t)W(s^\prime,t^\prime)).
\end{align*}
This proves (\ref{A8}) in the case $f=E$.

Finally, it remains to check that (\ref{A15}) holds. Since $\W$ is a left 
$\U_q(gl_2)$-module algebra as just shown, it suffices to do this 
for the generators $f$ of $\U_q(gl_2)$. Again we restrict ourselves to the case 
$f=\lambda E, z=W(s,t)$. Since 
$S(\lambda E)^\ast=\overline{\lambda} (-q)E^\ast =-\lambda E$ by (\ref{A13}) and 
$W(s,t)^\ast=W(-\overline{s},-\overline{t})$, we have
\begin{align*}
\qquad\qquad&(\lambda E\ang W(s,t))^\ast =
(e^{-2\pi\gamma s}-e^{2\pi\gamma s}) W(s+i,t-i))^\ast\\
&=(e^{-2\pi\gamma \overline{s}}-e^{2\pi\gamma\overline{s}}) 
W(-\overline{s} +i,-\overline{t}-i)=-\lambda E\ang W(-\overline{s}, -\overline{t})\\
&=S(\lambda E)^\ast\ang W(s,t)^\ast.\hskip8cm \Box
\end{align*}

The $\ast$-algebra $\W$ consists of Hilbert space 
operators and formulas (\ref{B13})--(\ref{B15}) have been derived by 
using formal operator calculus. However, the content of 
Proposition 17 is purely algebraic: It is obvious that the complex vector 
space $\W$ with basis 
$W(s,t), s,t \in \bbbc$ is a $\ast$-algebra with 
multiplication and involution defined by (\ref{B5}) and (\ref{B6}). 
Proposition 17 says that $\W$ is a 
$\U_q(gl_2(\bbbr))$-module $\ast$-algebra with respect to the left action 
defined by (\ref{B13})--(\ref{B15}). 

Recall that $\cO(\bbbr^2_q)$ is a 
$\U_q(gl_2(\bbbr))$-module $\ast$-subalgebra of $\W$. 
By definition, the products $xw$ and $yw$ for $w \in \W$ are the 
operator products $Xw$ and $Yw$, respectively, in the 
Hilbert space $L^2(\bbbr)$. From (\ref{B4}) and (\ref{B8}) we obtain \begin{align}
\label{B1b}
xW(s,t) = e^{-\pi \gamma t} W(s{-}i,t),~~
yW(s,t) = e^{\pi \gamma s} W(s,t{-}i),\\
\label{B1bb}
W(s,t)x = e^{\pi \gamma t} W(s{-}i,t),~~
W(s,t)y = e^{-\pi \gamma s} W(s,t{-}i).
\end{align}

{\it 3.2 Covariant differential calculus on $\W$}
\sn

In this subsection we 
extend the differential calculus $\Gamma_-$ of $\cO(\bbbr^2_q)$ to $\W$. 
In order to do so, we use the approach given in 2.2 with $\Z=\W$ and  
write $\Gamma,d,\omega$ for $\Gamma_-,d_-,\omega_-$, respectively.

As in 2.2, we set $\omega=q^2x^2y^{-2}e_1+y^{-2} e_2$ and define
$$
dz=\omega z-z\omega,~z\in\W.
$$
Obviously, $\Gamma:=\W{\cdot} d\W{\cdot} \W$ is a first 
order differential calculus over $\W$ with differentiation $d$ such that the differentials 
$dx, dy$ form a free left $\W$-module basis of $\Gamma$. Because of 
this property, the partial derivatives $\partial_x(z)$ and $\partial_y(z)$ are well-defined by (\ref{F0}). In order to compute the latter 
for $z=W(s,t)$, we use the commutation rules 
$xW(s,t)=e^{-2\pi\gamma t}W(s,t)x$ and 
$yW(s,t)= e^{2\pi\gamma s}W(s,t)y$ (by (\ref{B1b}) and (\ref{B1bb})) 
and the expressions (\ref{F4}) for $dx$ and $dy$. Comparing coefficients 
in (\ref{F0}), we obtain for $s,t\in\bbbc$,
\begin{align*}
\partial_x(W(s,t))={\frac{1-e^{4\pi\gamma s}}{1-q^{-2}}} 
e^{\pi\gamma t} W(s{+}i,t),~~
\partial_y(W(s,t))={\frac{1-e^{4\pi\gamma t}}{1-q^{-2}}} 
e^{3\pi\gamma s} W(s,t{+}i.)
\end{align*}
 

\mn
{\bf 4. The ${\bf \U_q(gl_2(\bbbr))}$-module 
${\bf \ast}$-algebra ${\bf \A(\bbbr^{++}_q)}$}

\sn
{\it 4.1} In the preceding section we extended the action of the 
Hopf $\ast$-algebra $\U_q(gl_2(\bbbr))$ on $\cO(\bbbr^2_q)$ 
to the larger $\ast$-algebra $\W$ such that 
$\W$ is a left module $\ast$-algebra of $\U_q(gl_2(\bbbr))$. We 
now go one step further and make the $\ast$-algebra 
$\gA(\bbbr^2)$ into a left $\U_q(gl_2(\bbbr))$-module $\ast$-algebra. 
In order to do so we use the formulas 
(\ref{B13})--(\ref{B15}) in order to derive the corresponding formulas for the action of the generators $E,F,K_1,K_2$ on $Op(a)$.
Suppose that $a\in\gA(\bbbr^2)$. For the generator $E$ we obtain
\begin{equation*}
\begin{split}
&\lambda E\ang Op(a)=\gamma {\iint} \hat{a} (\alpha s,\beta t)(\lambda 
E\ang W (s,t)) ds dt\\
&\qquad\qquad =\gamma {\iint} \hat{a} (\alpha s,\beta t)
(e^{-2\pi\gamma s}{-}e^{2\pi\gamma s}) W(s{+}i,t{-}i) ds dt\\
&\qquad\qquad=\gamma{\iint}\hat{a}(\alpha(s{-}i),\beta(t{+}i))(e^{-2\pi
\gamma(s-i)}{-}e^{2\pi\gamma(s-i)})W(s,t) ds dt.
\end{split}
\end{equation*}
Let us explain the steps of this computation. The first equality is only a formal interchanging of integrals 
and left action, while the second follows from formula (\ref{B13}). 
The third equality is obtained by the formal replacements 
$s \rightarrow s{+}i$ and $t \rightarrow  t{-}i$. These substitutions are 
justified by a standard 
argument from complex analysis which has been used already in the proof of Lemma 11: 
The integral of the holomorphic 
operator-valued function 
$$ 
s\rightarrow\hat{a}(\alpha s,\beta t)
(e^{-2\pi\gamma s}{-}e^{2\pi\gamma s})W(s+i,t-i)
$$ 
along the boundary of the rectangle $-R, R,R{-}i,-R{-}i$ for fixed 
$t\in\bbbc$ and $R>0$ is zero. By Lemma 9,  the integrals from $R$ to $R-i$ and from $-R-i$ to $-R$ tend to zero as 
$R\rightarrow +\infty$. Arguing similarly for the variable $t$, 
the third equality is obtained. In order to complete this reasoning, we note that the function 
$$
\hat{a}(\alpha(s{-}i),\beta(t{+}i))(e^{-2\pi\gamma(s-i)}-
e^{2\pi\gamma(s-i)})
$$ 
is the Fourier transform of the function $a_E\in\gA(\bbbr^2)$ defined by
\begin{equation}
a_E(x_1,x_2):=e^{2\pi(\beta x_2-\alpha x_1)}(a(x_1+\beta i,x_2)-a(x_1-\beta i,x_2)).
\end{equation}
Thus, we have seen that $\lambda E\ang Op(a)=Op(a_E)$. Using
(\ref{B14})
and (\ref{B15}) instead of (\ref{B13}) a similar reasoning shows that 
$\lambda F\ang Op(a)=Op(a_F)$ and $K_j\ang Op(a)=Op(a_{K_j})$, 
where the symbol $a_F,a_{K_j}\in\gA(\bbbr^2)$ are given by
\begin{equation*}
\begin{split}
a_F(x_1,x_2)&=e^{2\pi(\alpha x_1-\beta x_2)}(a(x_1,x_2+\alpha i)-
a(x_1,x_2-\alpha i)),\\
a_{K_1}(x_1,x_2)&= a(x_1-{\tfrac{\beta}{2}} i, x_2),\\
a_{K_2}(x_1,x_2)&= a(x_1,x_2-{\tfrac{\alpha}{2}} i).
\end{split}
\end{equation*}

Summarizing, in terms of the symbol we have derived the following 
formulas for the actions of the generators $E,F,K_1,K_2$ of 
$\U_q(gl_2(\bbbr))$:
\begin{align}
\label{C6}
(E\ang a)(x_1,x_2)&=\lambda^{-1} e^{2\pi(\beta x_2-\alpha x_1)} 
(a(x_1{+}\beta i, x_2)-a(x_1{-}\beta i, x_2)),\\
\label{C7}
(F\ang a)(x_1,x_2)&=\lambda^{-1} e^{2\pi(\alpha x_1-\beta x_2)} 
(a(x_1,x_2{+}\alpha i)-a(x_1,x_2{-}\alpha i)),\\
\label{C8}
(K_1\ang a)(x_1,x_2)&= a(x_1{-}{\tfrac{\beta}{2}} i, x_2),~
(K_2\ang a)(x_1,x_2)= a(x_1, x_2 {-}{\tfrac{\alpha}{2}} i).
\end{align}
The derivation of these formulas is rigorous except for the 
justification of the interchanging of integrals and actions. This could 
be made rigorous by introducing appropriate locally convex topologies. We 
shall not proceed this way, because we shall use formulas (\ref{C6})--(\ref{C8}) only as definitions of the action of $U_q(gl_2(\bbbr))$ on $\gA(\bbbr^2)$ and prove the corresponding properties directly in 4.3.

Note that formulas (\ref{C6})--(\ref{C8}) and also formulas (\ref{prodx}) and (\ref{prody}) below are meaningful for larger classes of symbols rather than $\gA (\bbbr^2)$. For instance, for the function $a(x_1,x_2)=e^{2\pi i(\alpha sx_1+\beta tx_2)}$ (which is of course not in $\gA(\bbbr^2)$ ) we have $Op(a)=e^{2\pi i(\alpha s\Q+\beta t\cP)}=W(s,t)$. In this case formulas (\ref{C6})--(\ref{C8}) reduces to the equations (\ref{B13})--(\ref{B15}) derived in the preceding section. If we allow the symbols to be distributions, then we recover also formulas (\ref{A5})--(\ref{A7}).

In a similar manner the product of the operators $Op (a)$ with operators $X$ and $Y$ can be computed by using formulas (\ref{B1b}) and (\ref{B1bb}) (or (\ref{C4})). We then obtain $X Op(a)=Op({_xa})$, $Op(a) X=Op(a_x)$, $Y Op(a)=Op({_ya})$ and $Op(a) Y=Op(a_y)$, where the symbol ${_xa},a_x,{_ya}, a_y\in\gA(\bbbr^2)$ are given by
\begin{align*}
&{_xa}(x_1,x_2)=e^{2\pi\alpha x_1}a(x_1,x_2{+}{\tfrac{\alpha}{2}} i),~~
a_x(x_1,x_2)= e^{2\pi\alpha x_1}a(x_1,x_2{-}{\tfrac{\alpha}{2}} i),\\
&{_ya}(x_1,x_2)=e^{2\pi\beta x_2}(a(x_1{-}{\tfrac{\beta}{2}} i),x_2),~~
a_y(x_1,x_2)= e^{2\pi\beta x_2}a(x_1{+}{\tfrac{\beta}{2}} i,x_2).
\end{align*}

Let $\A(\bbbr^{++}_q)$ denote the direct sum of vector spaces $\cO(\bbbr^2_q)$ and $\gA(\bbbr^2)$. 

\mn
{\bf Lemma 18.} {\it There is a unique structure of a $\ast$-algebra on 
$\A(\bbbr^{++}_q)$ such that $\cO(\bbbr^2_q)$ and $\gA(\bbbr^2)$ are 
$\ast$-subalgebras of $\A(\bbbr^{++}_q)$ and the products of the 
generators $x, y$ of $\cO(\bbbr^2_q)$ and symbols $a \in \gA(\bbbr^2)$
are given by}
\begin{align}\label{prodx}
xa(x_1,x_2)
=e^{2\pi\alpha x_1}a(x_1,x_2{+}{\tfrac{\alpha}{2}} i),~
ax(x_1,x_2)
=e^{2\pi\alpha x_1}a(x_1,x_2{-}{\tfrac{\alpha}{2}} i),\\
\label{prody}
ya(x_1,x_2)
=e^{2\pi\beta x_2}a(x_1{-}{\tfrac{\beta}{2}} i,x_2),~
ay(x_1,x_2)
= e^{2\pi\beta x_2}a(x_1{+}{\tfrac{\beta}{2}} i,x_2).
\end{align}
{\bf Proof.} We first note that the maps $z \rightarrow \rho_{++}(z)$ 
(see (\ref{B1a})) and $a \rightarrow Op (a)$ (see (\ref{oppro})) are 
faithful $\ast$-representations of the $\ast$-algebras $\cO(\bbbr^2_q)$ and $\gA(\bbbr^2)$ on the domain $\gA(\bbbr)$ on the Hilbert space $L^2(\bbbr)$. Since $\rho_{++}(z)Op(a)=Op({_za})$ and $Op(a)\rho_{++}(z) = Op(a_z)$ for $z=x,y$ by the above definitions, it is clear that the sum  $\rho_{++}(\cO(\bbbr^2_q))+ Op(\gA(\bbbr^2))$ of the images of these $\ast$-algebras is a $\ast$-algebra of unbounded operators on the domain $\gA(\bbbr)$. It is easily seen that the only bounded operators in $\rho_{++}(\cO(\bbbr^2_q))$ are the multiples of the identity operator and that no operator in $Op(\gA(\bbbr^2))$ is a multiple of the identity. Thus,  $\rho_{++}(\cO(\bbbr^2_q)) \cap Op(\gA(\bbbr^2)) = \{0 \}$.  Hence the map $\J:(z,a) \rightarrow \rho_{++}(z) + Op(a)$ of $\A(\bbbr^{++}_q)$ to $\rho_{++}(\cO(\bbbr^2_q))+ Op(\gA(\bbbr^2))$ is bijective. The unique $\ast$-algebra structure on $\A(\bbbr^{++}_q)$ for which $\J$ is a $\ast$-homomorphism has obviously the desired properties. \hfill$\Box$ 

\sn
We shall show by Theorem 21 below that $\A(\bbbr^{++}_q)= 
\cO(\bbbr^2_q)+\gA(\bbbr^2)$ is even 
a left $\U(gl_2(\bbbr))$-module $\ast$-
algebra. We call this left $\U(gl_2(\bbbr))$-module $\ast$-
algebra $\A(\bbbr^{++}_q)$ the 
{\it $\ast$-algebra of functions on the 
quantum quarter plane}. Obviously, the $\ast$-subalgebra $\cO(\bbbr^2_q)$ is 
considered as the algebra generated by the two coordinate functions 
$x$ and $y$ of the quantum 
quarter plane. The elements of $\gA(\bbbr^2)$ can be 
interpreted as ``functions on the quantum quarter plane which go 
rapidly to zero at the boundary of the quantum quarter plane''. Note that 
$\gA(\bbbr^2)$ is a two-sided $\ast$-ideal 
of the $\ast$-algebra $\A(\bbbr^{++}_q)$. 

\mn
{\it 4.2}  In this subsection  we introduce two useful algebra homomorphisms 
in order to understand the algebraic content behind 
formulas (\ref{C6})--(\ref{C8}).  Let $\cB_q$ denote the complex unital 
algebra with generators 
$x_1,x_1^{-1},y_1,y^{-1}_1$, $x_2$, $x^{-1}_2$,$y_2,y^{-1}_2$ and 
defining relations
\begin{gather}\label{C9}
x_jy_j=q^{1/8} y_jx_j, x_jx_j^{-1}=x^{-1}_j x_j=1, y_jy^{-1}_j=y^{-1}_j y_j=1~\text{for}~j=1,2,\\
\label{C10}
x_1x_2=x_2x_1,y_1y_2=y_2y_1,x_1y_2=y_2x_1,x_2y_1=y_1x_2,
\end{gather}
where we set $q^{1/8}:=e^{\pi \gamma i / 4}$. The subalgebra 
$\cB_j$, $j=1,2$, generated by $x_j,x_j^{-1},y_j$, $y^{-1}_j$ is 
nothing but the localization of the algebra 
$\cO(\bbbc^2_{q^{1/8}})$ at the elements $x_j$ and $y_j$, and 
$\cB_q$ is just the tensor product of the algebras $\cB_1$ and $\cB_2$.

\mn
{\bf Lemma 19. } {\it There are injective algebra homomorphisms 
$\psi:\U_q(gl_2)\rightarrow\cB_q$ and $\psi:\cO(\bbbr^2_q)\rightarrow\cB_q$ 
such that}
\begin{align}\label{C11}
\psi(E)&=\lambda^{-1} x^2_2 x^{-2}_1 (y^{-4}_1-y^4_1),\\
\label{C12}
\psi(F)&=\lambda^{-1}x^2_1 x^{-2}_2 (y^{-4}_2-y^4_2),\\
\label{C13}
\psi(K_1)&=y_1^2,~~ \psi(K_2)=y_2^2,\\
\label{C14}
\psi(x)&=x^2_1 y^{-2}_2,~~ \psi(y)=x_2^2 y^2_1.
\end{align}
{\bf Proof.} In oder to prove the assertion for $\U_q(gl_2)$ it suffices 
to check that the operators $\psi (E), \psi(F), \psi(K_1)$ and $\psi(K_2)$ 
satisfy the defining relations of the algebra $\U_q(gl_2)$. Using the 
relations (\ref{C9})--(\ref{C10}) of the algebra $\cB_q$ we obtain
\begin{align*}
\lambda\psi(E)\psi(K_1)&=
x^2_2x^{-2}_1(y^{-4}_1-y^4_1)y^2_1 =x_2^2(x^{-2}_1 y^2_1)(y^{-4}_1-y^4_1)\\
&=x_2^2(q^{1/4})^{-2} y^2_1x^{-2}_1 (y^{-4}_1-y^4_1)
=q^{-1/2}y_1^2 x^2_2x^{-2}_1(y^{-4}_1-y^4_1)\\
&=q^{-1/2}\lambda\psi(K_1)\psi(E).
\end{align*}
The relations $\psi(E)\psi(K_2)=q^{1/2} \psi(K_2)\psi(E), \psi(F)\psi(K_1)=q^{1/2} \psi(K_1)\psi(F),$\\ 
$\psi(F)\psi(K_2)=q^{-1/2} \psi(K_2)\psi(F)$ 
and $\lambda (\psi(E)\psi(F)-\psi(F)\psi(E))=
\psi(K)^2-\psi(K)^{-2}$ are 
verified by similar computations. Obviously we have 
$\psi(x)\psi(y)=q\psi(y)\psi(x)$. Hence the above formulas 
define indeed algebra 
homomorphisms of $\U_q(gl_2)$ and $\cO(\bbbr^2_q)$ into $\cB_q$. 
Since the sets $\{E^kK_1^nK_2^mF^l;k,l \in \bbbn_0, n,m \in \bbbzz \}$, 
$\{x^k y^n; k,n \in \bbbn_0 \}$ and 
$\{x_1^k x_2^l y_1^n y_2^m; k,l,n,m \in \bbbn_0 \}$ are vector space bases of 
$\U_q(gl_2)$, $\cO(\bbbr^2_q)$ and $\cB_q$, respectively, it follows easily 
from formulas (\ref{C6})-(\ref{C14}) that the mappings 
$\psi:\U_q(gl_2)\rightarrow\cB_q$ and 
$\psi:\cO(\bbbr^2_q)\rightarrow\cB_q$ are injective.  \hfill$\Box$

\sn
Since $|q| = 1$, $\cB_q$ is 
a $\ast$-algebra with involution determined by $x_j^\ast := x_j$ and 
$y_j^\ast := y_j$, $j=1,2.$ The algebra homomorphism 
$\psi:\U^{tw}_q(gl_2(\bbbr))\rightarrow\cB_q$ does not preserve the 
involution. The next lemma shows that $\psi$ is similar to a 
$\ast$-homomorphism.

\mn
{\bf Lemma 20.} {\it For $z\in \U_q(gl_2(\bbbr))$ and 
$z \in \cO(\bbbr^2_q)$, define 
$$\varphi(z)=x_1x_2y^{-1}_1y_2 
\psi(z)(x_1x_2y^{-1}_1y_2)^{-1}.$$  
Then $\varphi:\U^{tw}_q(gl_2(\bbbr))\rightarrow\cB_q$ and 
$\varphi:\cO(\bbbr^2_q)\rightarrow\cB_q$ are injective $\ast$-homomorphisms of the corresponding $\ast$-algebras. In fact, we have}
\begin{align}\label{C14a}
\varphi(E^\prime) &= x^2_2 x^{-1}_1 (y^{-4}_1-y^4_1)x^{-1}_1
= x^2_2 x^{-2}_1 (q^{-1/2}y^{-4}_1-q^{1/2}y^4_1)\nonumber\\
&=x^2_2 (q^{1/2} y^{-4}_1-q^{-1/2}y^4_1)x^{-2}_1 ,\\
\label{C14b}
\varphi(F^\prime) &=x^2_1 x^{-1}_2 (y^{-4}_2-y^4_2)x^{-1}_2
=x^2_1 x^{-2}_2 (q^{-1/2}y^{-4}_2- q^{1/2}y^4_2),\nonumber\\
&=x^2_1(q^{1/2} y^{-4}_2- q^{-1/2}y^4_2)x^{-2}_2,\\
\label{C14bb}
\varphi(q^{-1/4}K_1)&= \psi(K_1) = y^2_1,~\varphi
(q^{- 1/4} K_2)= \psi(K_2) = y^2_2,\\
\label{C14c}
\varphi(x) &= \psi(x) = x_1^2 y_2^{-2},~~ \varphi(y) = \psi(y) = x_2^2 y^2_1. 
\end{align}
{\bf Proof.}  
Clearly, $\varphi:\U_q(gl_2(\bbbr)) \rightarrow\cB_q$ and 
$\varphi:\cO(\bbbr^2_q)\rightarrow\cB_q$ are injective homomorphisms, 
because $\psi$ are by Lemma 19. Therefore, it is sufficient to prove that 
$\varphi(z^\ast)=\varphi(z)^\ast$ for the four generators 
$z= E^\prime, F^\prime,q^{-1/4}K_1, 
q^{-1/4}K_2$ of $\U^{tw}_q(gl_2(\bbbr))$ and the two generators $z=x, y$ 
of $\cO(\bbbr^2_q)$. Since all these generators $z$ and their 
images $\varphi(z)$ are hermitean, 
it suffices to check formulas (\ref{C14a})--(\ref{C14c}). The 
latter formulas follow by straightforward computations from 
(\ref{C6})-(\ref{C14}) combined with the relations 
(\ref{C9})-(\ref{C10}) of the algebra $\cB_q$. As a sample we 
verify (\ref{C14a}) and compute
\begin{align*}
\varphi(\lambda E)&=x_1x_2y_1^{-1}y_2 \psi(\lambda E) 
(x_1x_2y_1^{-1}y_2)^{-1}\\
&=x_1x_2y_1^{-1}y_2x_2^2x_1^{-2}(y^{-4}_1-y^4_1)y_2^{-1}
y_1x^{-1}_2x^{-1}_1\\
&=x_1x_2q^{-1/4}x_1^{-2}y_1^{-1}q^{-1/4}x_2^2y_2(y^{-4}_1- y^{4}_1) 
y_2^{-1}y_1x_2^{-1}x_1^{-1}\\
&=q^{-1/2}x^2_2 x^{-1}_1(y^{-4}-y_1^4)x_1^{-1}
\end{align*} 
which gives the first formula of(\ref{C14a}). The second and third formulas of (\ref{C14a}) follow by applying 
once more the commutation rules (\ref{C9}) and (\ref{C10}). \hfill$\Box$ 

\sn
The $\ast$-homomorphisms $\varphi$ of $\U^{tw}_q(gl_2(\bbbr))$ and 
$\cO(\bbbr^2_q)$ are crucial in what follows. 

\mn
{\it 4.3} Let us return to the left action of the Hopf $\ast$-algebra 
$\U(gl_2(\bbbr))$ on $\gA(\bbbr^2)$ given by the formulas (\ref{C6})--(\ref{C8}). 
We define a $\ast$-representation $\rho_0$ of the 
$\ast$-algebra 
$\cB_q$ on the invariant dense domain $\gA(\bbbr^2)$ of the Hilbert 
space $L^2(\bbbr^2)$ by
\begin{equation}\label{C15}
\rho_0(x_1)=e^{\pi\alpha\Q_1},
\rho_0(y_1)=e^{{\tfrac{\pi}{2}}\beta\cP_1},
\rho_0(x_2)=e^{\pi\beta\Q_2},
\rho_0(y_2)=e^{{\tfrac{\pi}{2}}\alpha\cP_2},
\end{equation}
where $q^{1/4}=e^{\pi\gamma i/2}$ and as always 
$\alpha\beta=\gamma$. It is obvious that 
these operators satisfy the relations of the $\ast$-algebra 
$\cB_q$, so (\ref{C15}) defines indeed a $\ast$-representation of 
$\cB_q$. Inserting (\ref{C15}) into (\ref{C11})--(\ref{C13}) we see that 
equations (\ref{C6})--(\ref{C8}) can be expressed as
\begin{equation}\label{C16}
f\ang a=\rho_0 ( \psi (f))a,~a\in \gA(\bbbr^2),
\end{equation}
for the generators $f=E,F,K_1,K_2$ of $\U_q(gl_2)$. We now take 
this equation as a definition for arbitrary elements 
$f\in\U_q(gl_2)$. Since $\rho_0 \circ \psi$ is an algebra 
homomorphism, (\ref{C16}) gives a well-defined left action of the 
algebra $\U_q(gl_2)$ on $\gA(\bbbr^2)$. Recall from 2.1 that we have also 
a left action $\ang$ of $\U_q(gl_2)$ on $\cO(\bbbr^2_q)$. Hence the equation 
\begin{equation}
f\ang (z+a) := f\ang z + f \ang a,~ 
f\in\U_q(gl_2), z \in \cO(\bbbr^2_q),a\in \gA(\bbbr^2),
\end{equation}
defines a left action of $\U_q(gl_2)$ on the direct sum
$\A(\bbbr^{++}_q) = \cO(\bbbr^2_q)+ \gA(\bbbr^2)$. 
In terms of the $\ast$-representations $\rho_0$ formulas (\ref{prodx}) and 
(\ref{prody}) can be written as
\begin{align}\label{xaax}
xa =\rho_0\circ \psi (x)a = \rho_0(x_1^2y_2^{-2})a,~
ax=\rho_0(x_1^2y_2^2)a,\\
\label{yaay} 
ya =\rho_0 \circ \psi (y)a = \rho_0(x_2^2y_1^2)a,~ay=\rho_0(x_2^2y_1^{-2})a. 
\end{align}
The main result of this section is the following theorem.

\mn
{\bf Theorem 21.} {\it With the preceding definitions, the 
$\ast$-algebra $\A(\bbbr^{++}_q)$ of functions on the quantum 
quarter plane is a left $\U_q(gl_2(\bbbr))$-module $\ast$-algebra.}

\mn
{\bf Proof.} We already noticed that $\ang$ 
is a left action of the algebra $\U_q(gl_2)$ on $\A(\bbbr^{++}_q)$. 
It remains to show that 
conditions (\ref{A8}) and (\ref{A15}) are fulfilled for arbitrary elements $z,z^\prime \in \A(\bbbr^{++}_q)$ and $f \in \U_q(gl_2)$. 

We first prove that $\gA(\bbbr^2)$ is a $\U_q(gl_2)$-left module algebra. 
Since $\ang$ is a left action of $\U_q(gl_2)$, it suffices to prove 
(\ref{A8}) for the generators $f=\lambda E, \lambda F, K_1, K_2$. 
These verifications are 
lengthy but straigthforward. We restrict ourselves to the 
case $f=\lambda E$. Then we compute
\begin{align*}
&(\lambda E\ang a){\scriptstyle{\#}} (K\ang b)(x_1,x_2)+(K^{-1}\ang a){\scriptstyle{\#}} (\lambda E\ang b)(x_1,x_2)\\
=\  &4  \iiiint du_1 du_2 dv_1 dv_2 e^{4\pi i[(x_1-u_1)(x_2-v_2)-(x_1-v_1)(x_2-v_2)]}\\
& \Bigm\{  e^{2\pi(\beta u_2-\alpha u_1)}(a(u_1+\beta i,u_2)-a(u_1-\beta i,u_2))b(v_1-{\textstyle\tfrac{\beta}{2}}i, v_2+{\textstyle\tfrac{\alpha}{2}} i)\\
&+a(u_1+\tfrac{\beta}{2} i, u_2-\tfrac{\alpha}{2}i) 
e^{2\pi(\beta v_2-\alpha v_1)}(b(v_1+\beta i,v_2)-b(v_1-\beta i,v_2))\} \\
=\ &4 \iiiint du_1 du_2 dv_1 dv_2 a(u_1,u_2)b(v_1,v_2)\\
&\Bigm\{ e^{2\pi(\beta u_2-\alpha(u_1-\beta i))+
4\pi i[(x_1-u_1+\beta i)(x_2-v_2-i \alpha/2)-(x_1- v_1-i \beta/2) (x_2-u_2)]}\\
&-e^{2\pi(\beta u_2-\alpha(u_1+\beta i))+
4\pi i[(x_1-u_1+\beta i)(x_2-v_2+i \alpha/2)-
(x_1- v_1-i\beta/2) (x_2-u_2)]}\\
&+e^{2\pi(\beta v_2-\alpha(v_1-\beta i))+
4\pi i[(x_1-u_1+i\beta/2)(x_2-v_2)-(x_1-v_1+\beta i) 
(x_2-u_2-i\alpha/2)]}\\
&-e^{2\pi(\beta v_2-\alpha(v_1+\beta i))
+4\pi i[(x_1-u_1+i\beta/2)(x_2-v_2)-
(x_1- v_1-\beta i)(x_2-u_2-i\alpha/2)]}\Bigr\}\\
=\ &4 \iiiint du_1 du_2 dv_1 dv_2 a(u_1,u_2)b(v_1,v_2)\\
&\Bigm\{ -e^{2\pi(\beta x_2-\alpha x_1)
+4\pi i[(x_1-\beta i- u_1)(x_2-v_2)-(x_1- \beta i-v_1) (x_2-u_2)]}\\
&+e^{2\pi(\beta v_2-\alpha v_1)
+4\pi i[(x_1+\beta i-u_1)(x_2-v_2)-(x_1+\beta i-v_1)(x_2-u_2)]}\Bigr\}\\
=\ &(\lambda E\ang (a{\scriptstyle{\#}} b))(x_1,x_2).
\end{align*}
The first equality is obtained by inserting the formulas (\ref{C6}) and 
(\ref{C8}) for the actions of $E$ and $K$ and (\ref{C5}) for the product 
${\scriptstyle{\#}}$ of the algebra $\gA(\bbbr^2)$. The second equality follows by the substitution $u_1\rightarrow u_1+\beta i, v_1\rightarrow v_1-\tfrac{\beta}{2}i, v_2\rightarrow v_2+\tfrac{\alpha}{2} i$ of the first 
summand and similar replacements of the other three summands. As noted in 
the considerations preceding (\ref{C5a}), these substitutions are 
justified because of Lemma 9. 
Next let us consider the expressions in the four exponentials after the 
second equality sign. By regrouping these terms we see that the first and 
the fourth exponentials cancel, while the second and third ones can be
 reexpressed as the exponentials after the third equality sign. The 
fourth equality follows by applying once more formulas (\ref{C5}) and 
(\ref{C6}). By a similar reasoning condition (\ref{A8}) can be checked 
for the other generators $f= \lambda F, K_1, K_2$. Thus, 
$\gA(\bbbr^2)$ is a $\U_q(gl_2)$-left module algebra.

Recall from 2.1 that $\cO(\bbbr^2_q)$ is also a $\U_q(gl_2)$-left 
module algebra. Therefore, in order to prove that the sum 
$\A(\bbbr^{++}_q) = \cO(\bbbr^2_q)+ \gA(\bbbr^2)$ 
is a $\U_q(gl_2)$-left module algebra, it remains to show that 
\begin{align}\label{wacon}
f \ang (wa) &= (f_{(1)} \ang w)(f_{(2)} \ang a),\\
\label{awcon}
f \ang (aw) &= (f_{(1)} \ang a)(f_{(2)} \ang w)
\end{align}
for $f\in\U_q(gl_2), w \in \cO(\bbbr^2_q)$ and $a\in \gA(\bbbr^2)$. 
It is easily seen that equation (\ref{wacon}) holds 
for the product $fg$ and arbitrary $w$ and $a$ provided that 
(\ref{wacon}) holds for $f$ and arbitrary $w$ and $a$ and also for $g$ 
and arbitrary $w$ and $a$. Hence it suffices to check condition 
(\ref{wacon}) for elements $f$ from a set $M$ of generators of the algebra 
$\U_q(gl_2)$ and for arbitrary $w$ and $a$. Suppose in addition that 
$M$ is a vector space 
such that $\Delta (M) \subseteq \U_q(gl_2) \otimes M$. Let 
$ w, w^\prime \in \cO(\bbbr^2_q)$ such that (\ref{wacon}) holds for $w$
and all $f \in M$ and $a$ and also for $w^\prime$ and all 
$f \in M$ and $a$. We show that then (\ref{wacon}) holds for the product 
$ww^\prime$ and arbitrary $f \in M$ and $a$ by computing 
\begin{align*}
(f_{(1)} \ang (ww^\prime))(f_{(2)} \ang a)
{=} (f_{(1)} \ang w)(f_{(2)} \ang w^\prime)(f_{(3)} \ang a) 
{=} (f_{(1)} \ang w)(f_{(2)} \ang (w^\prime a))
{=} f \ang (w w^\prime a).
\end{align*}
Note that for the second equality we used that 
$\Delta (f)\subseteq \U_q(gl_2) \otimes M$ by
 assumption and so (\ref{wacon}) is valid for the elements in the second tensor factor of $\Delta (f)$. Applying the preceding 
with $M ={\rm Lin} \{E,F,K_1,K_2 \}$ we conclude that 
condition (\ref{wacon}) is fulfilled provided that it holds for 
$f= E,F,K_1,K_2$, $w=x,y$ and 
arbitrary $a \in \gA(\bbbr^2)$. Arguing in a similar manner with 
condition (\ref{awcon}) it follows that it is sufficient to verify 
(\ref{awcon}) for the generators $f= E,F,K_1,K_2$ and $w=x,y$.
As a sample, we prove equation (\ref{wacon}) for 
$f=E$ and $w=x$. Using the formulas (\ref{A3}), (\ref{C11}), 
(\ref{C13}) and (\ref{xaax}) we obtain 
\begin{align*}
(\lambda E \ang x)(K \ang a) + & (K^{-1} \ang x)(\lambda E \ang a)\\
&= \lambda y (K \ang a) + q^{1/2} x (\lambda E \ang a) \\
&= \lambda \rho_0(x_2^2y_1^2) \rho_0(y_1^2y_2^{-2})a
+ q^{1/2} \rho_0(x_1^2y_2^{-2}) \rho_0(x_1^{-2}x_2^2(y_1^{-4}-y_1^4))a\\
&=\rho_0(\lambda y_1^4x_2^2y_2^{-2} + 
q^{1/2}y_2^{-2}x_2^2(y_1^{-4}-y_1^4)) a\\
&=\rho_0((qy_1^{-4} -q^{-1} y_1^4)x_2^2y_2^{-2} )a\\
&=\rho_0 (x_1^{-2}x_2^2(y_1^{-4}-y_1^4)) \rho_0 (x_1^2y_2^{-2}) a\\
&=\lambda E \ang (xa).
\end{align*}
The other verifications are carried out in a similar manner. Thus we have 
shown that $\A(\bbbr^{++}_q)$ is a $\U_q(gl_2)$-left module algebra.

Finally, we turn to condition (\ref{A15}). Because 
$\A(\bbbr^{++}_q)$ is a left $\U_q(gl_2)$-module algebra, it is 
enough to prove (\ref{A15}) for the generators 
$f=\lambda E,\lambda F, K_1, K_2$. We verify (\ref{A15}) 
for $f=\lambda E$ and $z=a \in \gA(\bbbr^2)$. Since 
$S(\lambda E)^\ast=-\lambda E$, we obtain
\begin{equation*}
\begin{split}
\qquad\qquad(\lambda E\ang a)^\ast (x_1,x_2)&=e^{2\pi(\beta x_2{-}\alpha x_1)} 
(\overline{a(x_1{+}\beta i, x_2)}- \overline{a(x_1{-}\beta i, x_2)})\\
&=e^{2\pi(\beta x_2{-}\alpha x_1)}(\overline{a}(x_1{-}\beta i, x_2)- 
\overline{a}(x_1{+}\beta i, x_2))\\
&=(S(\lambda E)^\ast\ang\overline{a})(x_1,x_2).
\end{split}
\end{equation*}
By similar computations we check that 
condition (\ref{A15}) is satisfied for $f= \lambda F, K_1, K_2$. 
Hence it follows that (\ref{A15}) holds for $f \in \U_q(gl_2(\bbbr))$ and 
$z=a \in \gA(\bbbr^2)$. Since $\cO(\bbbr^2_q)$ is a left 
$\U_q(gl_2(\bbbr))$-module $\ast$-algebra, (\ref{A15}) holds also 
for $z=a \in\cO(\bbbr^2_q)$. Thus it remains to show that 
(\ref{A15}) is satisfied for products $z=wa$ and $z=aw$, 
where $a \in \gA(\bbbr^2)$ and $w \in\cO(\bbbr^2_q)$. 
We carry out this for $z=wa$ and compute
\begin{align*}
S(f)^\ast \ang (wa)^\ast 
&= S(f_{(2)}^\ast  \ang a^\ast)
S(f_{(1)}^\ast  \ang w^\ast)\\
&=(f_{(2)} \ang a)^\ast (f_{(1)}  \ang w)^\ast
=((f_{(1)}\ang w)(f_{(2})\ang a))^\ast 
=(f \ang (wa))^\ast.
\end{align*}
This completes the proof of Theorem 21. $\hfill \Box$

\sn
We close this subsection by collecting some additional 
useful relations. From the above formulas and the 
defining relations of the algebra $\cB_q$ 
we derive the following cross commutation relations for 
elements $\psi(E), \psi(F), \psi(K)$ 
and $\psi(x), \psi(y)$ in the algebra $\cB_q$: 
\begin{align*}
&\psi(E)\psi(x)- q^{1/2}\psi(x)\psi(E) = \psi(y)\psi(K),~ 
\psi(E)\psi(y) =q^{-1/2}\psi(y)\psi(E),\\
&\psi(F)\psi(x)=q^{1/2}\psi(x)\psi(F),~
\psi(F)\psi(y)-q^{-1/2}\psi(y)\psi(F) = \psi(x)\psi(K),\\
&\psi(K)\psi(x) =q^{-1/2} \psi(x)\psi(K),~ 
\psi(K)\psi(y) = q^{1/2}\psi(y)\psi(K).
\end{align*}
The images of the 
algebras $\U_q(sl_2(\bbbr))$ and $\cO(\bbbr^2_q)$ under the 
algebra homomorphism $\psi$ do not generate the whole algebra $\cB_q$. 
But they are large enough such that fourth powers of the generators 
$x_1, x_2, y_1^{-1}, y_2$ can be expressed as
\begin{align*}
y_1^{-4} &=q^{-1} + q^{1/2}\lambda \psi(EK^{-1})\psi(x)\psi(y)^{-1},~
y_2^4 = q - q^{-1/2}\lambda \psi(FK^{-1})\psi(y)\psi(x)^{-1},\\
x_1^4 &=\psi(x^4)(q-q^{-1/2}\lambda \psi(FK^{-1}) \psi(y)\psi(x)^{-1}),\\
x_2^4 &=\psi(y^4)(q^{-1}+q^{1/2}\lambda \psi(EK^{-1}\psi(x)\psi(y)^{-1}).
\end{align*}

\bn
{\it 4.4. Covariant differential calculus on the quantum quarter plane}

\sn 
In this subsection we extend the differential calculus 
$\Gamma=\Gamma_-$ of $\cO(\bbbr^2_q)$ to the larger algebra 
$\A(\bbbr^{++}_q)$. As in 3.2, we use the approach developed in 2.2, 
but now with the algebra $\Z=\A(\bbbr^{++}_q)$. We briefly repeat the 
construction from 2.2. 
Let $V$ be a two-dimensional vector space with basis $\{e_1,e_2\}$. The 
vector space $\Omega=\A(\bbbr^{++}_q)\otimes V$ becomes a 
$\A(\bbbr^{++}_q)$-bimodule with bimodule structure defined by (\ref{bimod}).
The differentiation $d$ is defined by the commutator with the element
$\omega=q^2x^2y^{-2}e_1+y^{-2} e_2$ of the $\A(\bbbr^{++}_q)$-bimodule 
$\Omega$, that is, 
$$
dz=\omega z-z\omega,~~ z\in\A(\bbbr^{++}_q).
$$
It is clear that $\Gamma:=\Z{\cdot} d\Z{\cdot}\Z$ is a 
first order differential calculus over $\Z=\A(\bbbr^{++}_q)$ 
with differentiation $d$ such that $\{ dx, dy\}$ is a free left 
$\Z$-module basis of $\Gamma$.

Let us compute the partial derivatives $\partial_x(a)$ and 
$\partial_y(a)$ for $a\in\gA(\bbbr^2)$.  Using formulas (\ref{prodx}), 
(\ref{prody}) and (\ref{C15}) we obtain from the definition of $d$ that
\begin{align}\label{C20}
da&=q^2 x^2y^{-2}e_1 a+y^{-2} e_2a-q^2ax^2y^{-2} e_1-a y^{-2} e_2\nonumber\\
&=q^2x^2y^{-2}\rho_0(1-y^8_1 y^8_2) ae_1 +y^{-2}\rho_0(1-y^8_1)ae_2.
\end{align}
On the other hand, by (\ref{F4}) 
and (\ref{F0}) we have
\begin{equation}\label{C21}
da=(q^2-1)(q^2x^3y^{-2}\partial_x(a) e_1+
xy^{-2}\partial_x(a) e_2+x^2y^{-1}\partial_y(a)e_1).
\end{equation}
Comparing the coefficient of $e_1$ and $e_2$ in (\ref{C20}) and 
(\ref{C21}) we derive
$$
\partial_x(a)=(1-q^{-2})^{-1} x^{-1}\rho_0(1-y^8_1)a,~~ \partial_y(a)=
(1-q^{-2})^{-1}y^{-1}\rho_0(y^8_1(1-y^8_2))a
$$
or equivalently
\begin{align*}
&\partial_x(a)=(1-q^{-2})^{-1} e^{-2\pi\alpha x_1}
(a(x_1,x_2{-}{\tfrac{\alpha}{2}}i)-a(x_1{-}2\beta i,
x_2{-}{\tfrac{\alpha}{2}}i)),\\
&\partial_y(a)=(1-q^{-2})^{-1} e^{-2\pi\beta x_2}
(a(x_1{-}{\tfrac{3}{2}}\beta i,x_2)-a(x_1{-}{\tfrac{3}{2}}
\beta i, x_2{-}2\alpha i)).
\end{align*}
In terms of the $\ast$-representation $\rho_0$ 
defined by (\ref{C15}) and the actions of  
generators $E,F,K_1,K_2,x,y$ given by (\ref{C6})--(\ref{C8}), 
these formulas can be written as 
\begin{align*}
&\partial_x(a)=(1-q^{-2})^{-1}\rho_0(x^{-2}_1 y^2_2(1-y^8_1))a=
q^{\frac{3}{2}}y^{-1}EK^3_1K_2\ang a,\\
&\partial_y(a)=(1-q^{-2})^{-1}\rho_0(x^{-2}_2y^6_1(1-y^8_2)) a=
q^{\frac{1}{2}} x^{-1} FK^3_1K_2\ang a
\end{align*}
for $a\in\gA(\bbbr^2)$. In 2.2 we have shown that the two latter 
expressions of $\partial_x(a)$ and $\partial_y(a)$ hold also 
for elements of $\cO(\bbbr^2_q)$. Therefore, we have proved that
$$
\partial_x(a)=q^{3/2} y^{-1} EK^3_1K_2\ang a,~\partial_y(a)=
q^{1/2} x^{-1} FK^3_1K_2\ang a
$$
for all elements $a$ in the $\ast$-algebra $\A(\bbbr^{++}_q)=\cO(\bbbr^2_q)+\gA(\bbbr^2)$.

\mn
{\bf 5. Covariant linear functionals on the quantum quarter plane}

\sn
{\it 5.1} In the first subsection we construct for any 
$k=(k_1,k_2)\in\bbbzz^2$ a $\U_q(gl_2(\bbbr))$-covariant linear 
functional $h_k$ on $\gA(\bbbr^2)$ and show that it defines a 
scalar product $\langle a,b \rangle_k$ on the $\ast$-algebra 
$\gA(\bbbr^2)$. This functional and the associated scalar product 
can be considered as $q$-analogs 
of the state given by Lebesgue measure and 
the $L^2$-scalar product on the classical quarter plane.
 
From the defining relations of the algebra $\U_q(gl_2)$ it is clear that 
there exists a unique character $\tau$ on $\U_q(gl_2)$ such that
$$
\chi (K_1)=\chi(K_2)=q^{1/2}~{\rm and}~\chi(E)=\chi(F)=0.
$$
Since 
the restriction of $\chi$ 
to the Hopf subalgebra $\U_q(sl_2)$ is the counit, 
any $\U_q(sl_2)$-covariant linear functional with respect to $\chi$ is 
$\U_q(sl_2)$-invariant.

\mn
{\bf Proposition 22.} {\it For $k=(k_1,k_2)\in\bbbzz^2$ and 
$a\in\gA(\bbbr^2)$ we define }
\begin{equation}\label{H1}
h_k(a)=\iint e^{2\pi(\alpha_k x_1+\beta_k x_2)} a(x_1,x_2)dx_1 dx_2,
\end{equation}
{\it where}
$$
\alpha_k:=\alpha+2\beta^{-1}k_1,~\beta_k:=\beta+2\alpha^{-1}k_2.
$$
(i) {\it The linear functional $h_k$ on the 
$\U_q(gl_2(\bbbr))$-module $\ast$-algebra $\gA(\bbbr^2)$ is 
covariant with respect to the character $\chi$. \\
(ii) {\it For $s,t\in\bbbr$ and a $\in\gA (\bbbr^{2})$,}
\begin{equation}\label{H2}
h_k(a(x_1{+}\beta s, x_2{+}\alpha t))=
e^{-2\pi\gamma (s+t)-4\pi (k_1s+k_2t)} h_k(a).
\end{equation}
(iii) {\it $h_k$  is continuous on the Frechet space $\gA(\bbbr^2)[\tau]$. 
More precisely, for any 
$\varepsilon=(\varepsilon_1,\varepsilon_2),\varepsilon_1>0,
\varepsilon_2>0$, we have}
$$
|h_k(a)|\le{\frac{1}{2\pi\sqrt{\varepsilon_1\varepsilon_2}}}\parallel 
S(e^{\varepsilon\Q})e^{2\pi(\alpha_k\Q_1+\beta_k\Q_2)}a\parallel,~ 
a\in\gA(\bbbr^2).
$$
}
{\bf Proof.} (i): It suffices to verify condition (\ref{G1}) for the 
generators $f=K_1,K_2,E,F$. That is, we have to show that 
$$
h_k(K_1\ang a)=h_k(K_2\ang a)=q^{1/2} h_k(a)~{\rm and}~ 
h_k(E\ang a)=h_k(F\ang a)=0
$$
for $a\in\gA(\bbbr_q^2)$. By formulas (\ref{C6})-(\ref{C8}), the 
latter conditions are equivalent to the relations
\begin{align*}
\iint e^{2\pi(\alpha_k x_1+\beta_k x_2)}
&a(x_1{-}i \beta / 2,x_2) dx_1 dx_2\\
&=\iint e^{2\pi(\alpha_k x_1+\beta_k x_2)}a(x_1,x_2{-}i\alpha / 2) dx_1 dx_2\\
&=e^{\pi\alpha\beta i}\iint e^{4\pi(\beta^{-1} k_1x_1+\alpha^{-1} k_2 x_2)}a(x_1,x_2)dx_1 dx_2,\\
\iint e^{2\pi(\alpha_k x_1 + (\beta_k x_2)} 
&a(x_1{+}\beta i, x_2) dx_1 dx_2\\ 
&=\iint e^{4\pi(\beta^{-1}k_1x_1 + (\alpha^{-1}k_2 +\beta) x_2)}a(x_1{-}\beta i, x_2) dx_1 dx_2,\\
\iint e^{4\pi((\alpha + \beta^{-1}k_1)x_1 + \alpha^{-1}k_2x_2)} 
&a(x_1, x_2{+}\alpha i) dx_1 dx_2 \\
&=\iint e^{4\pi((\alpha + \beta^{-1}k_1)x_1 + \alpha^{-1}k_2x_2)} a(x_1, x_2{-}\alpha i) dx_1 dx_2.
\end{align*}
These identities follow by the formal replacements 
$(x_1,x_2)\rightarrow (x_1{+} i \beta/2,x_2)$, 
$(x_1,x_2)\rightarrow (x_1,x_2{-}i\alpha/2)$ and 
$(x_1,x_2)\rightarrow (x_1,x_2{-}2\alpha i)$,
respectively. Similarly as above, these substitutions are justified by 
integrating in the complex plane and using the asymptotic estimate 
of Lemma 9.\\
(ii): The formula follows by the substitution 
$(x_1,x_2)\rightarrow(x_1{-}\beta s, x_2{-}\alpha t)$.\\
(iii) follows from (\ref{H1}) 
and the Cauchy-Schwarz inequality.  \hfill $\Box$

\mn
Let $\langle\cdot,\cdot\rangle_k$ be the sesquilinear form defined by 
means of the functional $h_k$ on the $\ast$-algebra $\gA(\bbbr^2)$, that is,
\begin{equation}\label{H3c}
\langle a,b \rangle_k = h_k (b^\ast {\scriptstyle{\#}} a),~ a, b \in \gA (\bbbr^2).
\end{equation}
Recall that $(\cdot,\cdot)$ denotes the 
scalar product of the Hilbert space $L^2(\bbbr^2)$. For $k\in\bbbzz^2$, we abbreviate
\begin{equation}\label{H3b}
T_k:= e^{\pi(\alpha_k\Q_1-{\frac{\beta_k}{2}} \cP_1)}\otimes e^{ \pi(\beta_k\Q_2+{\frac{\alpha_k}{2}}\cP_2) }\equiv
e^{\pi\alpha_k\Q_1}e^{-{\tfrac{\pi}{2}}\beta_k\cP_1}\otimes 
e^{\pi\beta_k\Q_2}e^{{\tfrac{\pi}{2}}\alpha_k\cP_2}.
\end{equation}

\mn
{\bf Proposition 23.} {\it The sesquilinear form $\langle\cdot,\cdot\rangle_k$ is a 
scalar product and $T_k$ is an isometric linear isomorphism of the unitary space $(\gA(\bbbr^2),\langle\cdot,\cdot\rangle_k)$ on the unitary space 
$(\gA(\bbbr^2), (\cdot,\cdot))$. For $a,b\in\gA(\bbbr^2)$, we have}
\begin{equation}\label{H3a}
\langle a,b \rangle_k=\iint e^{2\pi(\alpha_kx_1{+}\beta_kx_2)} a(x_1{+}i\beta_k/4,
x_2{-}i\alpha_k/4)\overline{b}(x_1{-}i\beta_k/4, x_2{+}i\alpha_k/4) dx_1 dx_2.
\end{equation}
{\bf Proof.} For $a,b\in\gA(\bbbr^2)$, we compute
\begin{align}
\langle a,b\rangle_k&=\iint e^{2\pi(\alpha_{k}x_1+\beta_kx_2)}
(b^\ast{\scriptstyle{\#}} a)(x_1,x_2)dx_1 dx_2\nonumber\\
&=\iint \left( e^{2 \pi(\alpha_k \Q_1 + 
\beta_k \Q_2)} b \right)^\ast {\scriptstyle{\#}} \left( e^{\pi (\alpha_k \cP_2 -\beta_k\cP_1)} a  \right)(x_1, x_2) dx_1 dx_2\nonumber\\
&=\iint\left( e^{2\pi(\alpha_k\Q_1+\beta_k\Q_2)}b\right)^\ast (x_1,x_2){\cdot}\left(e^{\pi(\alpha_k\cP_2-\beta_k\cP_1)}a\right)(x_1,x_2) dx_1 dx_2\nonumber\\
&=\left(e^{\pi(\alpha_k \cP_2 - \beta_k \cP_1)} a, 
e^{2 \pi(\alpha_k\Q_1 + \beta_k\Q_2)}b\right)\nonumber\\
\label{H4}
&=(T^2_k a,b)=(T_k a, T_k b).
\end{align}
Here the first equality combines the definitions (\ref{H1}) and (\ref{H3c}) of $h_k$ and $\langle\cdot,\cdot\rangle_k$, respectively. The second 
equality follows from (\ref{F5}) and (\ref{F7}), 
while the third one follows from formula (\ref{scalar}). The fourth 
equality is just the definition of the scalar product ($\cdot ,\cdot$.), and 
the fifth and the sixth are easily derived from (\ref{C111}). Since the operator $T_k$ is a bijective 
linear mapping of $\gA(\bbbr^2)$, we conclude from 
formula (\ref{H4}) that$\langle\cdot,\cdot\rangle_k$ is indeed 
a scalar product on the vector space $\gA(\bbbr^2)$.
The expression in (\ref{H3a}) is obtained from (\ref{H4}) 
by inserting the actions of the operator $T_k$. Further, 
it follows from (\ref{H4}) that $T_k$ is an isometric linear 
isomorphism of the unitary space $\gA_k:=(\gA(\bbbr^2),\langle\cdot,
\cdot\rangle_k)$ onto the unitary space $\gA(\cdot,\cdot):=
(\gA(\bbbr^2),(\cdot,\cdot))$. \hfill$\Box$

\mn
{\it 5.2} In this subsection we investigate the left actions of 
the algebras $\U_q(gl_2)$ 
and $\cO(\bbbr^2_q)$ on the unitary space 
$\gA_k:=(\gA(\bbbr^2),\langle\cdot,\cdot \rangle_k)$. Among others, we shall 
transform these actions to the domain $\gA(\bbbr^2)$ in the Hilbert 
space $L^2(\bbbr^2)$ by means of the unitary operator $T_k$.

Recall from 4.2 that the map $\rho_0\circ\psi$ defines left actions
of the algebras $\U_q(gl_2)$ 
and $\cO(\bbbr^2_q)$ on $\gA(\bbbr^2)$. For the generators 
$E,F,K_1, K_2$ of $\U_q(gl_2)$ this action has been also given 
by formulas (\ref{C6})--(\ref{C8}), see also (\ref{C16}). 
For the algebra $\cO(\bbbr^2_q)$ the action 
$\rho_0\circ\psi$ is just the left multiplication in the larger algebra 
$\A(\bbbr^{++}_q)$, see Lemma 18 and formulas (\ref{xaax}) and 
(\ref{yaay}). Let $\psi_k$ denote the action $\rho_0\circ\psi$ of
$\U_q(gl_2)$ and $\cO(\bbbr^2_q)$
considered as representation on the unitary space $\gA_k$. 
Since $T_k$ is a unitary transformation of 
$\gA_k\equiv(\gA(\bbbr^2), \langle\cdot,\cdot \rangle_k)$ on 
$(\gA(\bbbr^2),(\cdot,\cdot))$ by Proposition 23, $\psi_k$ is unitarily 
equivalent to the representation
\begin{equation}\label{I4}
\Psi_k(\cdot):=T_k\psi_k(\cdot)T^{-1}_k
\end{equation}
on the domain $\gA(\bbbr^2)$ in the Hilbert space $L^2(\bbbr^2)$. 
Further, the compositions 
$\Phi :=\rho_0\circ\varphi$ of the $\ast$-homomorphisms 
$\varphi$ (defined in Lemma 20) and the $\ast$-representation 
$\rho_0$ of $\cB_q$ (defined by (\ref{C15})) are also $\ast$-representations 
of the $\ast$-algebras $\U^{tw}_q(gl_2(\bbbr))$ and  $\cO(\bbbr^2_q)$, 
respectively, on the domain $\gA(\bbbr^2)$ in $L^2(\bbbr^2)$.

Let $T$ denote the operator 
$T_k$ defined by (\ref{H3b}) for $k=(0,0)$, that is,
\begin{equation}\label{I1}
T=e^{\pi\alpha\Q_1} 
e^{-{\tfrac{\pi}{2}}\beta\cP_1}\otimes e^{\pi\beta\Q_2} 
e^{ {\tfrac{\pi}{2}} \alpha\cP_2}
=e^{\pi\alpha\Q_1 -{\tfrac{\pi}{2}}\beta\cP_1}\otimes 
e^{\pi\beta\Q_2 +{\tfrac{\pi}{2}} \alpha\cP_2}.
\end{equation}
Using the operator
\begin{align}\label{I2}
\begin{split}
C_k&:=e^{2\pi k_1\beta^{-1}\Q_1} e^{-\pi k_2\alpha^{-1}\cP_1}
\otimes e^{2\pi k_2\alpha^{-1}\Q_2} e^{\pi k_1\beta^{-1}\cP_2}\\
&~=e^{\pi(2k_1\beta^{-1}\Q_1-k_2\alpha^{-1}\cP_1)}\otimes e^{\pi(2k_2\alpha^{-1}\Q_2+k_1\beta^{-1}\cP_2)}
\end{split}
\end{align}
acting on the Hilbert space $L^2(\bbbr^2)$, we can write the operator 
$T_k$, as
\begin{equation}\label{I3}
T_k=i^{k_1-k_2} C_k T~.
\end{equation}
Comparing formulas (\ref{I1}) and (\ref{C15}) we see that 
$T=\rho_0(x_1x_2y^{-1}_1y_2)$. Therefore, by Lemma 13, we get
\begin{align*}
\Psi_0(z)=T\psi_0(z)T^{-1}_0=\rho_0((x_1x_2y^{-1}_1y_2)
\psi(z)(x_1x_2y^{-1}_1y_2)^{-1})
=\rho_0\circ\varphi (z)=\Phi(z)
\end{align*}
for $z\in \U_q(gl_2)$ and $z\in\cO(\bbbr^2_q)$. That is, we have 
$\Psi_0 = \Phi$ for both $\ast$-algebras $\U^{tw}_q(gl_2(\bbbr))$ 
and $\cO(\bbbr^2_q)$.
 
Next we relate the representations $\Psi_k$ and $\Psi_0 = \Phi$. 
Using the definitions of the operator $C_k$ and 
$\Phi (f)=\rho_0\circ\varphi(f), f=E^\prime, F^\prime, K_1,K_2, x, y,$ 
we compute
\begin{align}\label{I5}
C_k\Phi(f)C^{-1}_k&=(-1)^{k_1+k_2}\Phi(f)~{\rm for}~ 
f=E^\prime, F^\prime, K,\\
\label{I6}
C_k\Phi(K_j)C^{-1}_k&=(-1)^{k_j}\Phi(K_j)~{\rm for}~j=1,2,\\
\label{I7}
C_k\Phi (x)C^{-1}_k&=\Phi(x),~ C_k\Phi(y) C^{-1}_k =\Phi(y)
\end{align}
Because of the formulas (\ref{I3}) and (\ref{I4}) we therefore have 
\begin{align}\label{I8}
\Psi_k(f)&=(-1)^{k_1+k_2}\Phi(f)~{\rm for}~ f=E^\prime, F^\prime, K,\\
\label{I9}
\Psi_k(K_j)&=(-1)^{k_j}\Phi(K_j)~{\rm for}~j=1,2,\\
\label{I10}
\Psi_k(x)&=\Phi(x),~\Psi_k(y)=\Phi(y).
\end{align}
In order to complete the picture we collect the formulas for the 
operators $\Phi(f)$, where 
$f=E^\prime, F^\prime, K_1,K_2, x, y$. Recall from (\ref{f}) and (\ref{lrform}) 
that $L_\alpha$ denotes the operator 
$L_\alpha =\overline{f_\alpha} (\cP) e^{-2\pi\alpha Q}$, 
where $f_\alpha(x) = - 2\sinh\pi\beta(2x+\alpha i)$. 
Combining (\ref{C14a})-(\ref{C14c}) and (\ref{C15}) we obtain
\begin{align}\label{I11}
\Phi(E^\prime)&= L_\alpha \otimes e^{2\pi\beta\Q_2},
\Phi(F^\prime)=e^{2\pi\alpha\Q_1}\otimes L_{\beta},\\
\label{I12}
\Phi(q^{-1/4}K_1)&= e^{\pi\beta\cP_1}\otimes I,~
\Phi(q^{-1/4}K_2)= I\otimes e^{\pi\alpha\cP_2},~
\Phi(K)=e^{\pi\beta\cP_1}\otimes e^{-\pi\alpha\cP_2},\\
\label{I13}
\Phi(x)&=e^{2\pi\alpha\Q_1}\otimes e^{-\pi\alpha\cP_2}, \Phi(y)=
e^{\pi\beta\cP_1}\otimes e^{2\pi\beta\Q_2}.
\end{align}

\mn 
Let us briefly discuss the outcome of these considerations. Since 
the functional $h_k$ on the left $\U_q(gl_2(\bbbr))$-module 
$\ast$-algebra $\gA(\bbbr^2)$ is covariant with 
respect to $\chi$, it follows from 
Lemma 2,(i)$\rightarrow$(ii), and the definition of the 
involution of $\U^{tw}_q(gl_2(\bbbr))$ that 
$\psi_k=\rho\circ\psi$ is a $\ast$-representation of the 
$\ast$-algebra $\U^{tw}_q(gl_2(\bbbr))$ on the unitary space 
$\gA_k=(\gA(\bbbr^2),\langle \cdot,\cdot \rangle_k)$. The $\ast$-representation 
$\psi_k$ is unitarily equivalent to the $\ast$-representation 
$\Psi_k$ of $\U^{tw}_q(gl_2(\bbbr))$ on the domain $\gA(\bbbr^2)$ in 
the Hilbert space $L^2(\bbbr^2)$. The actions of the operators 
$\Psi_k(f)$ for the generators $f=E^\prime, F^\prime, K_1,K_2$ are 
explicitly given by the formulas (\ref{I8})--(\ref{I9}) and 
(\ref{I11})--(\ref{I12}). Note that the dependence of the 
operators $\Psi_k(f)$ on 
$k\in\bbbzz^2$ appears only in the signs in (\ref{I8})-(\ref{I9}). 
In particular, if $k_1$ and $k_2$ are both even, then the 
$\ast$-representation $\Psi_k$ of $\U^{tw}_q(gl_2(\bbbr))$ on $\gA_k$ is 
unitarily equivalent to the fixed $\ast$-representation $\Phi$ on the 
domain $\gA(\bbbr^2)$ in the Hilbert space $L^2(\bbbr^2)$. From 
(\ref{I10}) and (\ref{I13}) we see that $\Psi_k$ is a 
$\ast$-representation of 
$\cO(\bbbr^2_q)$ on the unitary space $\gA_k$. For any $k\in\bbbzz^2$, the 
$\ast$-representation $\Psi_k$ of $\cO(\bbbr^2_q)$ on $\gA_k$ is unitarily 
equivalent to the 
$\ast$-representation $\Phi$ of $\cO(\bbbr^2_q)$ on the domain 
$\gA(\bbbr^2)$ in $L^2(\bbbr^2)$.

\mn
Remark 3. The preceding derivation shows the reason for the non-uniqueness of covariant functionals on the left 
$\U_q(gl_2(\bbbr))$-module $\ast$-algebra 
$\gA(\bbbr^2)$ from the technical side: For even numbers $k_1$ and $k_2$ 
the unbounded positive self-adjoint operator $C_k$ 
commutes with all representation operator 
$\Phi(z)$, $z\in\U^{tw}_q(gl_2(\bbbr))$, so that the 
{\it unbounded} commutant of $\Phi(\U_q(gl_2))$ is non-trivial. However, 
it can be shown that the $\ast$-representation $\Phi$ of 
$\U^{tw}_q(gl_2(\bbbr))$ on $L^2(\bbbr^2)$ is irreducibel. Hence the 
bounded commutant (more precisely, the strong bounded commutant, see [S1]) is trivial. The existence of examples of that kind is 
a well-known phenomena for unbounded operator algebras.

\mn

By the preceding we have expressed the actions $\psi_k$ of the algebras 
$\U_q(gl_2)$ and $\cO(\bbbr^2_q)$ on the unitary space 
$\gA_k=(\gA(\bbbr^2),\langle\cdot,\cdot \rangle_k)$ by means of the 
$\ast$-representations $\Phi $ on the domain $\gA(\bbbr^2)$ in 
$L^2(\bbbr^2)$. Since the representation 
$\Phi = T \psi_0 T^{-1}$ is obtained from $\psi_0$ by 
the unitary operator $T$, it is natural to transform also 
the structure of the $\ast$-algebra $\A(\bbbr^{++}_q)$ 
and the covariant functional $h:= h_0$ 
under the bijective linear mapping $T$ of $\gA(\bbbr^2)$. 
That is, for $f, g \in \A(\bbbr^{++}_q)$ and $a \in \gA(\bbbr^2)$
we define
\begin{equation}\label{timesprod}
f \natural g = T(T^{-1} f \cdot T^{-1}g),~~ 
f^\star = T(T^{-1}(f)^\ast)~~{\rm and}~~
{\tilde h}(a) = h(T^{-1}a).~~ 
\end{equation}
Since $\A(\bbbr^{++}_q)$ is a $\ast$-algebra with product ${\cdot}$ and 
involution $f \rightarrow f^\ast$, the vector space $\A(\bbbr^{++}_q)$ is 
a $\ast$-algebra, denoted ${\tilde \A}(\bbbr^{++}_q)$, with product $\natural$ 
and involution $f \rightarrow f^\star$. Further, since $\A(\bbbr^{++}_q)$ 
is a left $\U_q(gl_2(\bbbr))$-module $\ast$-algebra (by Theorem 21) 
and $h$ is covariant with respect to the left action $\psi_0$ 
(by Proposition 22), it is clear that ${\tilde \A}(\bbbr^{++}_q)$ 
is a left $\U_q(gl_2(\bbbr))$-module $\ast$-algebra and the 
linear functional ${\tilde h}$ is covariant with respect to the 
left action $\Phi = T \psi_0 T^{-1}$. 

Let us make the transformed structures more explicit. 
Suppose that $f=z+a$ and $g=w+b$, where $z,w \in 
\cO(\bbbr^2_q)$ and $a,b \in \gA(\bbbr^2)$. From (\ref{prodx}), (\ref{prody}) 
and (\ref{I1}) it follows that $xTc=Txc$ and $yTc=Tyc$ and so $vTc=Tvc$ 
for all $v \in \cO(\bbbr^2_q)$ and $c \in \gA(\bbbr^2)$ . Hence we get 
$$
(z+a)\natural (w+b) = zw+zb+aw+ a\natural b ~~{\rm and}~~ 
(z+a)^\star = z^\ast + a^\star .
$$
That is, product and involution of the 
$\ast$-algebra $\cO(\bbbr^2_q)$ remain unchanged and also the products of 
elements of $\cO(\bbbr^2_q)$ and $\gA(\bbbr^2)$. It remains 
to describe the transformed product and involution of the $\ast$-subalgebra 
$\gA(\bbbr^2)$. From the definition 
(\ref{I1}) of the operator $T$ and the formulas (\ref{F5})--(\ref{F12}) 
we obtain  
\begin{align}\label{atimesb}
a \natural b &= a {\scriptstyle{\#}} (e^{-{\tfrac{\pi}{2}}\beta\cP_1} e^{-\pi\alpha\Q_1}
\otimes e^{ {\tfrac{\pi}{2}} \alpha\cP_2}e^{-\pi\beta\Q_2} b)
= (e^{{\tfrac{\pi}{2}}\beta\cP_1} e^{-\pi\alpha\Q_1}
\otimes e^{- {\tfrac{\pi}{2}} \alpha\cP_2}e^{-\pi\beta\Q_2}a) {\scriptstyle{\#}} b \nonumber\\ 
&= (e^{{\tfrac{\pi}{4}}\beta\cP_1} e^{- {\tfrac{\pi}{2}} \alpha\Q_1}
\otimes e^{-{\tfrac{\pi}{4}}\alpha\cP_2} 
e^{- {\tfrac{\pi}{2}} \beta\Q_2}) a {\scriptstyle{\#}} 
(e^{-{\tfrac{\pi}{4}}\beta\cP_1} e^{- {\tfrac{\pi}{2}} \alpha\Q_1}
\otimes e^{{\tfrac{\pi}{4}}\alpha\cP_2} e^{- {\tfrac{\pi}{2}} \beta\Q_2}) b,\\ 
\label{star}
&a^\star (x_1,x_2) = 
(e^{\pi (\beta \cP_1 - \alpha \cP_2)} a)(x_1,x_2)
=\bar{a} (x_1-{\tfrac{\beta}{2}} i,x_2+{\tfrac{\alpha}{2}} i).
\end{align}
The vector space $\gA(\bbbr^2)$ equipped with this transformed 
$\ast$-algebra 
structure is denoted by ${\tilde \gA}(\bbbr^2)$. 
Inserting the definition (\ref{H1}) of the functional $h=h_0$ and 
substituting 
$(x_1,x_2) \rightarrow (x_1+{\tfrac{\beta}{4}}i, x_2 - {\tfrac{\alpha}{4}} i)$ 
we get
\begin{equation}\label{htilde}
{\tilde h}(a) = \iint  e^{\pi(\alpha x_1 +\beta x_2)}
a(x_1,x_2) dx_1dx_2 ,~~ a \in \gA(\bbbr^2).
\end{equation}
Summarizing, we have transformed all structures obtained so far of the left 
$\U_q(gl_2(\bbbr))$-module $\ast$-algebra $\A(\bbbr^{++}_q)$ of 
``functions on the quantum quarter plane'' under the the mapping $T$. 
The main advantage of this new picture is 
that the scalar product 
$$
(a,b)_{\tilde h} := {\tilde h}(b^\star \natural a),~a,b \in \gA(\bbbr^2),
$$
derived from the transformed covariant functional 
${\tilde h}$ is just the $L^2$-scalar product $(\cdot,\cdot)$ on 
$\bbbr^2$. (This follows at once from the construction. It can be also 
verified directly by using formulas (\ref{atimesb}) and (\ref{htilde}).) 
The latter fact will be crucial when we are looking for 
self-adjoint extensions of the representation operators 
$\Phi(E^\prime)$ and $\Phi(F^\prime)$ in a larger Hilbert space.

\mn
{\it 5.3} In this last subsection of Section 5 we prove two 
uniqueness result for covariant linear functionals.
They say that under some technical assumptions equation 
(\ref{H2}) essentially characterizes the covariant functional $h_k$.

\mn
{\bf Proposition 24.} {\it Let $k=(k_1,k_2)\in\bbbzz^2$ and let $h$ be a 
faithful positive linear functional on the $\ast$-algebra 
$\gA_{pex}(\bbbr^2)$ which satisfies relation (\ref{H2}) for all 
$s,t\in\bbbr$ and $a\in \gA_{pex}(\bbbr^2)$. Suppose that $h$ 
is continuous relative to the topology $\tau$. Then there exists a 
positive number $\nu$ such that $h(a)=\nu h_k(a)$ for all 
$a\in\gA_{pex}(\bbbr^2)$.}

\mn
{\bf Proof.} The proof mimics some of the preceding considerations in 
reversed order. Since $h$ is a faithful positive linear functional 
and $T_k$ is a bijective linear mapping of $\gA_{pex}(\bbbr^2)$, the 
equation
\begin{equation}\label{D1}
(a,b)_h=h(T_k^{-1}b)^\ast{\scriptstyle{\#}}(T_k^{-1} a)),~ a,b\in\gA_{pex}(\bbbr^2)~,
\end{equation}
defines a scalar product on $\gA_{pex}(\bbbr^2)$. 
(In the case when $h=h_k$ this scalar product is just the $L^2$-scalar 
product.) Consider the following four one-parameter 
groups of operators acting on the unitary space 
$(\gA_{pex}(\bbbr^2), (\cdot,\cdot)_h)$:
\begin{align*}
&(U_1(t)f)(x_1,x_2)=e^{2\pi it\alpha x_1} f(x_1,x_2),~ 
(U_2(t)f)(x_1,x_2)=e^{2\pi it \beta x_2} f(x_1,x_2),\\
&(V_1(t)f)(x_1,x_2)=f(x_1{+}\beta t,x_2),~ 
(V_2(t)f)(x_1,x_2)=f(x_1,x_2{+}\alpha t)~.
\end{align*}
It is straightforward to verify the commutation relations
\begin{align*}
&T^{-1}_k U_1(t)=e^{{\frac{\pi}{2}} t\gamma+\pi tk_2} U_1(t)T^{-1}_k,~ T^{-1}_k U_2(t)=e^{-{\frac{\pi}{2}} t\gamma-\pi tk_1} U_2(t)T^{-1}_k,\\
&T^{-1}_k V_1(t)=e^{\pi t\gamma+2\pi tk_1} V_1(t)T^{-1}_k,~ T^{-1}_k V_2(t)=e^{\pi t\gamma+2\pi tk_2} V_2(t)T^{-1}_k
\end{align*}
for arbitrary $t\in\bbbr$. From these commutation rules and the assumption 
(\ref{H2}) we derive that the operators $U_1(t), U_2(t), V_1(t), V_2(t)$ 
preserve the scalar product $(\cdot,\cdot)_h$ for real $t$. Let us carry 
out this computation for $U_1(t)$. Using the definition of the scalar 
product $(\cdot,\cdot)_h$, the above relation for 
$T_k^{-1} U_1(t)$, formula (\ref{H3b}) and finally the 
assumption (\ref{H2}), we obtain
\begin{align}\label{D2}
(U_1(t)a,U_1(t)b)_h&=h((T^{-1}_kU_1(t)b)^\ast{\scriptstyle{\#}} (T^{-1}U(t)a))\nonumber\\
&=e^{\pi t\gamma+2\pi tk_2}h((e^{2\pi it\alpha \Q_1}T^{-1}_kb)^\ast{\scriptstyle{\#}}(e^{2\pi it\alpha \Q_1}T^{-1}_k a))\nonumber\\
&=e^{\pi t\gamma+2\pi tk_2}h((e^{-2\pi it\alpha \Q_1}(T^{-1}_k b)^\ast){\scriptstyle{\#}} (e^{2\pi it\alpha \Q_1}T^{-1}_k a))\nonumber\\
&=e^{\pi t\gamma+2\pi tk_2}h(e^{\pi it\alpha\cP_2}(T^{-1}_k b)^\ast{\scriptstyle{\#}} T^{-1}_k a))\nonumber\\
&=e^{\pi t\gamma+2\pi tk_2}e^{-\pi t\gamma -2\pi tk_2} h((T^{-1}_k b)^\ast{\scriptstyle{\#}}(T^{-1}_k a))\nonumber\\
&=(a,b)_h~.
\end{align}
The corresponding proofs for the other operators are similar. 
The continuous extensions of the unitary  operators $U_j(t)$ and 
$V_j(t)$ to the Hilbert space completion 
$\G$  of the unitary space $(\gA_{pex},(\cdot,\cdot)_h)$ are denoted by the 
same symbols.

Before we continue let us note an estimate for the Hilbert space norm 
$\|\cdot\|_h$. Since the functional $h$ is $\tau$-continuous by assumption, 
there is a $\tau$-continuous seminorm $\gr$ on $\gA_{pex}(\bbbr^2)$ such that
$$
\| a\|^2_h=h((T^{-1}_k b)^\ast{\scriptstyle{\#}}(T^{-1}_k a))\le 
\gr((T^{-1}_kb)^\ast{\scriptstyle{\#}}(T^{-1}_k a)),~ a\in\gA_{pex}(\bbbr^2)~.
$$
By the definition of the topology $\tau$, we can take $\gr$ to be a 
sum of norms
\begin{equation}\label{D3}
\|\cdot\|_{c,d}:=\| e^{2\pi(c_1\cP_1+c_2\cP_2)} e^{2\pi(d_1\Q_1+d_2\Q_2)}~.\|,
\end{equation}
where $c=(c_1, c_2), d=(d_1,d_2)\in\bbbr^2$. Using the formulas (\ref{F5})--
(\ref{F12}) we conclude that the Hilbert space norm $\|\cdot \|_h$ is $\tau$-continuous 
on $\gA_{pex}(\bbbr^2)$.

Now let $s\in \bbbr$ and let $a\in\gA_{pex}(\bbbr^2)$. 
For $c,d\in\bbbr^2$, we then have
\begin{align*}
&\|(U_1(t)-I-s 2\pi i\alpha\Q_1) a\|_{c,d}\\
=&\| e^{2\pi(c_1\cP_1+c_2\cP_2)}
e^{2\pi(d_1\Q_1+d_2\Q_2)}(U_1(t)-I-s 2\pi i\alpha\Q_1)a\|\\
=&\|(e^{2\pi t\alpha c_1} U_1(t)-I-s 2\pi i\alpha\Q_1 -s2\pi\alpha t)
e^{2\pi(c_1\cP_1+c_2\cP_2)}e^{2\pi(d_1\Q_1+d_2\Q_2)}a\|
\end{align*}
Since the norm $\|\cdot \|_h$ can be estimated by sums of 
norms $\|\cdot\|_{c,d}, c,d\in\bbbr^2$, it 
follows from the preceding equality, applied with $s=0$, that 
$U_1(t) a\rightarrow a$ in the Hilbert space $\G$ as $t\rightarrow 0$. 
Therefore, the one-parameter unitary group $t\rightarrow U_1(t)$ 
on the Hilbert space $\G$ is strongly continuous. Next we set $s=t$. 
Then we conclude from the preceding that 
$t^{-1}(U_1(t)-I)a\rightarrow 2\pi i\alpha \Q_1$ in the 
Hilbert space $\G$ as $t\rightarrow 0$. (These facts are obvious in case of the Hilbert space $L^2(\bbbr^2)$, but we do not yet know that 
$(\cdot ,\cdot)_h$ is a multiple of the $L^2$-scalar product 
$(\cdot,\cdot)$.) The corresponding assertions for the other unitary 
groups follow by a similar reasoning.

It is obvious that the operators $U_j(t)$ and $V_j(t)$ 
satisfy the commutation relations
\begin{align*}
&V_j(t)U_j(s)=e^{2\pi i\gamma ts}U_j(s)V_j(t),~V_1(t)U_2(s)=U_2(s)V_1(t),~V_2(t)U_1(s)=U_1(s)V_2(t),\\
&U_1(t)U_2(s)=U_2(s)U_1(t),~V_1(t)V_2(s)=V_2(s)V_1(t)
\end{align*}
for $j=1,2$ and $s,t\in\bbbr$. That is, the unitary one-parameter 
groups $U_j(t)$ and $V_j(t), j=1,2$ and $t\in\bbbr$, fulfill the 
Weyl relation. 
Therefore, by the Stone-von Neumann uniqueness theorem [Pu]
that there exist  a Hilbert space $\cK$ and a unitary 
transformation $\J$ of $\G$ on the Hilbert space 
$\tilde{\Hh}:=L^2(\bbbr^2)\otimes\cK$ such that
\begin{equation}\label{D4}
\J V_j(t)U_l(s)\J^{-1}=\tilde{V}_j(t)\tilde{U}_l(s)~{\rm for}~ j,l=1,2, s,t\in\bbbr,
\end{equation}
where $\tilde{V}_j(t)$ and $\tilde{U}_l(s)$ denote the 
unitary groups on the Hilbert space $\tilde{\Hh}$ defined 
by the same formulas as $V_j(t)$ and $U_l(s)$, respectively.

For $\delta=(\delta_1, \delta_2)\in\bbbr^2_{++}$ and 
$c=(c_1, c_2)\in\bbbc^2$, let $e_{c,\delta}$ denote the 
function \break$e^{2\pi(c_1x_1+c_2 x_2-\delta_1x^2_1-\delta_2 x^2_2)}$ on 
$\bbbr^2$. Let
$$
A_\delta:= (\cP_1+2i\delta_1\Q_1)(\cP_1-2i\delta_1\Q_1)+
(\cP_2+2i\delta_2\Q_2)(\cP_2-2i\delta_2\Q_2)
$$ 
be the operator on unitary space $(\gA_{pex}(\bbbr^2),(\cdot,\cdot)_h)$ and 
let $\tilde{A}_\delta$ denote the operator on $\tilde{\Hh}$ 
given by the same formula. Since $(\cP_j-2i\delta_j\Q_j) e_{0,\delta}=0$ 
for $j=1,2$, it is clear that $A_\delta e_{0,\delta} =0$. 
From (\ref{D3}) it follows that the unitary transformation $\J$ 
intertwines the generators of the unitary groups $V_j (t), U_l(s)$ 
and $\tilde{V}_j(t),\tilde{U}_l(s)$, respectively. Hence we also have 
$\tilde{A}_\delta\J(e_{0,\delta})=0$. It is not diffult to show that 
$\ker\tilde{A}_\delta=e_{0,\delta}\otimes \cK$. Hence 
we conclude that for arbitrary $\delta\in\bbbr^2_{++}$ there 
exist a vector $x_\delta\in\cK$ such that 
$\J(e_{0,\delta})=e_{0,\delta}\otimes x_\delta$.

We next show that the vector $x_\delta$ does not depend on 
$\delta \in\bbbr^2_{++}$. For $\varepsilon\in\bbbr$, let 
$\varepsilon_1:=(\varepsilon,0)$ and $\varepsilon_2:=(0,\varepsilon)$. 
It is not difficult to verify that
\begin{equation}\label{D5}
\varepsilon^{-1}(e_{0,\delta+\varepsilon_j}-
e_{0,\delta})-2\pi\Q^2_j e_{0,\delta}\rightarrow 0
\end{equation}
as $\varepsilon\rightarrow 0$ in the topology $\tau$ on $\gA_{pex}(\bbbr^2)$. 
Since the norm $\|\cdot\|_h$ is $\tau$-continuous as noted above, (\ref{D5}) 
holds also in the Hilbert space $\G$. Therefore, the image under 
$ \J$ of the expression in (\ref{D5}) tends to zero in the 
Hilbert space $\Hh$ as $\varepsilon\rightarrow 0$. This means that
$$
\big[ \varepsilon^{-1}(e_{0,\delta+\varepsilon_j}-e_{0,\delta})-2\pi\Q^2_je_{0,\delta}\big]\otimes x_\delta+e_{0,\delta+\varepsilon_j}\otimes \varepsilon^{-1}(x_{\delta+\varepsilon_j}-x_\delta)\rightarrow 0
$$
as $\varepsilon\rightarrow 0$ in $\Hh$. Using once more the fact that 
(\ref{D5}) holds in $L^2(\bbbr^2)$, we conclude that 
$\varepsilon^{-1}(x_{\delta+\varepsilon_j}-x_\delta)\rightarrow 0$ in $\cK$. 
That is, the partial derivatives of the $\cK$-valued function 
$\delta\rightarrow x_\delta$ vanish on $\bbbr^2_{++}$. Hence 
this function is constant with respect to $\delta\in\bbbr^2_{++}$, 
say $x_\delta= x$.

Recall that $\gA_{pex}(\bbbr^2)$ is defined as the linear span of 
function $x_1^n x_2^m e_{c,\delta}$, where $n,m\in\bbbn_0, c\in\bbbc^2$ 
and $\delta\in\bbbr^2_{++}$. From (\ref{D5}) and the fact that 
$\J(e_{0,\delta})=e_{0,\delta}\otimes x$ it follows that 
$\J(a)=a\otimes x$ for all $a\in\gA_{pex}(\bbbr^2)$. Since 
$\J$ is unitary, we obtain
$$
(a,b)_h=(\J(a),\J(b))_{\tilde{\Hh}}=(a,b)_{L(\bbbr^2)}(x,x)_\cK.
$$
Thus the scalar product $(\cdot,\cdot)_h$ is a positive multiple of 
the $L^2$-scalar product $(\cdot,\cdot)$ on $\gA_{pex}(\bbbr^2)$, that is, $(\cdot,\cdot)_h=\nu(\cdot,\cdot)$, where $\nu:=(x,x)$.

By the latter and the definition (\ref{D1}) of the scalar product 
$(\cdot,\cdot)_h$, we have 
\begin{align*}
h(b^\ast{\scriptstyle{\#}} a)&=(T_k a,T_k b)_h=\nu(T_k a, T_k b)=\nu(T^2_k a,b)\\
&=\nu(e^{2\pi\alpha_k\Q_1-\pi\beta_k\cP_1}\otimes 
e^{2\pi\beta_k\Q_2+\pi\alpha_k\cP_2}a,b)
\end{align*}
for $a,b\in\gA(\bbbr^2)$. Now we set $b=f_\varepsilon$, where 
$f_\varepsilon$ is the approximate identity from Proposition 14. 
Taking the limit $\varepsilon \rightarrow +0$ in the topology $\tau$ and using (\ref{apid}), we obtain
$$
h(a)=\nu\iint e^{2\pi\alpha_kx_1+2\pi\beta_kx_2} 
a(x_1+{\tfrac{\beta_k}{2}}i, x_2-{\tfrac{\alpha_k}{2}}i) dx_1 dx_2~.
$$
Arguing as in the proof of Proposition 22, we substitute $(x_1,x_2)(x_1-i\beta_k/2,x_2+i\alpha_k/2)$ in the latter formula and obtain 
$h(a)=\nu h_k(a),a\in\gA_{pex}(\bbbr^2)$.\hfill $\Box$

\mn
{\bf Theorem 25.} {\it Let $k=(k_1,k_2)\in\bbbzz^2$ and 
let $c=(c_1,c_2),d=(d_1,d_2)\in\bbbr^2$ be such that 
$8|c_1 d_1|<1$ and $8|c_2 d_2|<1$. Suppose that $h$ is a 
faithful positive linear functional on the $\ast$-algebra 
$\gA(\bbbr^2)$ such that
\begin{align}\label{D6}
&h(a(x_1{+}\beta s, x_2{+}\alpha s))=
e^{-2\pi\gamma(s+t)-4\pi(k_1s+k_2t)}h(a),\nonumber\\
&|h(a)|\le C\| S(e^{c\Q})S(e^{d\cP})a\|
\end{align}
for some positive constant $C$ and for all $s,t\in\bbbr$ 
and $a\in\gA(\bbbr^2)$. \\
Then there is a positive 
constant $\nu$ such that $h=\nu h_k$.}

\mn
{\bf Proof.} By Proposition 24, we have $h(a)=\nu h_k(a)$ for all 
$a\in\gA_{pex}(\bbbr^2)$. By Lemma 10, 
$\gA_{pex}(\bbbr^2)$ is dense in $\gA(\bbbr^2)$ relative to the norm 
$\| S(e^{c\Q})S(e^{d\cP}){\cdot}\|$. Since $h$ and $h_k$ are both 
continuous with respect to this norm, it follows that $h=\nu h_k$ on 
$\gA(\bbbr^2)$.\hfill $\Box$

\mn
Remark 4. It is likely that the assertion of Theorem 25 remains 
valid if we assume only the $\tau$-continuity of the functional $h$ 
instead of inequality (\ref{D6}) with $8|c_j d_j|<1, j=1,2$. For this 
it would be sufficient to know that $\gA_{pex}(\bbbr^2)$ is $\tau$-dense 
in $\gA(\bbbr^2)$.

\hfill\eject\noindent
{\bf 6. The real quantum plane}

\mn
{\it 6.1} In order to motivate the definitions given below, 
we first recall the following 
well-known representation of the Lie algebra $gl_2(\bbbr)$ on the 
$C^\infty$-functions of $\bbbr^2$: The action of the 
generators $e,f,h_1,h_2$ of $gl_2(\bbbr)$ satisfying the relations
$$
[h_1,e]=e, [h_1,f] =-f, [h_2,e]=-e, [h_2,f]=f, [e,f]=h_1-h_2
$$
is given by the formulas
\begin{equation}\label{E1}
e=-y {\frac{\partial}{\partial x}}, 
f=-x{\frac{\partial}{\partial y}}, 
h_1 = -x{\frac{\partial}{\partial x}}, 
h_2 = -y{\frac{\partial}{\partial y}}.
\end{equation}
The groups of relations (\ref{A5})--(\ref{A7}) 
and (\ref{C6})--(\ref{C8}) can be interpreted as quantum 
versions of the formulas (\ref{E1}). Likewise, we see that 
the linear mappings $\cD_x^q$ and $\cD^q_y$ can be thought as 
$q$-deformations of the partial derivatives ${\frac{\partial}{\partial x}}$ 
and ${\frac{\partial}{\partial y}}$, respectively. We shall not make this 
more explicit because we shall not need these details.

Let us explain the underlying idea of our construction of the real quantum plane in the classical situation. We consider the plane $\bbbr^2$ 
as the direct sum of the four quarter planes 
$\bbbr^{++},\bbbr^{-+},\bbbr^{+-},\bbbr^{--}$ that are 
glued together along the two coordinate axis.
On the level of functions this means that a continuous function $f$ on 
$\bbbr^2\backslash\{(0,0)\}$ is given by a 4-tuple $(f_{++},f_{-+},f_{+-},f_{--})$ 
of continuous functions on the quarter planes satisfying the boundary 
conditions
\begin{align*}
&f_{++}(+0,y)=f_{-+}(-0,y), ~f_{+-}(+0,y)=f_{--}(-0,y),\\
&f_{++}(x,+0)=f_{+-}(x,-0),~f_{-+}(x,+0)=f_{--}(x,-0).
\end{align*}
We now turn to the quantum case and consider the direct sum 
\begin{equation}\label{algebrab}
\gB(\bbbr^2)_4:= {\tilde \gA}(\bbbr^2)\oplus 
{\tilde \gA}(\bbbr^2)\oplus {\tilde \gA}(\bbbr^2) 
\oplus {\tilde \gA}(\bbbr^2)
\end{equation}
of the four $\ast$-algebras ${\tilde \gA}(\bbbr^2)$ corresponding to 
the four quarter planes. 
Recall that product and involution of the 
$\ast$-algebra ${\tilde \gA}(\bbbr^2)$ are given by formulas 
(\ref{atimesb}) and (\ref{star}). We could have also taken here the 
$\ast$-algebra $\gA(\bbbr^2)$. The reason why we prefer to work 
with ${\tilde \gA}(\bbbr^2)$ is that on direct sums of Hilbert 
spaces $L^2(\bbbr^2)$ it is easier 
to describe self-adjoint extensions of the operators $\Phi(E^\prime)$ and
$\Phi(F^\prime)$. Since the elements of ${\tilde \gA}(\bbbr^2)$ 
correspond to ``functions on the quantum 
quarter plane which vanish at the boundaries'',
no boundary condition occurs and we can just take the 
direct sum of the four $\ast$-algebras. 

Next we want to make the direct sum 
$\cB(\bbbr^2_q)=\cO(\bbbr^2_q) +\gB(\bbbr^2)_4$ of the 
vector spaces $\cO(\bbbr^2_q)$ and $\gB(\bbbr^2)_4$ into a 
left $\U_q(gl_2(\bbbr))$-module $\ast$-algebra. This means that 
we have to define the products of the generators $x$ and $y$ of 
$\cO(\bbbr^2_q)$ by a 4-tuple $a=(a_1,a_2,a_3,a_4)$, $a_j \in 
{\tilde \gA}(\bbbr^2)$, and left actions of the generators 
$E, F, K_1, K_2$ of $\U_q(gl_2)$ on $a$.  
Let us look for a moment at the classical case and consider a quarter 
plane $\bbbr^{\epsilon \epsilon^\prime}$, where 
$\epsilon, \epsilon^\prime \in \{+,-\}$. From the formulas (\ref{E1}) we see
that if we pass from $\bbbr^{++}$ to 
$\bbbr^{\epsilon \epsilon^\prime}$ then $x$ has to be replaced 
by $\epsilon x$, $y$ by $\epsilon^\prime y$, $e$ by 
$\epsilon \epsilon^\prime e$ and $f$ by $\epsilon \epsilon^\prime f$, 
while $h_1$ and $h_2$ remain unchanged. This suggests 
to take the following definitions in the quantum case:
\begin{align*}
xa&=(xa_1,xa_2,{-}xa_3,{-}xa_4),~ ya =(ya_1,{-}ya_2,ya_3,{-}ya_4),\\
\Phi(E^\prime) a &= (\Phi(E^\prime) a_1, {-}\Phi (E^\prime)a_2, {-}\Phi(E^\prime) a_3, \Phi(E^\prime)a_4),\\
\Phi(F^\prime)a &=  (\Phi(F^\prime)a_1,{-}\Phi(F^\prime)a_2,{-}\Phi(F^\prime)a_3,\Phi(F^\prime)a_4),\\
\Phi(K_j)a &= (\Phi(K_j)a_1,\Phi(K_j)a_2,\Phi(K_j)a_3,\Phi(K_j)a_4),~ j=1,2,
\end{align*}
where $a=(a_1,a_2,a_3,a_4)$, $a_j \in {\tilde \gA}(\bbbr^2)$. Reall that 
the action of the operators $x \equiv\Phi(x)$, 
$y \equiv \Phi (y)$ and $\Phi(f)$, $f=E^\prime, F^\prime, K_1, K_2$, 
are described by
formulas (\ref{I11})--(\ref{I13}). Since 
${\tilde \A}(\bbbr^{++}_q)=\cO(\bbbr^2_q) +\gA(\bbbr^2)$ 
is a left $\U_q(gl_2(\bbbr))$-module $\ast$-algebra with left action $\Phi$, 
it is easily verified that 
the preceding formulas define indeed a unique $\ast$-algebra structure 
on $\cB(\bbbr^2_q)_4$ and that $\cB(\bbbr^2_q)_4$ becomes a left 
$\U_q(gl_2(\bbbr))$-module $\ast$-algebra in this manner. 

From the preceding it is clear that the functional $\tilde{h}$ on 
$\gB(\bbbr^2)_4$ defined by
\begin{equation}\label{htildeb}
\tilde{h}(a)=\tilde{h}(a_1)+\tilde{h}(a_2)+\tilde{h}(a_3)+\tilde{h}(a_4),~ 
a=(a_1, a_2, a_3, a_4)\in \gB(\bbbr^2)_4,
\end{equation}
is $\U_q(gl_2(\bbbr))$-covariant with respect to $\chi$, where $\tilde{h}(a_j)$ is given by (\ref{htilde}).

\mn
{\it 6.2}  The vector space $\gB(\bbbr^2)_4$ is also a 
dense domain of the Hilbert space
\begin{equation}\label{hilberth}
\Hh =L^2(\bbbr^2)\oplus L^2(\bbbr^2)\oplus L^2 (\bbbr^2)\oplus L^2(\bbbr^2).
\end{equation}

\mn
From the formulas (\ref{I12}) and (\ref{I13}) it is clear 
that the operators $x$, $y$ and $\Phi(q^{-1/4}K_j)$, $j=1,2$, are 
essentially self-adjoint on the domain $\gB(\bbbr^2)_4$. From 
(\ref{I11}) and Lemma 5(ii) it follows that the operators 
$\Phi(E^\prime)$ and $\Phi(F^\prime)$ are symmetric but not essentially 
self-adjoint on the domain $\gB(\bbbr^2)_4$ in $\Hh$. (Indeed, by 
Lemma 5(ii) the adjoint of $L_\alpha$ is the operator 
$R_\alpha$ and $R_\alpha$ is easily seen to be 
a proper extension of $L_\alpha$.) That the operators 
$\Phi(E^\prime)$ and $\Phi(F^\prime)$ are not essentially self-adjoint on 
the domain $\gB(\bbbr^2)_4$ is not surprising, because $\gB(\bbbr^2)_4$ 
contains only ``functions which vanish at the boundaries of the four 
quantum quarter planes''. In the classical case 
the corresponding symmetric operators 
$ie=-y i\frac{\partial}{\partial x}$ and 
$if=-xi\frac{\partial}{\partial y}$ are also not essentially self-adjoint 
when the functions in the domain have boundary values zero at the 
$x$- and $y$-axis. Our next step is to ``glue 
together the four quarter quantum planes'' 
and to obtain self-adjoint extensions of 
the symmetric operators $\Phi (E^\prime)$ and $\Phi (F^\prime)$ 
in this manner. 

In order to explain the glueing procedure let us return to the classical 
situation and consider the direct sum of four quarter planes. 
Since the corresponding functions have the boundary values zero at the 
$x$- and $y$-axis, the operators 
$ie=-y i\frac{\partial}{\partial x}$ and 
$if=-xi\frac{\partial}{\partial y}$ are not essentially 
self-adjoint. 
The particular self-adjoint extensions of the symmetric 
operators $ie$ and $if$ 
we are interested in are determined by the boundary conditions
\begin{equation}\label{E5a}
f(+0,y)=f(-0,y)~{\rm and}~f(x,+0)=f(x,-0).
\end{equation}
For simplicity, let us first look how the two upper quarter 
planes $\bbbr^{-+}$ and $\bbbr^{++}$ are glued together along 
the positive $y$-axis. We identify a function $f_{-+}$ on $\bbbr^{-+}$ 
with the function $f_{++}$ on $\bbbr^{++}$ given by 
$f_{++}(x,y):=f_{-+}(-x,y)$, $y>0,x\in\bbbr$. Let $D_x$ denote the symmetric operator 
$-i\frac{\partial}{\partial x}$ on $L^2(\bbbr^{++})$ with 
boundary condition $f(+0,y)=0$, $y>0$. Then the operator 
$T_0=-i\frac{\partial}{\partial x}$ on $\bbbr^{-+}\cup\bbbr^{++}$ 
with boundary condition zero at the positive $y$-axis acts 
on the Hilbert space $L^2(\bbbr^{++})\oplus L^2(\bbbr^{++})$ and 
has the form 
$$
T_0= \left( \begin{matrix}  D_x &0\\  
                    0   &-D_x \end{matrix} \right)
$$
Let $T$ denote the self-adjoint extension of $T_0$ 
with boundary condition $f(+0,y)=f(-0,y)$, $y>0$, and 
let $J_0$ be the symmetry (that is, self-adjoint unitary)
$$
J_0={\frac{1}{\sqrt{2}}}~\left(\begin{matrix} 1&1\\
1&-1 \end{matrix} \right)~.
$$
One easily verifies that
$$
J_0 T J_0 = \left( \begin{matrix} 0 &D_x\\
D_x^\ast &0 \end{matrix}\right)~,
$$
where $D^\ast_x$ is the adjoint of the closed symmetric operator 
$D_x$ on the Hilbert space $L^2(\bbbr^{++})$.

In a similar manner we proceed with the four quarter 
planes $\bbbr^{++}$, $\bbbr^{-+}$, $\bbbr^{+-}$, $\bbbr^{--}$. As 
above, we identify functions on $\bbbr^{-+},\bbbr^{+-},\bbbr^{--}$ 
with the corresponding functions on $\bbbr^{++}$, and let 
$D_y=-i\frac{\partial}{\partial y}$ on $L^2(\bbbr^{++})$ with boundary 
condition $f(x,+0)=0, x>0$. Further, let 
$T_{10}=-i\frac{\partial}{\partial x}$ and 
$T_{20}=-i\frac{\partial}{\partial y}$ be the 
symmetric operators with boundary values zero at the positive 
$y$- resp. $x$-axis acting on the Hilbert space 
\begin{equation*}
\G=L^2(\bbbr^{++})\oplus L^2(\bbbr^{++})\oplus 
L^2(\bbbr^{++})\oplus L^2(\bbbr^{++}).
\end{equation*}
Consider the symmetries
$$
J_1={\frac{1}{\sqrt{2}}}\left(\begin{matrix} 1&~~1&0&~~0\\ 1&-1&0&~~0\\ 
0&~~0&1&~~1\\ 0&~~0&1&-1\end{matrix}\right)~,~~
J_2={\frac{1}{\sqrt{2}}}\left(\begin{matrix} 1&0&~~1&~~0\\ 0&1&~~0&~~1\\ 
1&0&-1&~~0\\ 0&1&~~0&-1\end{matrix}\right)
$$
on the Hilbert space $\G$. Since $J_1$ and $J_2$ commute, their product
\begin{equation}\label{E6}
J:=J_1J_2=\frac{1}{2}\left(\begin{matrix}1&~~1&~~1&~~1\\
1&-1&~~1&-1\\
1&~~1&-1&-1\\
1&-1&-1&~~1\end{matrix}\right)
\end{equation}
is again a symmetry. We shall transform all structures by means of the 
symmetry $J$. 

In order to save space, let us introduce  
abbreviations for some operator matrices. If $z_1$, $z_2$ and $z$ is an operator on a Hilbert space $\cK_0$, we set
\begin{align*}
&\sigma_1 (z_1,z_2)=\left( \begin{matrix} z_1&~~0&0&~~0\\
                                   0&z_2&0&~~0\\
                                   0&~~0&z_1&~~0\\
                                   0&~~0&0&z_2  
                   \end{matrix} \right)~, 
\quad \sigma_2 (z_1,z_2)=\left( \begin{matrix}  z_1&0&~~0&~~0\\
                                   0&z_1&~~0&~~0\\ 
                                   0&0&z_2&~~0\\ 
                                   0&0&~~0&z_2
                   \end{matrix} \right)~,
\\
\\
&\theta_1(z)=
                  \left( \begin{matrix} 0~&z&0~&0\\  
                                         z^\ast &0&0~&0\\ 
                                         0~&0&0~&z\\                                          0~&0&z^{\ast}&0
                  \end{matrix} \right)~,
\quad \theta_2(z)=\left( \begin{matrix}  0~&0~&z&0\\ 
                                         0~&0~&0&z\\ 
                                         z^\ast&0~&0&0\\                                          0~&z^\ast &0&0
                    \end{matrix} \right)~,
\\
\\
&\kappa_1 (z)=\left( \begin{matrix} 0~&0&0~&z\\ 
                                          0~&0&z^\ast &0\\                                           0~&z&0~&0\\                                           z^\ast &0&0 &0
                      \end{matrix} \right)~,
\quad 
\kappa_2 (z)=\left( \begin{matrix} 0~&0~&0&z\\ 
                                         0~&0~&z&0\\ 
                                         0~&z^\ast &0&0\\                                          z^\ast &0~&0&0
                     \end{matrix} \right)~
\end{align*}
and $\sigma_j(z):=\sigma_j(z,-z)$, $j=1,2$. The matrices $\sigma_j(z_1,z_2), \theta_j(z)$ and 
$\kappa_j (z)$ act as operators on the 
Hilbert space $\cK_0 \oplus \cK_0 \oplus K_0 \oplus \cK_0$. 

Let us now continue the above discussion. In terms of the preceding 
notation we have $T_{10}=\sigma_1(D_x)$ and $T_{20}=\sigma_2(D_y)$. 
It is not difficult to check that the  
self-adjoint extensions 
$T_1$ and $T_2$ of $T_{10}$ and $T_{20}$, respectively, with 
boundary condition (\ref{E5a}) satisfy the relations
\begin{align*}
J T_1J=J_2\kappa_1(D_x)J_2=\theta_1(D_x),~~
J T_2J=J_1\kappa_2(D_y)J_1=\theta_2(D_y).
\end{align*}
That is, the particular self-adjoint extensions $T_1$ and $T_2$ are 
characterized in a simple manner by means of the symmetry transformation $J$. 

Let $x$ and $y$ denote the multiplication operators by the 
coordinate functions $x$ and $y$, respectively, on $L^2(\bbbr^{++})$. 
Then the multiplication operators $M_x$ and $M_y$ by the coordinate 
functions on $L^2(\bbbr^2)$ can be expressed as
\begin{equation}\label{mxmy}
J M_xJ=J\sigma_1(x)J=\theta_1(x),~ J M_yJ=J\sigma_2(y)J=\theta_2(y)~.
\end{equation}
Further, the self-adjoint operators $i{\hat e}
=-iy{\frac{\partial}{\partial x}}$ and 
$i{\hat f}=-ix{\frac{\partial}{\partial y}}$ with boundary condition 
(\ref{E5a}) can be written as 
\begin{equation}\label{ef}
J i\hat{e} J=J M_x T_2J =\kappa_2 (x D_y),~~ 
J i\hat{f} J =J M_y T_1J =\kappa_1(y D_x).
\end{equation}

Finally, we also transform the structure of the $\ast$-algebra 
of functions on $\bbbr^2$ under the symmetry $J$. If 
we consider functions $f, g$ on $\bbbr^2$ as 4-tuples of 
functions on the quarter planes, then the product $f \cdot g$ 
is transformed under $J$ as  
\begin{equation}\label{twp}
f {\scriptscriptstyle{\circ}} g := J( J(f){\cdot} J(g)).
\end{equation}
More explicitly, for $f=(f_1,f_2,f_3,f_4)$ and $g=(g_1,g_2,g_3,g_4)$ 
we compute 
\begin{align}\label{twpcomp}
f {\scriptscriptstyle{\circ}} g =1/4 (f_1g_1+f_2g_2+f_3g_3+f_4g_4, 
f_1g_2+f_2g_1+f_3g_4+f_4g_3,\nonumber\\
f_1g_3+f_2g_4+f_3g_1+f_4g_2,
f_1g_4+f_2g_3+f_3g_2+f_4g_1),
\end{align}
where $f_jg_k$ and $g_kf_j$ mean the usual pointwise products of functions on the quarter plane. 
Obviously, the involution of functions is invariant under $J$, that is, 
we have $J(f^\star)=J(f)^\star$.

Thus, we have described the operators and the algebra of functions on the 
plane by means of the symmetry $J$ and the corresponding operators and 
algebras of the quarter planes. The advantage of this approach is that 
it gives a convenient algebraic form for the particular self-adjoint 
extensions of the first order differential operators. There is 
the disadvantage that the algebra product has to be changed too.
The preceding formulas and considerations will serve as guiding motivation 
for the construction of the real quantum plane in the next subsection. 

\mn 
{\it 6.3} In this subsection we develop the basics of the 
construction of the real quantum plane. Our starting point is the 
description of the quantum quarter plane by means of the 
left $\U_q(gl_2(\bbbr))$-module $\ast$-algebra ${\tilde \A}(\bbbr^{++}_q)$ developed in 4.2. In what follows we assume that
\begin{equation}\label{assump}
0 < |\gamma | < 1/3.
\end{equation}

 Let us begin by defining the left action of $\U_q(gl_2)$. 
Remembering the formulas (\ref{ef}), (\ref{I11}) and (\ref{I12}), 
we consider the self-adjoint operators on the Hilbert space $\Hh$ 
(see (\ref{hilberth})) defined by the $4{\times} 4$-operator matrices 
\begin{align*}
\sE :=\kappa_1(L_\alpha\otimes e^{2\pi\beta\Q_2}),~
\fF :=\kappa_2(e^{2\pi\alpha\Q_1}\otimes L_{\beta,\alpha}),~
\fK_1 := I(e^{\pi\beta\cP_1}),~
\fK_2 := I(e^{\pi\alpha\cP_2}).
\end{align*}
Here $I(z)$ denotes the diagonal matrix with diagonal entries equal to $z$. 
Note that the entries of these matrices are just the 
operators occuring in formulas (\ref{I11})--(\ref{I12}). 
Let $\gA(\bbbr^2)_4$ be the domain
\begin{equation}\label{ar2q}
\gA(\bbbr^2)_4:= \gA_{12}(\bbbr^2)\oplus\gA_2 (\bbbr^2) 
\oplus\gA_1 (\bbbr^2)\oplus \gA (\bbbr^2)
\end{equation}
in the Hilbert space $\Hh$, where
\begin{align*}
\gA_1(\bbbr^2)&=f_\alpha(\cP_1)^{-1}\gA(\bbbr^2),~~
\gA_2(\bbbr^2)=f_{\beta}(\cP_2)^{-1}\gA(\bbbr^2),\\
\gA_{12}(\bbbr^2)&=f_\alpha(\cP_1)^{-1}\gA(\bbbr^2)+
f_{\beta}(\cP_2)^{-1}\gA(\bbbr^2).
\end{align*}
It is clear that $\gA(\bbbr^2) \subseteq\gA_1(\bbbr^2)\cup \gA_2(\bbbr^2)
\subseteq\gA_{12}(\bbbr^2)$. Since 
$f_\alpha(\cP_j)\gA(\bbbr^2) \subseteq\gA(\bbbr^2)$ for $j=1,2$,
$\gA(\bbbr^2)_4$ contains the domain $\gB(\bbbr^2_q)$ 
defined by (\ref{algebrab}). In particular, $\gA(\bbbr^2)_4$ is dense 
in $\Hh$. Using the fact that the operator $L_\alpha$ in 
$L^2(\bbbr)$ has the adjoint operator 
$R_\alpha = e^{-2\pi\alpha\Q}f_\alpha(\cP)$ by Lemma 5(ii), 
it is clear that $\gA(\bbbr^2)_4$ is contained in the domains of the 
operators $\sE, \fF, \fK_j$, $j=1,2$. Define
\begin{equation}\label{leftac}
\Phi (E^\prime)=\sE{\upharpoonright} \gA(\bbbr^2)_4,~
\Phi (F^\prime)=\fF{\upharpoonright} \gA(\bbbr^2)_4,~
\Phi (q^{-1/4}K_j)=\fK_j{\upharpoonright} \gA(\bbbr^2)_4.
\end{equation}

\mn
{\bf Proposition 26.}  (i) {\it $\gA(\bbbr^2)_4$ is contained in the domains of the operator products $\Phi(E^\prime)\Phi(F^\prime)$,
 $\Phi(F^\prime)\Phi(E^\prime)$, $\Phi(E^\prime)\Phi(K_j)$, 
 $\Phi(K_j)\Phi(E^\prime)$, $\Phi(F^\prime)\Phi(K_j)$, 
 and $\Phi(K_j)\Phi(F^\prime)$ for $j=1,2$.}\\
(ii) {\it The operators $\Phi(E^\prime), \Phi(F^\prime)$ and $\Phi(K_j)$ 
satisfy the defining relations in terms of the generators 
$E^\prime, F^\prime$ and $K_j$ of the algebra $\U_q(gl_2)$.}\\ 
(iii) {\it The operators $\Phi(E^\prime), \Phi(F^\prime)$ and $\Phi(q^{-1/4}K_j)$ 
are essentially self-adjoint on the domain $\gA(\bbbr^2)_4$.}

\mn
{\bf Proof.} (i): We show that $\gA(\bbbr^2)_4$ is contained in the 
domain of the product $\Phi(E^\prime)\Phi(F^\prime)$. By the definition 
of operators $\Phi(E^\prime)$ and $\Phi(F^\prime)$ this means that 
\begin{align*}
&\gA_{12}(\bbbr^2) {\subseteq} \cD((L_\alpha  {\otimes} e^{2\pi\beta \Q_2})(e^{2\pi\alpha Q_1} {\otimes} R_{\beta})),~\gA_2(\bbbr^2) {\subseteq}\cD((R_\alpha {\otimes} e^{2\pi\beta\Q_2})(e^{2\pi\alpha\Q_1} {\otimes} R_{\beta})),\\
&\gA_1(\bbbr^2) {\subseteq} \cD((L_\alpha  {\otimes} e^{2\pi\beta \Q_2})
(e^{2\pi\alpha Q_1} {\otimes} L_{\beta})),~
\gA(\bbbr^2) {\subseteq} \cD((R_\alpha {\otimes} e^{2\pi\beta\Q_2})
(e^{2\pi\alpha\Q_1} {\otimes} L_{\beta})).
\end{align*}
The last inclusion for the fourth component is obvious. Let us verify the 
assertion for the third component of $\Phi(E^\prime)\Phi(F^\prime)$. 
That is, we have to show that each vector 
$\eta = f_\alpha(\cP_1)^{-1} \zeta \in \gA_1(\bbbr^2)$, 
where $\zeta \in \gA(\bbbr^2)$, belongs to the domain of the product 
$(L_\alpha  \otimes e^{2\pi\beta \Q_2})
(e^{2\pi\alpha Q_1} \otimes L_{\beta})$. In order to prove this, 
it suffices to check that for all $\xi \in \gA(\bbbr)$ the vector 
$f_\alpha(\cP)^{-1} \xi$ belongs to the domain of the operator 
$e^{2\pi\alpha \Q}$. Applying the Fourier transform we see that the 
latter is equivalent to the fact that for arbitrary $\xi \in \gA(\bbbr)$ the function $f_\alpha(x)^{-1} \xi (x)$ belongs to the 
domain of the operator $e^{-2\pi\alpha \cP}$. The 
assumption $0 < |\gamma | < 1/3$ implies that the infimum of the 
holomorphic function $f_\alpha(x)$ on the strip 
$0 < {\rm Im}~x < \alpha$ when $\alpha > 0$ resp. 
$0 > {\rm Im}~x > -\alpha$ when $\alpha < 0$ is positive. It 
follows from the characterization of the operator $e^{-2\pi\alpha \cP}$ 
given in Lemma 4 that the function $f_\alpha(x)^{-1} \xi$ belongs 
to the domain of this operator. The corresponding proofs for the first and the seond components of $\Phi(E^\prime)\Phi(F^\prime)$ and for the other 
operator products are similar. \\
(ii): We carry out the proof for the relation 
$E^\prime F^\prime-F^\prime E^\prime 
=\lambda(K_1^2K_2^{-2}-K_1^{-2}K_2^2)$ of the algebra $\U_q(gl_2)$. 
For the other defining relations 
these verifications are straigthforward and will be omitted. 
Computing the commutator of the operator matrices $\sE$ and $\fF$, 
we obtain a diagonal operator having the operators
\begin{equation*}
\begin{split}
A_1 &:= (L_\alpha  \otimes e^{2\pi\beta \Q_2})
(e^{2\pi\alpha Q_1} \otimes R_{\beta}) - 
(e^{2\pi\alpha \Q_1} \otimes L_{\beta})
(R_\alpha \otimes e^{2\pi\beta\Q_2}),\\
A_2 &:= (R_\alpha \otimes e^{2\pi\beta\Q_2})
(e^{2\pi\alpha\Q_1} \otimes R_{\beta}) - 
(e^{2\pi\alpha\Q_1}\otimes L_{\beta})
(L_\alpha \otimes e^{2\pi\beta\Q_2}),\\
A_3 &:=(L_\alpha  \otimes e^{2\pi\beta \Q_2})
(e^{2\pi\alpha Q_1} \otimes L_{\beta}) - 
(e^{2\pi\alpha \Q_1} \otimes R_{\beta})
(R_\alpha \otimes e^{2\pi\beta\Q_2}),\\
A_4 &:= (R_\alpha \otimes e^{2\pi\beta\Q_2})
(e^{2\pi\alpha\Q_1} \otimes L_{\beta}) - 
(e^{2\pi\alpha\Q_1}\otimes R_{\beta})
(L_\alpha \otimes e^{2\pi\beta\Q_2})
\end{split}
\end{equation*}
as diagonal entries. By (i), for any vector $\eta = (\eta_1,\eta_2,\eta_3,\eta_4) \in \gA(\bbbr^2)$ the $j$-th component $\eta_j$ belongs to the domain of each of the two summands of the operator $A_j$, $j=1,2,3,4$. Using this fact we compute the terms $A_j\eta_j$ and obtain 
$$A_j\eta_j  =\lambda (e^{2\pi\beta\cP_1}\otimes e^{-2\pi\alpha\cP_2}
-e^{-2\pi\beta\cP_1}\otimes e^{2\pi\alpha\cP_2}) \eta_j.
$$
Thus, $\Phi(E^\prime)\Phi(F^\prime)\eta-\Phi(F^\prime) 
\Phi(E^\prime) \eta  = \lambda(\Phi(K_1)^2\Phi(K_2)^{-2}-
\Phi(K_1)^{-2}\Phi(K_2)^2) \eta $. 

(iii): By Lemma 4(ii), the domain $\cD_\delta$ is a core 
for the self-adjoint operator $e^{c\cP}$, $c \in \bbbr$. 
Since $\cD_\delta \otimes \cD_\delta \subseteq \gA(\bbbr^2)$, the 
operator $\Phi(K_j)$ is essentially self-adjoint 
even on the smaller domain $\gB(\bbbr^2)$. For the operators 
$\Phi(E^\prime)$ and $\Phi(F^\prime)$ the assertion follows from 
statements (iii) and (iv) of 
Lemma 5. \hfill$\Box$

\sn
It is easily seen that the domain $\gA(\bbbr^2)_4$ is {\it not invariant} under the operators $\Phi (E^\prime)$ and $\Phi(F^\prime)$.
Therefore, we do not get a $\ast$-representation of the whole $\ast$-algebra 
$\U_q^{tw}(gl_2(\bbbr))$ on the domain  $\gA(\bbbr^2)_4$ of the 
Hilbert space $\G$. However, the actions of the generators 
$f= E^\prime, F^\prime, K_1, K_2$ satisfy the defining relations and they 
have the hermiticity properties of the $\ast$-algebra 
$\U_q^{tw}(gl_2(\bbbr))$. 

Our next step is to define a $\ast$-algebra structure on $\gA(\bbbr^2)_4$.
Recall that we used the $\ast$-algebra ${\tilde \gA}(\bbbr^{++}_q)$
with product $\natural$ and involution $\star$ (see 4.2) as $\ast$-algebra 
of functions on the quantum quarter plane. As in the case of functions 
on $\bbbr^2$ we transform the componentwise product of 4-tuples under 
the symmetry $J$ (see (\ref{twp}) and (\ref{twpcomp})). That is, for
$a=(a_1,a_2,a_3,a_4), b=(b_1,b_2,b_3,b_4)\in \gA(\bbbr^2)_4$ we define
\begin{align}
\begin{split}\label{twpa} 
a {\scriptscriptstyle{\circ}} b =1/4 (&a_1\natural b_1+a_2\natural b_2+a_3\natural b_3+a_4\natural b_4, a_1\natural b_2+a_2\natural b_1+a_3\natural b_4+a_4\natural b_3,
\\
&a_1\natural b_4+a_2\natural b_3+a_3\natural b_2+a_4\natural b_1, 
a_1\natural b_3+a_2\natural b_4+a_3\natural b_1+a_4\natural b_2),
\end{split} 
\\
\label{twia}
a^\star = (&a_1^\star,a_2^\star,a_3^\star,a_4^\star).
\end{align}
For the product $a_j \natural b_k$ we take the 
symmetrized version in formula (\ref{atimesb}): 
\begin{align}\label{atimesb1}
a_j \natural b_k = 
(e^{{\tfrac{\pi}{4}}\beta\cP_1} e^{- {\tfrac{\pi}{2}} \alpha\Q_1}
\otimes e^{-{\tfrac{\pi}{4}}\alpha\cP_2} 
e^{- {\tfrac{\pi}{2}} \beta\Q_2}) a_j {\scriptstyle{\#}} 
(e^{-{\tfrac{\pi}{4}}\beta\cP_1} e^{- {\tfrac{\pi}{2}} \alpha\Q_1}
\otimes e^{{\tfrac{\pi}{4}}\alpha\cP_2} 
e^{- {\tfrac{\pi}{2}} \beta\Q_2}) b_k .
\end{align}
The involution $a_j^\star$ of a $a_j$ is defined as in (\ref{star}) by 
\begin{equation}
\label{star1}
a_j^\star (x_1,x_2) = 
(e^{\pi (\beta \cP_1 - \alpha \cP_2)} \bar{a}_j)(x_1,x_2).
\end{equation}

\mn
{\bf Proposition 27.} {\it The vector space $\gA(\bbbr^2)_4$ is a 
$\ast$-algebra with product ${\scriptscriptstyle{\circ}}$ and involution $\star$ given 
by the formulas (\ref{twpa}), (\ref{atimesb1}), (\ref{twia}) and (\ref{star1}).}

\mn
{\bf Proof.} By arguing as in the proof of assertion (i) of Proposition 26 it follows from assumption (\ref{assump}) that the components 
$a_j, b_k$ in formula (\ref{atimesb1}) are in the corresponding 
operator domains and in a domain $\cD_{\mu,\nu}$ for certain 
$\mu, \nu \in \bbbr^2$. Thus, the product (\ref{atimesb1}) is indeed 
well-defined. That is the reason why we have chosen the symmetrized 
version in (\ref{atimesb}) rather than the two other formulas in 
(\ref{atimesb}). Note that for $a, b \in \gA(\bbbr^2)$ all three 
formulas in (\ref{atimesb}) coincide, but we dealing now with larger 
classes of symbols. 

Now we prove that for $a, b \in \gA(\bbbr^2)_4$ the components of the 
product $a \natural b$ belong again to the corresponding component space in (\ref{ar2q}). As a sample, we show that 
$\gA_{12}(\bbbr^2)\natural \gA_{12}(\bbbr^2) \subseteq \gA_{12}(\bbbr^2)$. Let $a, b \in f_\alpha(\cP_1)^{-1}\gA(\bbbr^2)$. From formulas 
(\ref{F9}) and (\ref{atimesb1}) we get
\begin{equation}\label{fp1}
f_\alpha(\cP_1) (a \natural b) = 
q^{1/2}  (f_\alpha(\cP_1) a) \natural e^{2 \pi \beta \cP_1} b
+ q^{-1/2} e^{-2 \pi \beta \cP_1} a \natural f_\alpha(\cP_1) b.
\end{equation}
Since $f_\alpha(\cP_1)a$ and $f_\alpha(\cP_1)b$ are in 
$\gA_2(\bbbr^2)$, it follows from formulas (\ref{F5})--(\ref{F12}) that 
$f_\alpha(\cP_1) (a \natural b)$ is in the domain of all operator 
$e^{2 \pi c \Q}e^{2 \pi d \cP}$, $c, d \in \bbbr^2$. Thus, we have 
$f_\alpha(c\cP_1) (a \natural b) \in \gA_2(\bbbr^2)$. If 
$a, b \in f_{\beta}(\cP_2)^{-1}\gA(\bbbr^2)$, then we use formula 
(\ref{F11}) rather than (\ref{F9}) and obtain the identity 
$$
f_{\beta}(\cP_2) (a \natural b) = 
q^{1/2}  (f_{\beta}(\cP_2) a) \natural e^{2 \pi \alpha \cP_2} b
+ q^{-1/2} e^{-2 \pi \alpha \cP_2} a \natural f_{\beta}(\cP_2) b
$$
which implies that $f_{\beta}(\cP_2) (a \natural b) \in \gA_2(\bbbr^2)$. 
Similar verifications can be done for the other cases and components. Thus,
we have shown that $a {\scriptscriptstyle{\circ}} b \in \gA(\bbbr^2)_4$ for 
$a, b\in \gA(\bbbr^2)_4$.

Recall that by construction the product ${\scriptscriptstyle{\circ}}$ and the 
involution $\star$ have been transformed from 
the products ${\scriptstyle{\#}}$ and the involution $\ast$, respectively, under the 
bijective mapping $J$. Hence $\gA(\bbbr^2)_4$ is a $\ast$-algebra, 
because the products ${\scriptstyle{\#}}$ and the involution $\ast$ satisfy all axioms of a $\ast$-algebra. \hfill $\Box$

\mn
Next we define the product of elements of the coordinate algebra  
$\cO(\bbbr^2_q)$ and the algebra $\gA(\bbbr^2)_4$. Because of (\ref{mxmy}) and (\ref{I13}), we consider the self-adjoint operators $\sx$ and $\sy$ on the Hilbert space $\G$ given by the the operator matrices
$$
\sx = \theta_1(e^{2 \pi \alpha \Q_1} \otimes e^{- \pi \alpha \cP_2}),~~
\sy = \theta_1(e^{\pi \beta \cP_1} \otimes e^{2 \pi \beta \Q_2}).
$$
and define 
\begin{equation}\label{prodxa}
x{\scriptscriptstyle{\circ}} a \equiv \Phi (x) a= \sx a,~~y {\scriptscriptstyle{\circ}} a \equiv \Phi (y)a = \sy a
\end{equation}
for $a \in \gA(\bbbr^2)_4$. Using once more assumption (\ref{assump}) it follows that each $a \in \gA(\bbbr^2)_4$ is contained in the domains of the
self-adjoint operators $\sx$ and $\sy$, so that (\ref{prodxa}) is 
well-defined. 

An arbitrary element $a \in \gA(\bbbr^2)_4$ is in general not in the domains
of the powers ${\sx}^n$ and ${\sy}^n$ for $n \in \bbbn$ and the expressions 
$x {\scriptscriptstyle{\circ}} a$ and $y {\scriptscriptstyle{\circ}} a$ 
defined by (\ref{prodxa}) are in general not in $\gA(\bbbr^2)_4$. Hence the direct sum $\A(\bbbr^2_q) = \cO(\bbbr^2_q)+\gA(\bbbr^2)_4$ of 
$\ast$-algebras $\cO(\bbbr^2_q)$ and $\gA(\bbbr^2)_4$ is not 
an algebra with respect to the product (\ref{prodxa}). 
However, for certain elements $z \in\cO(\bbbr^2_q)$ and 
$a \in \gA(\bbbr^2)_4$ the above definition (\ref{prodxa}) gives 
indeed a well-defined product $z {\scriptscriptstyle{\circ}} a$ and $\A(\bbbr^2_q)$ 
becomes a {\it partial} $\ast$-algebra with product (\ref{prodxa}) in 
this manner. We shall not pursue this further. 

We now replace $\gA(\bbbr^2)_4$ by its $\ast$-subalgebra
$$
\gA_0(\bbbr^2)_4:= \gA(\bbbr^2)\oplus \gA(\bbbr^2)\oplus \gA(\bbbr^2) 
\oplus \gA(\bbbr^2).
$$
According to our general picture, the elements of this subalgebra can be considered as functions on the real quantum plane which are vanishing on the coordinate axis. For $a \in \gA_0(\bbbr^2_q)$, the elements $\Phi(x)a$ and $\Phi(y)a$ are obviously again in $\gA_0(\bbbr^2)_4$, so equation 
(\ref{prodxa}) defines a $\ast$-representation $\Phi$ of the $\ast$-algebra $\cO(\bbbr^2_q)$ on the invariant dense domain $\gA_0(\bbbr^2)_4$ of the Hilbert space $\Hh$. By Lemma 4(ii), the operators $\Phi(x^n)$ 
and $\Phi(y^n)$ are essentially self-adjoint on $\gA_0(\bbbr^2)_4$. 
For $z \in \cO(\bbbr^2_q)$ and $a \in \gA_0(\bbbr^2)_4$ we define
\begin{equation*} 
z {\scriptscriptstyle{\circ}} a := \Phi (z)a.
\end{equation*}
Then the direct sum $\A_0(\bbbr^2_q):= \cO(\bbbr^2_q)+\gA_0(\bbbr^2)_4$ 
becomes a $\ast$-algebra with product (\ref{prodxa}) such that 
$\cO(\bbbr^2_q)$ and $\gA_0(\bbbr^2)_4$ are $\ast$-subalgebras. Indeed, 
$\A_0(\bbbr^2_q)$ is merely the transformation 
of the left $\U_q(gl_2(\bbbr))$-module $\ast$-algebra $\gB(\bbbr^2)_4$ 
defined in 6.1 (see (\ref{algebrab})) under the symmetry $J$. Therefore, 
$\A_0(\bbbr^2_q)$ is a $\ast$-algebra with product (\ref{prodxa}) 
and there is a left action of $\U_q(gl_2)$ given by 
formulas (\ref{leftac}) such that $\A_0(\bbbr^2_q)$ is a left 
$\U_q(gl_2(\bbbr))$-module $\ast$-algebra.

For the study of the coordinate functions $x$ and $y$ the $\ast$-algebra 
$\gA_0(\bbbr^2)_4$ is ``large enough'', since the operators $\Phi(x)$ and 
$\Phi(y)$ are essentially self-adjoint on $\gA_0(\bbbr^2)_4$. However, 
for the action of the generators of $\U_q(gl_2(\bbbr))$ it is not, 
because $\gA_0(\bbbr^2)_4$ is not a core 
for the essentially self-adjoint operators $\Phi(E^\prime)$ and 
$\Phi(F^\prime)$. The larger $\ast$-algebra $\gA(\bbbr^2)_4$ 
is needed in this case. Since we have only an action of the generators and 
not of the whole algebra $\U_q(gl_2))$, $\gA(\bbbr^2)$ cannot be a 
left $\U_q(gl_2(\bbbr))$-module algebra. But {\it for the generators 
$f = E^\prime, F^\prime, K_1, K_2$ the two conditions (\ref{A8}) and 
(\ref{A15}) are valid on the $\ast$-algebra $\gA(\bbbr^2)_4$.} 

The proof of the latter assertion requires a number of computations. 
We carry out this verification only for the generator $f = E^\prime$ and 
for elements of the form $a=(a_1,0,0,0), b=(b_1,0,0,0) \in \gA(\bbbr^2)_4$, where $a_1, b_1 \in f_\alpha(\cP_1)^{-1}\gA(\bbbr^2)$. First we note that from formulas (\ref{F5})--(\ref{F8}) and (\ref{atimesb1}) we derive that 
\begin{align}\label{eq1times}
e^{-2 \pi \alpha \Q_1} (c \natural d) &= 
q^{-1/4} (e^{-2 \pi \alpha \Q_1}c) \natural (e^{- \pi \alpha \cP_2}d) =
q^{1/4} (e^{ \pi \alpha \cP_2}c)\natural (e^{-2 \pi \alpha \Q_1}d),\\
\label{eq2times}
e^{2 \pi \beta \Q_2} (c \natural d) &= 
q^{-1/4} (e^{2 \pi \beta \Q_2}c) \natural (e^{- \pi \beta \cP_1}d) =
q^{1/4} (e^{ \pi \beta \cP_1}c)\natural (e^{2 \pi \beta \Q_2}d)
\end{align}
for arbitrary elements $c, d \in f_\alpha(\cP_1)^{-1}\gA(\bbbr^2)$. 
Using the definitions of the product ${\scriptscriptstyle{\circ}}$ and of 
the operators
$\Phi(E^\prime)$ and $\Phi(K)$ and formulas (\ref{fp1}), (\ref{eq1times}) and 
(\ref{eq2times}) we compute
\begin{align*}
(\Phi(E^\prime)( a \natural b))_1 &= 
R_\alpha\otimes e^{2\pi\beta\Q_2} (a_1 \natural b_1)
\\
&=e^{-2 \pi \alpha \Q_1} \otimes e^{2 \pi \beta \Q_2}
f_\alpha(\cP_1) (a_1 \natural b_1)
\\
& = q^{1/2} e^{-2 \pi \alpha \Q_1} \otimes e^{2 \pi \beta \Q_2}
 ((f_\alpha(\cP_1) a_1) \natural e^{2 \pi \beta \cP_1} b_1)\\
&\quad + q^{-1/2}e^{-2 \pi \alpha \Q_1} \otimes e^{2 \pi \beta \Q_2} 
(e^{-2 \pi \beta \cP_1} a_1 \natural f_\alpha(\cP_1) b_1)
\\
& =(e^{-2 \pi \alpha \Q_1} \otimes e^{2 \pi \beta \Q_2}
 f_\alpha(\cP_1) a_1) \natural 
(e^{ \pi \beta \cP_1} \otimes e^{- \pi \alpha \cP_2} b_1)\\
&\quad +(e^{- \pi \beta \cP_1} \otimes e^{ \pi \alpha \cP_2} a_1) \natural
(e^{-2 \pi \alpha \Q_1} \otimes e^{2 \pi \beta \Q_2} 
 f_\alpha(\cP_1) b_1)
\\ 
&=(\Phi(E^\prime)a)_1 \natural (\Phi(K)b)_1 
+(\Phi(K^{-1})a)_1 \natural (\Phi(E^\prime)b)_1,
\end{align*}
where the lower index $1$ always refers to the first component. This proves 
condition (\ref{A8}) for $f = E^\prime$ and the particular elements $a$ 
and $b$. The other cases can be treated in a similar manner. 

Next we want to have a counterpart on $\gA(\bbbr^2)_4$ of the covariant 
linear functional $h\equiv h_0$ on $\gA(\bbbr^2)$. Recall that this 
counterpart on the $\ast$-algebra $\cB(\bbbr^2)_4$ is the functional 
$\tilde{h}$ defined by (\ref{htildeb}). We have to transform this 
functional $\tilde{h}$ under the symmetry $J$ by setting 
$h(a):=\tilde{h}(Ja)$. By the definitions (\ref{E6}) of $J$ 
and (\ref{htildeb}) of $\tilde{h}$ we obtain $h(a)=2\tilde{h}(a_1)$, where $\tilde{h}(a_1)$ is given by (\ref{htilde}). Inserting formula (\ref{htilde}) we are lead to define
\begin{equation}\label{ha4}
h(a)=2\iint e^{\pi(\alpha x_1+\beta x_2)} a_1(x_1,x_2) dx_1 dx_2,~ a=(a_1,a_2,a_3, a_4)\in\gA(\bbbr^2)_4.
\end{equation}
Let us check first that $h(a)$ is well-defined for arbitrary elements 
$a\in\gA(\bbbr^2)_4$. Indeed, if $a\in\gA(\bbbr^2)_4$, then we have 
$a_1\in f_\alpha (\cP_1)^{-1}\gA(\bbbr^2)+f_{\beta} 
(\cP_2)^{-1}\gA(\bbbr^2)$. From assumption (\ref{assump}) 
we conclude that $a_1$ is in the domain 
$\cD(e^{\pi\alpha\Q_1}\otimes e^{\pi\beta\Q_2})$ and that 
$(e^{\pi\alpha\Q_1}\otimes e^{\pi\beta\Q_2})a_1\in\cD_{\nu,\mu}$ 
for some $\nu,\mu\in\bbbr^{++}$. By Lemma 9(ii), the 
latter implies that the function $e^{\pi(\alpha x_1+\beta x_2)} a_1(x_1,x_2)$ is in the Schwartz space $\sS(\bbbr^2)$. Hence the integral in (\ref{ha4}) exists.

From the construction it follows that the sesquilinear form 
$\langle\cdot,\cdot\rangle_h$ associated with $h$ by (\ref{G4}) is the 
scalar product of the Hilbert space $\Hh$ (see (\ref{hilberth})). This can 
be also verified directly by using formulas (\ref{twpa}), (\ref{twia}), 
(\ref{atimesb}), (\ref{star}) and (\ref{F5})-(\ref{F12}). That is, we have
\begin{equation}\label{hl2}
\langle a,b\rangle_h=h(b^\star{\scriptscriptstyle{\circ}} a)={\iint_{\bbbr^2}} (a_1\overline{b_1}+a_2\overline{b_2}+a_3\overline{b_3}+a_4\overline{b_4}) dx_1dx_2
\end{equation}
for $a=(a_1,a_2,a_3,a_4),b=(b_1,b_2,b_3,b_4)\in\gA(\bbbr^2)_4$. Thus we 
have reduced the scalar product $\langle\cdot,\cdot\rangle_h$ to 
$L^2$-scalar products on $\bbbr^2$. To achieve formula (\ref{hl2}) was 
the main aim of our considerations and it is the reason why we have 
transformed all structures by means of the operator $T$ and the symmetry 
$J$.

Finally, let us turn to the differential calculus on the quantum plane. The 
construction from 3.2 carries over almost verbatim to the algebra 
$\A_0(\bbbr^2_q):=\cO(\bbbr^2_q)+\gA_0(\bbbr^2)_4$ and yields a first 
order differential calculus $\Gamma$ over the algebra $\A_0(\bbbr^2_q)$ 
such that $\{dx, dy\}$ is a free left module basis of $\Gamma$. The 
corresponding partial derivatives $\partial_x$ and $\partial_y$ 
are of the form
\begin{align*}
&\partial_x(a)=\left(\begin{matrix} 0&\partial_x(a_1) &0& 0\\
                                    \partial_x(a_2)&0&0&0\\
                                    0&0&0&\partial_x(a_3)\\
                                    0&0&\partial_x(a_4)&0\end{matrix}\right).\\
\\
&\partial_y(a)=\left(\begin{matrix} 0&0&\partial_y(a_1) &0\\
                                    0&0&0&\partial_y(a_2)\\
                                    \partial_y(a_3)&0&0&0\\
                                    0&\partial_y(a_4)&0&0\end{matrix}\right)
\end{align*}
for $a=(a_1,a_2,a_3,a_4)\in\gA_0(\bbbr^2)_4$, 
where $\partial_x(a_j)$ and $\partial_y(a_j)$, $j=1,2,3,4$, are as in 3.2.

Recall that $\A(\bbbr^2_q)=\cO(\bbbr^2_q)+\gA(\bbbr^2)_4$ is not an 
algebra, because not all $a\in\gA(\bbbr^2)_4$ are in the domains of 
$\Phi(x)^n\Phi(y)^m$. Likewise, the differential $da=\omega a-a\omega$ and the partial derivatives $\partial_x(a)$ and $\partial_y(a)$ are well-defined only for those elements $a\in\gA(\bbbr^2)_4$ belonging to the corresponding operator domains.

We close this subsection by listing the formulas of the operator 
matrices for the generators $E^\prime, F^\prime, K_1, K_2, x$ and $y$. 
Recall that 
\begin{align*}
&L_\alpha=\overline{f_\alpha} (\cP) 
e^{-2\pi\alpha Q},~~ R_\alpha=
e^{-2\pi\alpha\Q}f_\alpha(\cP),\\
&f_\alpha(\cP) = 
- 2\sinh\pi\beta(2\cP{+}\alpha i) =
-q^{1/2} e^{2 \pi \beta \cP} + q^{-1/2} e^{-2 \pi \beta \cP}.
\end{align*} 
\begin{align*}
&\overline{\Phi(E^\prime)} = \sE =\left( \begin{matrix} 
                     0 &0 &0 &L_\alpha\otimes e^{2\pi\beta\Q_2}\\
                     0 &0 &R_\alpha\otimes e^{2\pi\beta\Q_2} &0\\
                     0 &L_\alpha\otimes e^{2\pi\beta\Q_2} &0 &0\\
                     R_\alpha\otimes e^{2\pi\beta\Q_2} &0 &0 &0
                 \end{matrix} \right)~,\\
\\
&\overline{\Phi (F^\prime)} = \sF =\left( \begin{matrix}
                    0 &0 &0 &e^{2\pi\alpha\Q_1}\otimes L_{\beta}\\
                    0 &0 &e^{2\pi\alpha\Q_1}\otimes L_{\beta} &0\\
                    0 &e^{2\pi\alpha\Q_1}\otimes R_{\beta} &0 &0\\
                    e^{2\pi\alpha\Q_1}\otimes R_{\beta} &0 &0 &0
                  \end{matrix} \right)~,\\
\\
&\overline{\Phi (q^{-1/4}K_1)} = \sK_1 =\left( \begin{matrix}
                  e^{\pi\beta\cP_1}\otimes I &0 &0 &0\\
                  0 &e^{\pi\beta\cP_1}\otimes I &0 &0\\
                  0 &0 &e^{\pi\beta\cP_1}\otimes I &0\\
                  0 &0 &0 &e^{\pi\beta\cP_1}\otimes I
                \end{matrix} \right)~,
\end{align*}
\begin{align*}
&\overline{\Phi (q^{-1/4}K_2)} = \sK_2 =\left( \begin{matrix}
                  I \otimes e^{\pi\alpha\cP_2} &0 &0 &0\\
                  0 &I \otimes e^{\pi\alpha\cP_2} &0 &0\\
                  0 &0 &I \otimes e^{\pi\alpha\cP_2} &0\\
                  0 &0 &0 &I \otimes e^{\pi\alpha\cP_2}
                \end{matrix} \right)~,\\
\\
&\overline{\Phi (x)} = \sx = \left( \begin{matrix}
                0 &e^{2\pi\alpha\Q_1}\otimes e^{{-}\pi\alpha\cP_2} &0 &0\\
                e^{2\pi\alpha\Q_1}\otimes e^{{-}\pi\alpha\cP_2} &0 &0 &0\\
                0 &0 &0 &e^{2\pi\alpha\Q_1}\otimes e^{{-}\pi\alpha\cP_2}\\
                0 &0 &e^{2\pi\alpha\Q_1}\otimes e^{{-}\pi\alpha\cP_2} &0\\
               \end{matrix} \right)~,\\
\\
&\overline{\Phi (y)} =\sy = \left( \begin{matrix}
                0 &0 &e^{\pi\beta\cP_1}\otimes e^{2\pi\beta\Q_2} &0\\
                0 &0 &0 &e^{\pi\beta\cP_1}\otimes e^{2\pi\beta\Q_2}\\
                e^{\pi\beta\cP_1}\otimes e^{2\pi\beta\Q_2} &0 &0 &0\\
                0 &e^{\pi\beta\cP_1}\otimes e^{2\pi\beta\Q_2} &0 &0
              \end{matrix} \right)~.
\end{align*}

\mn
{\it 6.4} In this last subsection we introduce three unitary operators 
$\F^q_x,\F^q_y$ and $\F^q$ on the 
Hilbert space $\Hh$ which can be considered as quantum 
analogs of the partial Fourier transforms and the Fourier transform 
on $\bbbr^2$, respectively.

First let us note that the counterparts of the $q$-deformed partial 
derivatives $\cD^q_x$ and 
$\cD^q_y$ (see (\ref{dxy})) on the algebra $\A(\bbbr^2_q)$ are defined by
\begin{align*}
&\cD^q_x{:=}\sK\sy^{-1}\sE{=}\left( \begin{matrix} 0&L_{\alpha}\otimes e^{-\pi\alpha\cP_2} &0&0\\
R_{\alpha}\otimes e^{-\pi\alpha\cP_2} &0&0&0\\
0&0&0&L_\alpha{\otimes} e^{-\pi\alpha\cP_2}\\
0&0&R_\alpha{\otimes} e^{-\pi\alpha\cP_2}&0\end{matrix}\right)~,\\
&\\
&\cD^q_y {:=}\sK \sx^{-1}\sF{=}\left( \begin{matrix} 0&0 &e^{\pi\beta\cP_1}{\otimes} L_\beta &0\\
0&0&0&e^{\pi\beta\cP_1}{\otimes} L_\beta\\
e^{\pi\beta\cP_1}{\otimes} R_\beta &0&0&0\\
0&e^{\pi\beta\cP_1}{\otimes} R_\beta &0&0 \end{matrix} \right)
\end{align*}
Clearly, $\cD^q_x$ and $\cD^q_y$ are self-adjoint operators on the 
Hilbert space $\Hh$.

Let $u$ be the unitary operator on $L^2(\bbbr)$ 
given by $(uf)(x)=f(-x)$ and let $w_\alpha,v_\alpha$ resp. 
$w_\beta,v_\beta$ be the holomorphic functions from Lemma 6. Then
\begin{align*}
\F^q_x:=\sigma_1(u\overline{w}_\alpha(\cP_1)\otimes I, 
u\overline{v}_\alpha(\cP_1)\otimes I),~ 
\F^q_y:=\sigma_2(I\otimes u\overline{w}_\beta(\cP_2),
I\otimes u\overline{w}_\beta(\cP_2)).
\end{align*}
are commuting unitaries on the Hilbert space $\Hh$. Set 
$\F^q{:=}\F^q_x\F^q_y$, that is,
\begin{equation*}
\F^q{=}\left( \begin{matrix} u\overline{w}_\alpha(\cP_1){\otimes} u\overline{w}_\beta(\cP_2) &0&0&0\\
0&\makebox[1cm][c]{$u\overline{v}_\alpha(\cP_1){\otimes} u\overline{w}_\beta(\cP_2)$} &0&0\\
0&0&u\overline{w}_\alpha(\cP_1){\otimes} u\overline{v}_\beta(\cP_2) &0\\
0&0&0 &u\overline{v}_\alpha(\cP_1){\otimes} 
u\overline{v}_\beta(\cP_2)\end{matrix}\right).
\end{equation*}

We call the unitaries $\F^q_x$ and $\F^q_y$ 
{\it quantum partial Fourier transforms} and  
$\F^q$ {\it quantum Fourier transform} of the 
real quantum plane. The reason for this terminology stems from the fact 
that, roughly speaking, these unitaries interchange the coordinate 
functions $x,y$ and the $q$-deformed partial derivatives 
$\cD^q_x,\cD^q_y$, respectively. 
More precisely, we have the following relations.

\mn
{\bf Proposition 27.} (i) $\F^q_x \sx(\F^q_x)^{-1} = -\cD^q_x, 
\F^q_x \cD^q_x (\F^q_x)^{-1}{=} \sx, 
\F^q_x \sy (\F^q_x)^{-1} {=} \sK^{-2}_1 \sy,\\ 
\F^q_x \cD^q_y (\F^q_x)^{-1}{=}\sK^{-2}_1 \cD^q_y.$\\
(ii) $\F^q_y \sy (\F^q_y)^{-1}{=}-\cD^q_y, \F^q_y \cD^q_y (\F^q_x)^{-1} {=} 
\sy, \F^q_y \sx (\F^q_y)^{-1} \sK^2_2 \sx, \
F^q_y \cD^q_x (\F^q_y)^{-1} {=} \sK^2_2 \cD^q_x$.\\
(iii) $ \F^q \sx (\F^q)^{-1} {=} -\sK^2_2 \cD^q_x, \F^q \sy (\F^q)^{-1} {=} -\sK^{-2}_1 \cD^q_y, \F^q \cD^q_x (\F^q)^{-1} {=} \sK^2_2 \sx,\\
\F^q \cD^q_y (\F^q)^{-1}{=} \sK^{-2}_2 \sy.$

\mn
{\bf Proof.} (i): By Lemma 7, written in terms of matrix entries, we have 
\begin{align}\label{uwva}
&u\overline{w}_\alpha (\cP)L_\alpha v_\alpha(\cP)u=e^{2\pi\alpha\Q}, u\overline{v}_\alpha(\cP) R_\alpha w_\alpha(\cP)u=e^{2\pi\alpha\Q},\\
\label{uwvb}
&u\overline{w}_\alpha(\cP)e^{2\pi\alpha\Q} v_\alpha(\cP)u{=}uL_{-\alpha} u{=}-L_\alpha,
u\overline{v}_\alpha(\cP)e^{2\pi\alpha\Q} w_\alpha(\cP)u{=}uR_{-\alpha} u{=}-R_\alpha.
\end{align}
Because of the two relations (\ref{uwva}) the matrix entries of 
$\F^q_x\cD^q_x(\F^q_x)^{-1}$ and $\sx$ coincide, while (\ref{uwvb}) 
implies that $\F^q_x\sx(F^q_x)^{-1}$ and $-\cD^q_x$ have the same matrix 
entries. This proves the first two relations of (i). The two other 
relations follow at once from the corresponding definitions combined with 
the fact that $ue^{\pi\beta\cP} u=e^{-\pi\beta\cP}$. 

The proof of (ii) is similar to the proof of (i). 
The relations of (iii) follow easily from (i) and (ii).\hfill$\Box$

\mn
Remark 5. The holomorphic functions $w_{-\alpha}$ and $v_{-\alpha}$ coincide with the functions $w_1$ and $w_2$ in [S4], where a closely related quantum Fourier transform of a $q$-deformed Heisenberg algebra appeared. Holomorphic functions of similar kind 
have been used by S.L. Woronowicz [W] in another context 
as quantum exponential functions for the quantum $ax+b$--group. 

\mn


\end{document}